\documentclass[11pt,reqno]{amsart}  
\usepackage{amsmath,amssymb,amsthm}
\theoremstyle{plain}
\newtheorem{theorem}{Theorem}[section]
\newtheorem{lemma}[theorem]{Lemma}
\newtheorem{proposition}[theorem]{Proposition}

\newtheorem{corollary}[theorem]{Corollary}

\numberwithin{equation}{section}
\allowdisplaybreaks[1]

\theoremstyle{definition}

\newtheorem{definition}[theorem]{Definition}
\newtheorem{remark}[theorem]{Remark}
\newtheorem{example}[theorem]{Example}

\theoremstyle{remark}

\newcommand{\bbC}{{\mathbb C}}

\newcommand{\fQ}{{\mathfrak Q}}
\newcommand{\fS}{{\mathfrak S}}

\newcommand{\cA}{{\mathcal A}}
\newcommand{\cB}{{\mathcal B}}
\newcommand{\cC}{{\mathcal C}}
\newcommand{\cD}{{\mathcal D}}
\newcommand{\cE}{{\mathcal E}}

\newcommand{\cG}{{\mathcal G}}
\newcommand{\cH}{{\mathcal H}}
\newcommand{\cI}{{\mathcal I}}

\newcommand{\cL}{{\mathcal L}}
\newcommand{\cM}{{\mathcal M}}
\newcommand{\cN}{{\mathcal N}}

\newcommand{\cR}{{\mathcal R}}
\newcommand{\cS}{{\mathcal S}}

\newcommand{\cU}{{\mathcal U}}
\newcommand{\cV}{{\mathcal V}}

\newcommand{\cX}{{\mathcal X}}
\newcommand{\cY}{{\mathcal Y}}

\newcommand{\fA}{{\mathfrak A}}
\newcommand{\fX}{{\mathfrak X}}
\newcommand{\bGamma}{{\boldsymbol{\Gamma}}}

\newcommand{\norm}[1]{||#1||}

\newcommand{\bbZ}{{\mathbb Z}}

\newcommand{\sbm}[1]{\left[\begin{smallmatrix} #1
		\end{smallmatrix}\right]}

\newcommand{\ip}[1]{\langle #1 \rangle}

\begin{document}

\title[Free noncommutative hereditary kernels]
{Free noncommutative hereditary kernels:  Jordan decomposition, Arveson extension,  kernel domination}
\author[J.A.\ Ball]{Joseph A. Ball}
\address{Department of Mathematics,
Virginia Tech,
Blacksburg, VA 24061-0123, USA}
\email{joball@math.vt.edu}
\author[G.\ Marx]{Gregory Marx}
\address{Department of Mathematics,
Virginia Tech,
Blacksburg, VA 24061-0123, USA}
\email{marxg@vt.edu}
\author[V.\ Vinnikov]{Victor Vinnikov}
\address{Department of Mathematics, Ben-Gurion University of the 
Negev, Beer-Sheva, Israel, 84105}
\email{vinnikov@cs.bgu.ac.il}

\begin{abstract} 
We discuss a (i) quantized version of the Jordan decomposition theorem  for a complex Borel measure on a compact Hausdorff space,
namely, the more general problem of decomposing a general noncommutative kernel (a quantization of the standard notion of 
kernel function) as a linear combination of completely positive noncommutative kernels (a quantization of the standard notion of 
positive definite kernel).  Other special cases of (i)  include: the problem of decomposing a general operator-valued kernel function 
as a linear combination
of positive kernels (not always possible), of decomposing a general  bounded linear Hilbert-space operator
as a linear combination of positive linear operators (always possible), of decomposing a completely bounded
linear map from a $C^*$-algebra $\cA$ to an injective $C^*$-algebra $\cL(\cY)$  as a linear combination of completely positive maps
from $\cA$ to $\cL(\cY)$ (always possible).  We also discuss (ii) a noncommutative kernel generalization of the Arveson extension 
theorem (any completely positive map $\phi$ from a operator system ${\mathbb S}$ to an injective $C^*$-algebra $\cL(\cY)$ 
can be extended to a completely positive map $\phi_e$ from a $C^*$-algebra containing ${\mathbb S}$ to $\cL(\cY)$), and (iii)
a noncommutative kernel version of a Positivstellensatz (i.e., finding a certificate
to explain why one kernel is positive at points where another given kernel is positive).

\end{abstract}


\maketitle

\tableofcontents

\section{Introduction}  \label{S:Intro}
\setcounter{equation}{0}

The idea of {\em quantized functional analysis} came out of the attempt to understand intrinsically  spaces of operators
beyond the category of concrete $C^*$-algebras (a closed subalgebra of $\cL(\cH)$ closed under taking of adjoints), e.g.
subalgebras of $\cL(\cH)$ (operator algebras),  a linear subspace of $\cL(\cH)$ (operator space), a unital subspace
closed under taking adjoints (operator system).  To get an intrinsic characterization of such objects, unlike the $C^*$-algebra
case where the axioms for a $C^*$-algebra does the job, it was found that one needs to study not only
the subspace as an abstract Banach space, but also a system of compatible norms on matrices over the subspace  and one 
must study such objects up to completely isometric isomorphism (preserving not only all the structure on the primordial space
$\cX$ but also  on all the matricial spaces $\cX^{n \times n}$ for all positive integers $n \in {\mathbb Z}_+$ (see
\cite{Paulsen, BlM, ER} for systematic treatments).  The same idea  has now invaded function theory with impressive results
(see  \cite{KVV-book}).  The purpose of this contribution is to continue this line of research with a study of the
quantized version of kernels begun in \cite{BMV1} and continued in \cite{BMV2}.  

By a classical kernel on a set $\Omega$ we mean a function $K$ on the Cartesian product set $\Omega \times \Omega$
with values in some space, usually a linear space, e.g., the space of operators on some Hilbert space
$\cL(\cY)$:  $(x, y) \mapsto K(x, y)$.  Let us say that the kernel is \textbf{Hermitian} if $k(x,y)^* = k(y,x)$.  
A \textbf{positive kernel}  (in the sense of Aronszajn \cite{Aron}) is one for which
$$
  \sum_{i,j=1}^N \langle K(z_i, z_j) y_j , y_i \rangle \ge 0 
$$
for all $z_1, \dots, z_N \in \Omega$ and $y_1, \dots, y_n \in \cY$.  An equivalent characterization of positive kernels  
is that they all have a \textbf{Kolmogorov decomposition}, i.e., one can find
operator-valued functions $H \colon \Omega \to \cL(\cX, \cY)$ (where $\cX$ is some Hilbert state space) so that
$K$ has the factorization $K(z,w) = H(z) H(w)^*$.  

A profound generalization of positive kernel is that of \textbf{completely positive kernel} given by 
Barreto-Bhat-Liebscher-Skeide
\cite{BBLS} (simplified here to the Hilbert space/Hilbert-space-operator context rather than 
Hilbert-module/$C^*$-correspondence
setting in \cite{BBLS}) whereby $K$ is still a function on $\Omega \times \Omega$ but takes values in the space
$\cL(\cA, \cL(\cY))$ of linear operators from a $C^*$-algebra $\cA$ to the space of Hilbert-space operators 
$\cL(\cY)$ and 
is required to satisfy the more elaborate positivity condition
\begin{equation}  \label{BBLS-pos}
 \sum_{i,j = 1}^N \langle  K(z_i, z_j)(a_i^* a_j) y_j, y_i \rangle \ge 0
\end{equation}
for all choices of $z_1, \dots, z_N \in \Omega$, $a_1, \dots, a_N \in \cA$, and $y_1, \dots, y_N \in \cY$.  The main result
in \cite{BBLS} concerning such kernels (\textbf{cp BBLS kernels} for short) is that they are characterized by having the
following more elaborate  \textbf{Kolmogorov decomposition}:
{\em there exist a Hilbert state space $\cX$, an operator-valued function $H \colon \Omega \to \cL(\cX, \cY)$ and a 
$*$-representation $\pi \colon \cA \to \cL(\cX)$ so that}
\begin{equation}   \label{BBLS-Kol}
    K(z,w)(a) = H(z) \pi(a) H(w)^*.
\end{equation}
More generally one can consider BBLS kernels without the positivity condition \eqref{BBLS-pos}:  we say that
any function $K$ from $\Omega \times \Omega \to \cL(\cA, \cL(\cY))$ is a \textbf{BBLS-kernel}.  If it is the case that
$K(z,w)(a)^* = K(w,z)(a^*)$, we will say that $K$ is a \textbf{Hermitian BBLS-kernel}.

Before discussing noncommutative kernels, we discuss noncommutative functions.
 A classical function $f$ is defined on some point set $\Omega$
with values in some target space which we take to be a linear space $\cV_0$.  If $\Omega$ has some 
additional structure
(topological and analytic), then we speak about $f$ being continuous or holomorphic.  The idea of a 
free noncommutative
function (\textbf{nc function} for short) $f$ is the quantization of a classical function: the domain $\Omega$ is 
partitioned up into levels
 $\Omega = \amalg_{n=1}^\infty \Omega_n$
where $\Omega_n$ consists of $n \times n$ matrices over the ambient vector space $\cV$ and similarly for the
target space:  $\cV_0 = \amalg_{n=1}^\infty (\cV_0)^{n \times n}$. The nc function $f$
is required to be \textbf{graded}
\begin{equation}  \label{ncfunc-graded}
Z \in \Omega_n \Rightarrow f(Z) \in (\cV_0)^{n \times n}
\end{equation}
and to \textbf{respect intertwining conditions:}
\begin{equation}  \label{ncfunc-intertwine}
Z \in \Omega_n, \, \widetilde Z \in \Omega_{\widetilde n}, \, \alpha \in {\mathbb C}^{n \times \widetilde n}, \,
 \alpha \widetilde Z = Z \alpha \Rightarrow 
\alpha f(\widetilde Z) = f(Z) \alpha.
\end{equation}
It is shown in \cite{KVV-book} how this algebraic condition along with
some weak topological conditions implies holomorphic structure for a nc function.

A  noncommutative kernel (\textbf{nc kernel} for short) $K$ is the quantization of a BBLS kernel as we now explain.
 The domain $\Omega$ for a nc kernel
$K$ again is partitioned up into levels $\Omega = \amalg_{n=1}^\infty \Omega_n$ where $\Omega_n$
consists of those elements of $\Omega$ which are in $\cV^{n \times n}$ (where $\cV$ is the ambient vector 
space for the domain).  
Then the target domain for $K$ is partitioned up as
$$
   \amalg_{n,m = 1}^\infty \cL(\cA^{n \times m}, \cL(\cY)^{n \times m}) =: \cL(\cA_{\rm nc}, \cL(\cY)_{\rm nc}) \text{ (for short)}
 $$
 and $K$ is \textbf{graded} in the following sense:  
 \begin{equation}   \label{ker-graded}
 Z \in \Omega_n, W \in \Omega_m \Rightarrow
 K(Z,W) \in \cL(\cA^{n \times m}, \cL(\cY)^{n \times m}).
 \end{equation}
   A completely positive noncommutative kernel (\textbf{cp nc kernel} for short) $K \colon
 \Omega \times \Omega \to \cL(\cA_{\rm nc}, \cL(\cY)_{\rm nc})$ is characterized by having a {\em quantized
 BBLS-Kolmogorov representation} as follows:  {\em there exists a Hilbert space $\cY$, a
 nc function $H \colon \Omega \to \cL(\cX, \cY)_{\rm nc}$, a unital $*$-representation $\pi \colon \cA \to \cL(\cX)$
 so that, for $Z \in \Omega_n$ and $W \in \Omega_m$ and $P = [ P_{ij} ] \in \cA^{n \times m}$,}
 \begin{equation}   \label{nc-Kol}
   K(Z,W)([ P_{ij}]) = H(Z) [ \pi(P_{ij}) ] H(W)^*.
 \end{equation}
Given such a representation \eqref{nc-Kol},  one can deduce \textbf{respects intertwining conditions} for $K$ 
from both the left and the right
sides (see Section \ref{S:nc-ker} below).  These conditions make sense without the Kolmogorov representation \eqref{nc-Kol}
holding but just the gradedness condition \eqref{ker-graded} holding; these conditions then become the definition of a
{\em noncommutative kernel}  (\textbf{nc kernel} for short) (not necessarily completely positive). A nc kernel  $K$ is said to be a \textbf{Hermitian nc kernel} if in addition
$K(Z,W)(P)^* = K(W,Z)(P^*)$.  A cp nc kernel can also be defined via a quantization of condition \eqref{BBLS-pos} being
imposed on a nc kernel: {\em a nc kernel
$K$ is cp if and only if }
\begin{equation}  \label{nc-pos-ker}
\sum_{i,j=1}^N \langle K(Z_i, Z_j)(P_i P_j^*) y_j, y_i \rangle \ge 0
\end{equation}
for all $Z_i \in \Omega_{n_i}$, $P_i \in \cA^{n_i \times 1}$, $y_i \in \cY^{n_i}$, $n_i \in {\mathbb N}$ for $i=1, \dots, N$, 
$N \in {\mathbb N}$.

Let us note various special cases of cp nc kernels:  a cp nc kernel with $\Omega = \Omega_1$ amounts to a  cp BBLS-kernel;
if in addition $\cA = {\mathbb C}$ and one identifies $K(Z,W)$ with $K(Z,W)(1)$, a cp BBLS-kernel becomes a Aronszajn positive
kernel.  For more complete details, see \cite{BMV1}.  In another direction, a completely positive map $\varphi$ between $C^*$-algebras
$\cA$ and $\cB$ can be identified with a cp nc kernel $K \colon \Omega \times \Omega \to \cL( \cA_{\rm nc}, \cB_{\rm nc})$,
with $\Omega$ having the special form that $\Omega_n$ consists of a single $n \times n$ diagonal matrix $\sbm{ z_0 & & \\ & \ddots & \\ & & z_0}$,
where $z_0$ is the unique point at level-1  of $\Omega$:  $\Omega_1 = \{ z_0 \}$. 
 This connection generalizes to the case where the point set $\Omega$ consists of finitely many points, say $N$,
 located at various levels;  for this case a cp nc kernel corresponds to a cp map from $\cA^M$ to $\cL(\cY)^M$
 which is also a $(\cS,\cS^*)$-bimodule map, where $\cS$ is a certain subalgebra of ${\mathbb C}^M$ determined by 
 the point set $\Omega$ (see Section \ref{S:encoding} below).  In fact this connection is one of the main tools for our work to follow.

Here we discuss three issues for general nc kernels.

\smallskip
\noindent
\textbf{Problem A.}  (Jordan decomposition for nc kernels.)  {\sl Given a free nc kernel
$$
K \colon \Omega \times \Omega \to \cL(\cA_{\rm nc}, \cL(\cY)_{\rm nc})
$$
on a  set of nc points $\Omega$, write $K$ as a linear combination of  
four cp nc  kernels $K_1$, $K_2$, $K_3$, $K_4$:
\begin{equation}   \label{4fold}
K(Z,W)(P) = K_1(Z,W)(P) - K_2(Z,W)(P) + i (K_3(Z,W)(P) - K_4(Z,W)(P))
\end{equation}
for all $Z \in \Omega_n$, $W \in \Omega_m$, $P \in \cA^{n \times m}$.}

\smallskip
We shall show that Problem A has a solution if (i) the set $\Omega$ is finite, (ii) each of the maps
$\{K(Z, Z) \colon \cA^{n_Z \times n_Z} \to \cL(\cY^{n_Z})$ is completely bounded, and (iii)
the  intertwining matrix subalgebra $\cS$ associated with $\Omega$ mentioned above is a $C^*$-algebra (i.e., 
$\Omega$  is {\em admissible} in the sense of Definition \ref{D:admissible} to come).
We also show that this result covers all the previously known particular cases where the conclusion is known to hold
while also excluding various counterexamples where the result is known to fail (see the discussion
around Examples \ref{E:1} and \ref{E:2} as well as Section \ref{S:KuroshDecom}).  In particular,
the special case $\cA = C(X)$ (continuous functions on a compact Hausdorf space $X$) and $\cY = {\mathbb C}$
gives us the classical Jordan decomposition for a complex Borel measure on $X$ (see Example \ref{E:1} below).

\smallskip

A consequence of the \textbf{respects intertwining conditions} property for a nc kernel is the 
\textbf{respects direct sums} property
(see Section \ref{S:nc-ker} below).  If $\Omega$ is closed under the taking of direct sums, then the positivity condition
\eqref{nc-pos-ker} can be reformulated more succinctly simply as
\begin{equation}   \label{cp-ker-sys}
Z \in \Omega_N, \, P \succeq 0 \text{ in } \cA^{N \times N} \Rightarrow K(Z,Z)(P) \succeq 0 \text{ in } \cL(\cY^N).
\end{equation}
Furthermore one can cut down on the number of points $Z$ one needs to consider if one insists that the map
$$  
P \in \cA^{N \times N} \mapsto K(Z,Z)(P) \in \cL(\cY^N)
$$
be cp from $\cA^{N \times N}$ to $\cL(\cY^N)$.
The advantage of this formulation is that it makes sense when the $C^*$-algebra $\cA$  is replaced by a unital selfadjoint
linear subspace ${\mathbb S}$ of some $C^*$-algebra $\cA$, or more abstractly, by an {\em operator system} 
${\mathfrak S}$ (see e.g.~\cite{Paulsen} for additional background).   This leads to the second problem to be discussed in this paper:

\smallskip
\noindent
\textbf{Problem B.} (Arveson extension theorem for completely positive noncommutative kernels.)
{\sl Given a cp nc kernel  ${\mathbb K}$ on a nc point set $\Omega$ with values mapping ${\mathbb S}_{\rm nc}$ into
$\cL(\cY)_{\rm nc}$ where ${\mathbb S}$ is an operator system  (a $*$-closed unital linear subspace 
of a $C^*$-algebra $\cA$),
find an extension of $K$ to the $C^*$-algebra $\cA$
$$
K \colon \Omega \times \Omega \to \cL(\cA, \cL(\cY))
$$
which is also a cp nc  kernel.}

\smallskip

Let us mention that for the case where the point set $\Omega$ is a single point at level 1 so the kernel amounts to a
cp map from ${\mathbb S}$ to $\cL(\cY)$, a result of Arveson (the Arveson extension theorem \cite{Arv1, Paulsen})
resolves Problem B in the affirmative.  We show how our technique of converting kernels to cp maps which are also bimodule maps
with respect to a certain matrix algebra $\cS$ associated with the point set $\Omega$ leads to a solution of Problem B
for the case where the point set is finite.  Combining this result with a general procedure of Kurosh (see \cite{Kurosh} as
well as \cite{AP}) leads to a solution of Problem B for the general case via a reduction of the infinite-point case to the finite-point
case.

\smallskip

Our next problem has the flavor of a Positivstellensatz for free  nc kernels, i.e., the problem of characterizing the form that a
free nc kernel must have if it is constrained to be positive at those points where another given kernel is positive.
We begin with a simplified version of the problem;  to actually solve the problem there are some additional 
hypotheses which must be incorporated. 

\smallskip

\noindent
\textbf{Problem C.} (Kernel dominance problem for  noncommutative kernels.)  {\sl Given HIlbert spaces $\cE$ and 
$\cG$, a full nc subset $\Xi$ of $\cV_{\rm nc}$, and  a  Hermitian nc kernel
$$
\fQ \colon \Xi \times \Xi \to \cL({\mathbb C}_{\rm nc}, \cL(\cS)_{\rm nc}),
$$
let ${\mathbb P}_\fQ$ be the strict positivity domain for $\fQ$ as defined by
$$
  {\mathbb P}_\fQ = \{ Z \in \Xi \colon \fQ(Z, Z)(I_n) \succ 0 \}.
$$
Suppose that $\Omega$ is a subset of ${\mathbb P}_\fQ$ 
and that $\fS$ is a Hermitian nc kernel defined on  $\Omega$ 
$$
\fS \colon \Omega \times \Omega \to \cL({\mathbb C}_{\rm nc}, \cL(\cY)_{\rm nc})
$$
which is positive semidefinite on $\Omega$:
$$
   \fS(Z, Z)(1_\cA) \succeq 0 \text{ for all } Z \in \Omega.
 $$
 Then we seek to find two completely positive nc kernels on $\Omega$
 $$
\Gamma_1 \colon \Omega \times \Omega  \to \cL(\cL(\cS)_{\rm nc}, {\mathbb C}_{\rm nc}),  \quad 
 \Gamma_2 \colon \Omega \times \Omega \to \cL({\mathbb C}_{\rm nc}, \cL(\cY)_{\rm nc})
$$
so that, for all $Z \in \Omega_n$, $W \in \Omega_m$, $P \in {\mathbb C}^{n \times m}$ we have the kernel decomposition}
\begin{equation}   \label{certificate}
\fS(Z,W)(P) = \Gamma_1(Z,W)(\fQ(Z,W)(P)) + \Gamma_2(Z, W)(P).
\end{equation}

\smallskip

Note that the representation \eqref{certificate} can be viewed as a certificate which explains why ${\mathfrak S}$
is dominated by the kernel ${\mathfrak Q}$ in the sense mentioned in the statement, i.e.,  it is immediate 
from the representation
\eqref{certificate} that $\fQ(Z,Z)(I) \succ 0$ on $\Omega$ then leads to $\fS(Q,Q)(I) \succeq 0$ on $\Omega$.

\smallskip

This formulation does not have an affirmative solution in general without imposing some additional hypotheses.  
The first additional hypothesis
is that the kernel $\fS$ should actually be defined on a somewhat larger set $\Omega'$ containing $\Omega$, namely
$\Omega' = [\Omega]_{\rm full} \cap {\mathbb P}_\fQ$ where $[\Omega]_{\rm full}$ is the {\em full envelope} of 
$\Omega$ (see Section \ref{S:tax} below for the precise definition).  In addition we need to assume  that each of t
he kernels $\fQ$ and $\fS$ 
has the property that its restriction to a finite subset of its domain ($\Xi$ or $\Omega'$ respectively) is decomposable.
In particular this is true if $\fS$ and $\fQ$ are assumed to be decomposable as kernels on their respective domains 
($\Xi$ or $\Omega'$ respectively)
 at the start.  The precise result is given as Theorem \ref{T:ker-dom} below.  

An interesting special case is the case where we take $\Omega = {\mathbb P}_\fQ$.  Then 
$\Omega' = {\mathbb P}_\fQ$
and  we assume that $\fS$ as well as $\fQ$ are free nc kernels defined on all of $\Xi$  and one asks that the decomposition
\eqref{certificate} holds on all of $\Xi$.  We are not able to get such a result in general;  however in case one assumes
that $\Xi = ({\mathbb C}^d)_{\rm nc}$ and  $\fQ$ and $\fS$ are nc polynomial kernels or (more generally)
nc rational kernels, then the result does hold if one imposes an additional Archimedean hypothesis (see \cite{HMcC2004, 
PascoePAMS}); a result of this type is called a {\em noncommutative Positivstellensatz}.  We refer to Remark 
\ref{R:PSS} below for additional discussion.

\smallskip

The paper is organized as follows. After the present Introduction, we present in Section \ref{S:prelim} preliminaries 
on nc functions and nc kernels
and our key tool giving the connection between a completely positive map from a $C^*$-algebra $\cA^M$ to a 
$C^*$-algebra of the form $\cL(\cY^M)$ which is also a $(\cS, \cS^*)$-bimodule map with respect to the action of 
a certain subalgebra $\cS$ of ${\mathbb C}^{M \times M}$
on the one hand, and cp nc kernels defined on a finite set of nc points $\Omega$ on the other.  Sections \ref{S:HerKer},
\ref{S:Arv-ext}, and \ref{S:ker-dom} present our results on Problems A, B, C respectively, for the special case where
the point set $\Omega$ consists of only finitely many points.  The final section \ref{S:Kurosh} shows how to extend the
finite-point results from Sections \ref{S:Arv-ext} and \ref{S:ker-dom}  to the case of a general nc point set $\Omega$
by using the abstract results of Kurosh concerning nonemptiness of inverse limits for a inverse spectrum of compacta.
For the case of the Jordan-decomposition problem for nc kernels (Problem A) treated in Section \ref{S:HerKer}, 
we also show in Section \ref{S:Kurosh} where the Kurosh formalism breaks down, as well as examples showing
that such decomposability results are not possible in general.  A final remark (Remark \ref{R:PSS}) makes precise
the connections of our kernel-domination result with some Positivstellens\"atze for the noncommutative setting
which have appeared relatively recently in the literature.

\section{Preliminaries} \label{S:prelim}
In this section we review some preliminaries from \cite[Section 2]{BMV2} concerning noncommutative (nc) functions 
and completely positive noncommutative (cp nc) kernels which will be needed in the sequel.

\subsection{Taxonomy of noncommutative sets}  \label{S:tax}
We let $\cV$ be a vector space, $\cV_{\rm nc} = \amalg_{n=1}^\infty \cV^{n \times n}$ the set of all square matrices 
over $\cV$ of arbitrary size.  Note that $\cV^{n \times m}$ is a left module over ${\mathbb C}^{n \times n}$ and a right module over ${\mathbb C}^{m \times m}$ by making use of ordinary matrix multiplication combined with the
bimodular structure of $\cV$ over ${\mathbb C}$ as a vector space over ${\mathbb C}$.  We say that a subset $\Xi$
of $\cV_{\rm nc}$ is a \textbf{nc subset} if $\Xi$ is closed under direct sums.  Following Definition 2.4
from \cite{BMV2}, we say that the nc subset $\Xi$ is a 
\textbf{full}  nc subset of $\cV_{\rm nc}$  if, in addition to being closed under direct sums, 
 $\Xi$ is \textbf{invariant under left injective intertwinings}, i.e., $Z \in \Xi_n$, $\widetilde Z \in \cV^{m \times m}$ such that
 $\cI \widetilde Z = Z \cI$ for some injective $\cI \in {\mathbb C}^{n \times m}$ (so $n \ge m$) implies that
 $\widetilde Z \in \Xi_m$.  An equivalent statement of this latter property is that $\Xi$ is \textbf{closed under restriction to 
 invariant subspaces:}  {\em whenever there is an invertible $\alpha \in{\mathbb C}^{n \times n}$ and a $Z \in \Xi$ of size
 $n \times n$ such that $\alpha^{-1} Z \alpha = \sbm{ \widetilde Z & Z_{12} \\ 0 & Z_{22} }$ with $\widetilde Z$ of size 
 $m \times m$, then $\widetilde Z$ is in $\Xi$.}
 The case $n=m$ is not excluded:  this special case of the condition gives us that any full nc subset $\Xi$ is also
 invariant under similarities.
 
Given an arbitrary subset $\Omega$ of $\cV_{\rm nc}$,  we have the following distinct notions of \textbf{envelopes}:
\begin{itemize}
\item The \textbf{nc envelope} $[\Omega]_{\rm nc}$ is the smallest superset of $\Omega$ in $\cV_{\rm nc}$ which is closed
under direct sums.

\item The \textbf{nc similarity envelope} $[\Omega]_{\rm sim}$ is the smallest superset of $\Omega$ in
$\cV_{\rm nc}$ which is closed under direct sums and similarity transforms.

\item The \textbf{full nc envelope}  $[\Omega]_{\rm full}$ is the smallest superset of $\Omega$ in $\cV_{\rm nc}$
which  is a full nc subset as defined above.
\end{itemize}

Note that each of the three properties {\em nc set}/{\em similarity-invariant nc set}/{\em full nc set} is closed
under intersections, so the each of the notions {\em smallest nc superset}/{\em smallest similarity-invariant nc superset}/{\em smallest full nc superset} containing a given subset is well defined.  For brevity we now focus on notions related
to full nc subsets as this is all that will be needed in the sequel; we leave the parallel notions concerning nc sets and similarity-invariant nc sets to the interested reader.

If ${\mathbb D}$ is another subset of $\cV_{\rm nc}$ and $\Omega'$ is a subset of ${\mathbb D}$, we say that
$\Omega'$ is a \textbf{${\mathbb D}$-relative full nc set} if
the full nc envelope $[\Omega']_{\rm full}$ intersected with ${\mathbb D}$
is again just  $\Omega'$:  $\Omega' = [\Omega']_{\rm full} \cap {\mathbb D}$.  If $\Omega$ is any subset of 
${\mathbb D}$, then the smallest ${\mathbb D}$-relative full nc set containing $\Omega$ is
$\Omega' = {\mathbb D} \cap [\Omega]_{\rm full}$.  In this case we say that the ${\mathbb D}$-relative full nc set
$\Omega'$ \textbf{is generated by} $\Omega$.  We shall be particularly interested in the case when $\Omega$
is a finite subset of ${\mathbb D}$.

\subsection{Noncommutative functions on a subset $\Omega$}

Let $\Omega$ be a subset of $\cV_{\rm nc}$ and let $\cV_0$ be another vector space.  We say that $f \colon \Omega
\to \cV_{0, {\rm nc}}$ is a \textbf{($\cV_0$-valued) nc function} if 
\begin{itemize}
\item $f$ is \textbf{graded}:  $f(Z) \in \cV^{n \times n}_{0}$ for $Z \in \Omega_n$,

\item $f$ \textbf{respects intertwinings}:  If $\alpha \in {\mathbb C}^{n \times m}$, $Z \in \Omega_n$,
$\widetilde Z \in \Omega_m$ are such that $ Z \alpha = \alpha \widetilde Z$, then
$ f(Z) \alpha = \alpha f(\widetilde Z)$.
\end{itemize}

It is known (see \cite{KVV-book}) that the "respects intertwinings" condition can be replaced by the pair of conditions:
\begin{itemize}
\item $f$ \textbf{respects direct sums}:  If $Z$, $\widetilde Z$, and $\sbm{Z & 0 \\ 0 & \widetilde Z}$ are all in $\Omega$, 
then $f\left( \sbm{ Z & 0 \\ 0 & \widetilde Z} \right) = \sbm{ f(Z) & 0 \\ 0 & f(\widetilde Z) }$, and
\item $f$ \textbf{respects similarities}:  If $Z$, $\widetilde Z$ are in $\Omega_n$, and $\alpha$ is invertible in 
${\mathbb C}^{n \times n}$
with $\widetilde Z = \alpha Z \alpha^{-1}$, then  $f( \widetilde Z ) = \alpha f(Z) \alpha^{-1}$.
\end{itemize}
Note that in these definitions we do not insist that $\Omega$ have additional structure as a subset of $\cV_{\rm nc}$
(e.g., being a nc subset, a nc similarity-invariant subset, or a full nc subset). 

\subsection{Noncommutative kernels}  \label{S:nc-ker}

Suppose $\cA$ and $\cB$ are $C^*$-algebras and that $K$ is a function from $\cV_{\rm nc} \times \cV_{\rm nc}$ into
$$
\cL(\cA, \cB)_{\rm nc} : = \amalg_{n,m \ge 1} \cL(\cA^{n \times m}, \cB^{n \times m}).
$$
We say that $K$ is a \textbf{nc kernel} if
\begin{itemize}
\item $K$ is \textbf{graded}:  $Z \in \Omega_n$, $W \in \Omega_m$ $\Rightarrow$ 
$K(Z,W) \in \cL(\cA^{n \times m}, \cB^{n \times m})$,

\item $K$ \textbf{respects intertwinings}:  
\begin{align*}
& Z \in \Omega_n,\,  \widetilde Z \in \Omega_{\widetilde n}, \, \alpha \in {\mathbb C}^{\widetilde n \times n} \text{ such that }
\alpha Z = \widetilde Z \alpha, \\
&W \in \Omega_m, \, \widetilde W \in \Omega_{\widetilde m},\, \beta \in {\mathbb C}^{\widetilde m \times m} \text{ such that }
\beta W = \widetilde W  \beta, \\
& P \in \cA^{n \times m} \Rightarrow \alpha K(Z,W)(P) \beta^*  = K(\widetilde Z, \widetilde W) (\alpha P \beta^*).
\end{align*}
\end{itemize}
As explained in \cite{BMV2} and in \cite{KVV-book}, the "respects intertwinings" condition can be replaced by a pair of 
conditions: "respects direct sums" and "respects similarities":
\begin{itemize}
\item $K$ \textbf{respects direct sums}:  for $Z \in \Omega_n$ and $\widetilde Z \in \Omega_{\widetilde n}$ such that
$\sbm{Z & 0 \\ 0 & \widetilde Z } \in \Omega_{n+m}$, $W \in \Omega_m$ and $\widetilde W \in \Omega_{\widetilde m}$
such that $\sbm{ W & 0 \\ 0 & \widetilde W } \in \Omega_{m + \widetilde m}$, and 
$P = \sbm{ P_{11} & P_{12} \\ P_{21} & P_{22} } \in \cA^{(n+m) \times (\widetilde n + \widetilde m)}$, it then holds that
$$
K\left( \sbm{ Z & 0 \\ 0 & \widetilde Z}, \, \sbm{ W & 0 \\ 0 &  \widetilde W} \right) \left( \sbm{ P_{11} & P_{12} \\ P_{21} &
P_{22} } \right) = \begin{bmatrix} K(Z,W)(P_{11}) & K(Z, \widetilde W)(P_{12}) \\ K(\widetilde Z, W)(P_{21}) &
K(\widetilde Z, \widetilde W)(P_{22}) \end{bmatrix},
$$
\item $K$ \textbf{respects similarities}: 
\begin{align*}  & Z, \widetilde Z \in \Omega_n, \alpha \in {\mathbb C}^{n \times n} \text{ invertible with } 
\widetilde Z = \alpha Z \alpha^{-1}, \\
& W, \widetilde W \in \Omega_m, \, \beta \in {\mathbb C}^{m \times m} \text{ invertible with }
\widetilde W = \beta W \beta^{-1}, \\
& P \in \cA^{n \times m} \Rightarrow K(\widetilde Z, \widetilde W)(P) =
\alpha K(Z,W)(\alpha^{-1} P \beta^{-1 *}) \beta^*.
\end{align*} 
\end{itemize}

Let us say that the nc kernel $K$ is \textbf{completely positive} (\textbf{cp}) if, for each $Z_1, \dots, Z_K \in 
\Omega$, with say $Z_i \in \Omega_{n_i}$,  the map
$ [ a_{ij} ] \mapsto [ K(Z_i, Z_j)(a_{ij}) ]$ is a completely positive map between the $C^*$-algebras 
 $\cA^{N \times N}$ and $\cB^{N \times N}$, where $N = \sum_{i=1}^K n_i$, or equivalently and more explicitly,
for all $a_i \in \cA^{n_i}$ and $b_i \in \cB^{n_i}$, it is the case that
\begin{equation}   \label{cp-nc-ker}
  \sum_{i,j = 1}^K  b_i^* K(Z_i, Z_j)(a_i^* a_j) b_j  \succeq 0.
\end{equation}
In case $\Omega$ is a nc set, one can use the "respects direct sums" property to express this last condition
more succinctly as (see \cite{BMV2}):  {\em for each $Z \in \Omega_n$,  $K(Z,Z)$ is a completely 
positive map between the $C^*$-algebras $\cA^{n \times n}$ and $\cB^{n \times n}$.} We shall be primarily interested
in the case where $K \colon \Omega \times \Omega \to \cL(\cA, \cB)_{\rm nc}$ with the $C^*$-algebra $\cB$ 
assumed to have the form $\cL(\cY)$ for some Hilbert space $\cY$;  the theory can be pushed more more generally 
by getting into a Hilbert $C^*$-module rather than Hilbert space setting (see \cite{BBLS,  Marx-dis},  but for our purposes
here the Hilbert space setting is sufficient.
With this assumption in place,  we rewrite \eqref{cp-nc-ker} as: for all  $Z_i \in \Omega_{n_i}$,
$a_i \in \cA^{n_i}$ and $y_i \in \cY^{n_i}$ for $1 \le i \le K$ we have
\begin{equation}   \label{cp-nc-ker'}
\sum_{i,j = 1}^K \langle K(Z_i, Z_j)(a_i^* a_j) y_j, y_i \rangle_{\cY^{n_i}} \ge 0.
\end{equation}
it is then known (see \cite{BMV1}) that
$K$ has a \textbf{Kolmogorov decomposition}, i.e., there is a Hilbert space $\cX$, and unital $*$-representation
$\pi \colon \cA \to \cL(\cX)$ and a nc function $H \colon \Omega \to \cL(\cX, \cY)_{\rm nc}$ so that
\begin{equation}  \label{Kol-decom}
K(Z,W)(P) = H(Z) \left( ({\rm id}_{{\mathbb C}^{n \times m}} \otimes \pi)(P) \right) H(W)^*
\end{equation}
 for $Z \in \Omega_n$, $W \in \Omega_m$ and $P \in \cA^{n \times m}$. 
 
 The following is a useful fact concerning cp nc kernels. 
 
 \begin{lemma} \label{L:zeroker}
	A cp nc kernel $K:\Omega \times \Omega \to \cL(\cA,\cB)_{\rm nc}$ is equal to the zero kernel if and only if 
	$K(Z,Z)(1_{\cA^{n \times n}})=0$ for any $Z \in \Omega_n$, $n=\bbZ_+$. 
\end{lemma}

\begin{proof}
The result follows from the fact that for any positive map $\phi$ between unital $C^*$-algebras $\norm{\phi}=\norm{\phi(1)}$ (see \cite{Paulsen}) .
\end{proof}

\subsection{Encoding of nc-kernel property via module structure}  \label{S:encoding}

Suppose that we are given a point set $\Omega = \{ Z_1, \dots, Z_N \}$ contained in an ambient universal nc set $\cV_{\rm nc}$
together with a function $K \colon \Omega \times \Omega \to {\mathcal L}({\mathcal A}_{\rm nc}, {\mathcal L}({\mathcal Y}))_{\rm nc}$.
Set $Z^{(0)}$ equal to the direct sum of all the points in $\Omega$:
$$
   Z^{(0)}  = \bigoplus_{i=1}^N Z_i.
$$
 If each point $Z_i$ say has size $n_i \times n_i$, we then see that $Z^{(0)}$ has size
$N_0 \times N_0$ where $N_0 = \sum_{i=1}^N n_i$. 
It is useful to observe that $N_0 \times N_0$ matrices  
(be they over ${\mathbb C}$, $\cA$ or $\cL(\cY)$) can be organized as block $N \times N$ matrices with block $(i,j)$ having size
$n_i \times n_j$.  We associate with the nc kernel $K$ the linear map $\phi_K$ from 
${\mathbb A}:= {\mathcal A}^{N_0 \times N_0}$
to ${\mathbb L}:= {\mathcal L}({\mathcal Y})^{N_0 \times N_0}$ given by 
\begin{equation}   \label{ker/map}
\phi_K \colon P \mapsto  K(Z^{(0)}, Z^{(0)})(P): =  [K(Z_i, Z_j)(P_{ij}) ]_{1 \le i, j \le N}
\end{equation}
if  $P = [ P_{ij}]_{1 \le i,j \le N_0}$  with  $P_{ij} \in \cA^{n_i \times n_j}$.  
We would like to understand how one can encode the property of $K$ being a nc kernel (or of being a cp nc kernel)
as a property of the map $\phi_K$ between $C^*$-algebras ${\mathbb A}$ and ${\mathbb L}$.  Toward this end,  let
$\cS$ denote the subalgebra of $N_0 \times N_0$ matrices consisting of $\alpha \in {\mathbb C}^{N_0 \times N_0}$
which intertwine the point $Z^{(0)}$ with itself: 
\begin{equation}  \label{algebraS}
\cS = \{ \alpha \in {\mathbb C}^{N_0 \times N_0} \colon \alpha Z^{(0)} = Z^{(0)} \alpha \}.
\end{equation}
Note that $\cS$ is a subalgebra but not necessarily a $*$-subalgebra of ${\mathbb C}^{N_0 \times N_0}$.  
Given a mapping $\phi \colon {\mathbb A} \to {\mathbb L}$, we say that $\phi$ is a \textbf{$(\cS, \cS^*)$-bimodule map}
if
\begin{equation}   \label{S-S*bimodulemap}
\phi( \alpha \cdot P \cdot \beta^*) = \alpha \cdot \phi(P) \cdot \beta^*
\end{equation}
for all $\alpha, \beta \in \cS$.  Note that the left-hand side in \eqref{S-S*bimodulemap} uses the ${\mathbb C}^{N_0 \times N_0}$-module
action on ${\mathbb A}$ which the right-hand side uses the ${\mathbb C}^{N_0 \times N_0}$-module action on ${\mathbb L}$.
As we are identifying $\cL(\cY)^{N \times N}$ with $\cL(\cY^N)$, it will be convenient to introduce the notation $L_\alpha$ and $\beta^*$
for the representations $\alpha \mapsto L_\alpha$, $\beta^* \mapsto L_\beta^*$ for the representations of $\cS$ and of $\cS^*$
respectively on $\cY^N$ given by 
$$
    L_\alpha y = \alpha \cdot y, \quad L_{\beta^*} y = \beta^* \cdot y \text{ for } y \in \cY^N.
$$ 
With these conventions in place we rewrite \eqref{S-S*bimodulemap} as
\begin{equation}  \label{S-S*bimodulemap'}
\phi( \alpha \cdot P \cdot \beta^*) = L_\alpha \,  \phi(P) \, L_{\beta^*}.
\end{equation}
Then we have the following result.

\begin{proposition}   \label{P:ker/map} {\rm (See \cite{Marx-dis}.)} Suppose that the map $K$ from $\Omega \times \Omega$ to
$\cL(\cA_{\rm nc}, \cL(\cY)_{\rm nc})$ and $\phi_K$ from ${\mathbb A}$ to  ${\mathbb L}$  are related as in \eqref{ker/map}.
Then:

\smallskip

\noindent
{\rm \textbf{(1)}} $K$ is a nc kernel on $\Omega$ if and only if $\phi_K$ is a $(\cS, \cS^*)$-bimodule map.
Conversely, if $\phi \colon \cA^{N \times N} \to \cL(\cY)^{N \times N}$ is a $(\cS, \cS^*)$-bimodule map, there is a
uniquely determined nc kernel $K$ so that $\phi$ has the form $\phi = \phi_K$.

\smallskip

\noindent
{\rm \textbf{(2)} }  $K$ is a cp nc kernel on $\Omega$ if and only if $\phi_K$ is a $(\cS, \cS^*)$-bimodule map
which is also cp.  Conversely, if  $\phi \colon {\mathbb A} \to {\mathbb L}$ is a cp $(\cS, \cS^*)$-bimodule map,
there there is a uniquely determined cp nc kernel $K$ on $\Omega$ so that $\phi = \phi_K$.
\end{proposition}

\begin{proof}
Suppose first that $K$ is a nc kernel on $\Omega$.  
Define $\widetilde K(Z^{(0)}, Z^{(0)})$ as in \eqref{ker/map} and also define
\begin{align*}
\widetilde K(Z^{(0)}, Z_j)(P) & = \begin{bmatrix} K(Z_1, Z_j)(P_{1}) \\ \vdots \\ K(Z_N, Z_j)(P_N) \end{bmatrix} \text{ if } 
P = \begin{bmatrix} P_1 \\ \vdots \\ P_N \end{bmatrix},  \\
\widetilde K(Z_i, Z^{(0)})(P) & = \begin{bmatrix} K(Z_i, Z_1)(P_1) & \cdots & K(Z_i, Z_N)(P_N)  \end{bmatrix} \\
&  \quad  \text{ if }
P = \begin{bmatrix} P_1 & \cdots & P_N \end{bmatrix}, \\
 \widetilde K(Z,W)(P)  & = K(Z, W)(P) \text{ if } Z \in \Omega_n, \, W \in \Omega_m,  \, P \in \cA^{n \times m}.
\end{align*}
Then by construction $\widetilde K$ as a function on $\widetilde \Omega \times \widetilde \Omega$ is an extension of the function $K$
defined on $\Omega \times \Omega$.
It is clear from the construction that $\widetilde K$ is a graded kernel.  

Our next goal is to show that $\phi$ is a left $\cS$-module map if and only if $K$ satisfies the left intertwining condition:
\begin{align}
& \alpha \in {\mathbb C}^{i_0 \times k_0},\, \alpha  Z_{k_0} = Z_{i_0}  \alpha, \, P \in \cA^{n_{k_0} \times n_j}  \notag \\
& \quad  \Rightarrow
L_\alpha K(Z_{k_0}, Z_j)(P) = K(Z_{i_0}, Z_j) (\alpha  P) \text{ for }  1 \le i_0, k_0, j \le N.
\label{leftintertwine}
\end{align}
Toward this end, let us assume first that $\phi$ is a left $\cS$-module map.  Let us write a $N_0 \times N_0$ complex matrix
$\alpha$ as a block $N \times N$ matrix $\alpha = [ \alpha_{ij}]_{i,j=1, \dots, N}$ where the entry $\alpha_{ij}$ has size
$n_i \times n_j$. Similarly we write  a matrix $P \in \cA^{N_0 \times N_0}$ as $P = [P_{ij}]_{1 \le i,j \le N}$ where
the block entry $P_{ij}$ is in $\cA^{n_i \times n_j}$.  Then we see that 
\begin{align*}
& \phi( \alpha \cdot P) = \bigg[ K(Z_i, Z_j) \big( \sum_{k=1}^N \alpha_{ik} P_{kj} \big) \bigg]_{1 \le i,j \le N},  \\
& L_\alpha \phi(P) = \bigg[ \sum_{k=1}^N \alpha_{ik} K(Z_k, Z_j) (P_{kj}) \bigg]_{1 \le i,j \le N}.
\end{align*}
By definition $\phi$ being a left $\cS$-module map means that
$$
\alpha Z^{(0)} = Z^{(0)}  \alpha \Rightarrow \phi( \alpha  P ) = L_\alpha \phi(P),
$$
or, in more detail,
\begin{align}
&  \alpha_{ik} Z_k = Z_i \alpha_{ik} \text{ for all } k \Rightarrow  \notag  \\
& \quad
\bigg[ K(Z_i, Z_j)\big(\sum_{k=1}^N \alpha_{ik} P_{kj}\big) \bigg]_{ij} =
\bigg[ \sum_{k=1}^N \alpha_{ik} K(Z_k, Z_j)(P_{kj}) \bigg]_{ij}.
\label{phi=leftmodmap}
\end{align}
Let us specialize this to the case where $\alpha$ has only one nonzero block-entry:
$$
 \alpha_{i j} = \delta_{i i_0} \delta_{k k_0} \alpha_{i_0 k_0} \text{ for some } n_{i_0} \times n_{k_0}\text{-block matrix }
 \alpha_{i_0k_0}
$$
where $\delta_{ii_0}$ and $\delta_{k k_0}$ are Kronecker deltas.  Then we see that the intertwining condition
$\alpha_{ik} Z_k = Z_i \alpha_{ik}$ is satisfied as long as 
\begin{equation}   \label{intertwiner}
 \alpha_{i_0 k_0} Z_{k_0} = Z_{i_0} \alpha_{i_0k_0}
\end{equation}
since the remaining conditions ($\alpha_{ij} Z_j= Z_i \alpha_{ij}$) are automatic in the form $0 = 0$
when $(i,j) \ne (i_0, j_0)$.  In this case  condition \eqref{phi=leftmodmap} works out to be
$$
\left[ \delta_{i i_0} K(Z_{i_0}, Z_j) (\alpha_{i_0 k_0} P_{k_0 j)}\right]_{ij}
 = \left[ \delta_{i i_0} \alpha_{i_0 k_0} K(Z_{k_0}, Z_j)(P_{k_0 j} \right]_{ij}.
$$
In particular we have equality of the $(i_0, j)$-entries:
$$
K(Z_{i_0}, Z_j) (\alpha_{i_0 k_0} P_{k_0 j}) =  \alpha_{i_0 k_0} K(Z_{i_0}, Z_j)(P_{k_0 j})
$$
As the indices $i_0, k_0$ are arbitrary and the matrix $\alpha_{i_0 k_0}$ is arbitrary subject  to the condition
\eqref{intertwiner}, we can now conclude that indeed $K$ satisfies the left intertwining conditions \eqref{leftintertwine}.
That $K$ also satisfies the right intertwining conditions
\begin{align} 
& \beta \in {\mathbb C}^{k_0 \times j_0}, \, \beta Z_{k_0} = Z_{j_0} \beta, \, P \in \cA^{n_i \times n_{j_0}} \notag \\
& \quad \Rightarrow K(Z_i, Z_{j_0})(P) \beta^* = K(Z_i, Z_{k_0})(P \beta^*)
\label{rightintertwine}
\end{align}
follows similarly by using the assumption that $\phi$ is also a right $\cS^*$-module.

Conversely suppose that $K$ satisfies the left intertwining conditions \eqref{leftintertwine} and we seek to verify that $\phi$
is a left $\cS$-module map.  Thus by assumption we know that the implication \eqref{leftintertwine} holds for each 
triple of indices $i_0, k_0, j$ and we seek to verify condition \eqref{phi=leftmodmap}.  The premise for \eqref{phi=leftmodmap}
is that  we are given $\alpha = [\alpha_{ik}]_{i,k}$ such that $\alpha_{ik} Z_k = Z_i \alpha_{ik}$.  As by assumption
$K$ satisfies \eqref{leftintertwine}, this implies that $\alpha_{ik} K(Z_k, Z_j)(P_{kj}) =K(Z_i, Z_j)( \alpha_{ik} P_{kj})$ for
each fixed $i,j,k$.  It now suffices to sum over $k$ from $1$ to $N$ to arrive at  \eqref{phi=leftmodmap} as wanted.
One can show that $K$ satisfying the right intertwining conditions \eqref{rightintertwine} implies that
$\phi$ is a right $\cS^*$-module map by a dual argument.

We now suppose that we are only given $\phi \colon {\mathbb A} \to {\mathbb L}$ which is a $(\cS, \cS^*)$-bimodule map.
Let us let $Q_i$ be the $N \times N$ matrix over ${\mathbb C}$ corresponding to the projection to the $i$-th block, i.e.,
$Q_i$ is the block diagonal matrix with only nonzero block diagonal entry equal to the identity matrix $I_{n_i}$ in the $i$-th block
Then it is easily checked that $Q_i \in \cS \cap \cS^*$, and hence $\phi(Q_i  \cdot  P) = L_{Q_i} \phi(P)$
and $\phi( P \cdot Q_j) = \phi(P) L_{Q_j}$ for $ 1 \le i,j \le N$.  From this property one can deduce that
the map $\phi$ then must  have the form
$$
\phi \big( [ P_{ij} ]_{1 \le i,j \le N} \big) = [ \phi_{ij} (P_{ij}) ]_{1 \le i,j \le N}
$$
for linear maps $\phi_{ij} \colon \cA^{n_i \times n_j} \to \cL(\cY)^{n_i \times n_j}$.  We then define
$K \colon \Omega \times \Omega \to \cL(\cA_{\rm nc}, \cL(\cY)_{\rm nc})$ by
\begin{equation}  \label{defK}
   K(Z_i, Z_j)(P_{ij}) = \phi_{ij}(P_{ij}).
\end{equation}
Now it is a simple bookkeeping exercise to check that the $(\cS, \cS^*)$-bimodule property of $\phi$ is exactly
what is needed for $K$ so defined to be a nc kernel on the finite set $\Omega$.

\smallskip

We next analyze the claim regarding complete positivity.
Suppose that $K$ is a nc kernel on $\Omega$.
When it is assumed that $\Omega$ is a nc set (i.e., invariant under formation of direct sums), it is often taken as the definition
of $K$ being completely positive simply that the map $K(Z,Z)$ is a positive map for all $Z \in \Omega$
(see Proposition 2.2 in \cite{BMV1}).  In case $\Omega$ is  finite set $\{ Z_1, \dots, Z_N\}$ augmented by the single point
$Z^{(0)} = \bigoplus_{1}^N Z_i$,  the adjustment of these observations is that $K$ is a c.p. kernel if and only if
$K(Z^{(0)}, Z^{(0)})$ is a cp map from ${\mathbb A}$ to ${\mathbb L}$, i.e., if and only if $\phi$ is a completely positive map.

Conversely, if $\phi$ is a cp $(\cS, \cS^*)$-bimodule map and $K$ is defined as in \eqref{defK},  one can check that the
complete positivity of $\phi$ is all that is required to guarantee the complete positivity of $K$ as a nc kernel.
\end{proof}

The next result characterizes the $(\cS, \cS^*)$-bimodule property for a cp map $\phi$ from $\cA^{M\times M}$
to $\cL(\cY^{M})$ in terms of a Stinespring representation
\begin{equation}   \label{Stinespring-rep}
\phi(P) = V \pi(P) V^* \text{ for } P \in \cA^{N \times N}.
\end{equation}
Here $V$ is an operator from $\cX$ to $\cY^{N_0}$ and $\pi \colon \cA^{N_0 \times N_0} \to \cL(\cX)$ is a $*$-representation
of $\cA^{N_0 \times N_0}$ on a Hilbert space $\cX$.

\begin{theorem}  \label{T:bimodule-Stinespring}
Suppose that  $\cS$ be a subalgebra of ${\mathbb C}^{M \times M}$;  use the natural ${\mathbb C}^{M \times M}$-bimodule
structure of $\cA^{M \times M}$ to also view $\cA^{M \times M}$ as a $(\cS, \cS^*)$-bimodule, and similarly for 
$\cL(\cY)^{M \times M} \cong \cL(\cY^M)$.  Let $\phi$ be a cp map from $\cA^{M \times M}$ to $\cL(\cY^M)$ with
Stinespring representation \eqref{Stinespring-rep}.
Then $\phi$ is a $(\cS, \cS^*)$-bimodule map if and only if
\begin{equation}  \label{Vmod}
   \pi(\beta^*) V^* =  V^* L_{\beta^*} \text{ for all } \beta \in \cS,
\end{equation}
or equivalently,
\begin{equation}  \label{Vmod'}
L_\alpha V = V \pi(\alpha) \text{ for all } \alpha \in \cS.
\end{equation}
\end{theorem}

\begin{proof}
Note first that \eqref{Vmod} and \eqref{Vmod'} follow from each other by taking adjoints.

Suppose next that \eqref{Vmod} and hence also \eqref{Vmod'} hold.  Then we compute
\begin{align*}
\phi(\alpha \cdot P \cdot \beta^*) & = V \pi(\alpha \cdot P \cdot \beta^*) V^*  = V \pi(\alpha) \pi(P) \pi(\beta^*) V^*  \\
& = L_\alpha V \pi(P) V^* L_{\beta^*} \text{ (by \eqref{Vmod} and \eqref{Vmod'})} \\
& = L_\alpha \phi(P) L_{\beta^*}
\end{align*}
and we conclude that $\phi$ is a $(\cS, \cS^*)$-bimodule map.

Conversely, suppose that $\phi$ is a $(\cS, \cS^*)$-bimodule map.  Then we compute
\begin{align*}
& \| \pi(\beta^*) V^* y - V^* L_{\beta^*}  y \|^2  = 
\langle y, V \pi(\beta \beta^*) V^* y \rangle - \langle y, V \pi(\beta) V^* L_{\beta^*} y \rangle \\
& \quad \quad  - \langle y, L_\beta  V \pi(\beta^*) V^* y \rangle + \langle y, L_\beta  VV^*L_{\beta^*} y) \rangle \\
& = \langle y, \phi(\beta \beta^*)y - \phi(\beta) L_{\beta^*}y) - L_\beta \phi(\beta^*) y + L_\beta 
\phi(1_{\cA^{N_0 \times N_0}}) L_{\beta^*} y \rangle \\
& = \langle y, \phi(\beta \beta^*) y - \phi(\beta \beta^*) y - \phi(\beta \beta^*) y + \phi(\beta \beta^*)y \rangle  = 0
\end{align*}
for all $y \in \cY^{N_0}$, and \eqref{Vmod} follows.

Alternatively, if we assume that $\cS$ is the intertwining algebra associated with a finite nc point set $\Omega =
\{ Z_1, \dots, Z_N\}$ and $\phi$ has the form $\phi_K$ coming from a cp nc kernel $K$ on $\Omega$ with values
equal to operators from $\cA_{\rm nc}$ to $\cL(\cY)_{\rm nc}$, we may arrive at \eqref{Vmod} as follows.
Here $M = N_0$ where $N_0 = \sum_{j=1}^N n_j$ where $N$ is the number of points in $\Omega$ and $Z_i$ has size
$n_i \times n_i$.
We let
$$
  K(Z,W)(P) = H(Z) ( {\rm id}_{n \times m} \otimes \pi) (P) H(W)^*
$$
be the Kolmogorov decomposition for the cp nc kernel $K$, where $\pi \colon \cA \to \cL(\cX)$ is a unital $*$-representation
of $\cA$.  Then we get a Stinespring representation for the
cp map $\phi_K$ via
\begin{align*}
   \phi_K([P_{ij}]) &  = H(Z^{(0)}) ({\rm id}_{N \times N} \otimes \pi)([P_{ij}])  H(Z^{(0)})^* \\
   & = V \Pi([P_{ij}]) V^*
\end{align*}
where $[P_{ij}]$ is a $N \times N$-block matrix with $(i,j)$-block entry of size $n_i \times n_j$ over $\cA$,
$\Pi = {\rm id}_{N \times N} \otimes \pi$ is a unital representation of $\cA^{N_0 \times N_0}$ on $\cX^{N_0}$
partitioned as
\begin{equation}   \label{partitioning}
\cX^{N_0} = \begin{bmatrix}  \cX^{n_1} \\ \vdots \\ \cX^{n_N}  \end{bmatrix}.
\end{equation}
As $H$ is a nc function, it follows that $H(Z^{(0)})^* L_{\beta^*} = L_{\beta^*} H(Z^{(0)})^*$ where 
$L_{\beta^*}$ on the right-side of this equality indicates the action of $\cS \subset {\mathbb C}^{N \times N}$
on $\cX^{N_0}$ via matrix multiplication using the partitioning \eqref{partitioning}.  Now it is a matter of bookkeeping to check that
this is exactly the right structure required to verify that, since $\pi$ is a unital representation, 
$$
    L_{\beta^*} = ({\rm id}_{N \times N} \otimes \pi)(\beta^* \cdot I_\cA).
$$
where $\beta^* \mapsto \beta^* \cdot I_\cA$ is the identification of $\cS$ with an element of $\cA^{N \times N}$.
and the formula \eqref{Vmod} follows.  One can argue that this second approach to proving \eqref{Vmod} gives a model
for how one can construct $(\pi, V)$ for which \eqref{Vmod} holds.
\end{proof}

Observing the proof of the converse direction in Theorem \ref{T:bimodule-Stinespring} leads to the following corollary.

\begin{corollary}  \label{C:bimodule-Stinespring}
Suppose that ${\mathbb S}$ is an operator system contained in $\cA^{M \times M}$ which is also a $(\cS, \cS^*)$-sub-bimodule
of $\cA^{M \times M}$ and $\phi_0 \colon {\mathbb S} \to \cL(\cY^M)$ is a cp map which is also a $(\cS, \cS^*)$-bimodule
map.  Let $\phi \colon \cA^{M \times M} \to \cL(\cY^M)$ be any cp map which extends $\phi_0$:
$$
 \phi(P) = \phi_0(P) \text{ if } P \in {\mathbb S}.
$$
Then $\phi$ is also a $(\cS, \cS^*)$-bimodule map.
\end{corollary}

\begin{proof}  The proof of the converse direction in Theorem \ref{T:bimodule-Stinespring} ($\phi = $ a bimodule map implies
the intertwining condition \eqref{Vmod}) only uses that $\phi_0:= \phi|_{\mathbb S}$ is a $(\cS, \cS^*)$-bimodule map.
Thus the result of the calculation (the intertwining condition \eqref{Vmod}) holds for the Stinespring representation for any cp map
$\phi$ extending $\phi_0$.  By the calculation used to prove that conditions \eqref{Vmod}-\eqref{Vmod'} imply the $(\cS, \cS^*)$-bimodule
property for $\phi$,  it then follows that any cp extension $\phi$ of $\phi_0$ to $\cA^{M \times M}$
is in fact also a $(\cS, \cS^*)$-bimodule map.
\end{proof}

\section{Problem A: The decomposability problem for Hermitian noncommutative kernels/Hermitian maps} \label{S:HerKer}
We are given a nc kernel $K$ and seek to show that it has a decomposition \eqref{4fold} with each $K_j$ ($j=1,2,3,4$) a
cp nc kernel.    We shall first go through a series of reductions.

Given a nc kernel $K \colon \Omega \times \Omega \to \cL(\cA_{\rm nc} , \cL(\cY)_{\rm nc})$, we define the
adjoint kernel $K^*$ also from $\Omega \times \Omega$ to $\cL(\cA_{\rm nc}, \cL(\cY)_{\rm nc})$ by
$$
K^*(Z,W)(P) = K(W,Z)(P^*)^*.
$$
We have several observations concerning this adjoint operator on nc kernels.

\begin{proposition}  \label{adjoint-ker}
\textbf{\rm (1)}  If $K$ is nc kernel, then $K^*$ is also a nc kernel.

\smallskip

\textbf{\rm (2)}  If $K$ is a cp nc kernel, then $K$ is Hermitian.

\smallskip

\textbf{\rm (3)}  Any nc kernel $K$ can be decomposed as $K = K_R + i K_I$ where $K_R$ and $K_I$ are Hermitian kernels.
Hence, to prove Conjecture A, it suffices to show that any Hermitian nc kernel  $K$ can be written as the difference of two
cp kernels:
$$
  K = K^* \Rightarrow K = K_+ - K_- \text{ with } K_+, \, K_- \text{ equal to cp kernels.}
$$
\end{proposition}

\begin{proof}   \textbf{(1)}  Suppose that $Z \in \Omega_n$, $\widetilde Z \in \Omega_{\widetilde n}$,
$\alpha \in {\mathbb C}^{\widetilde n \times n}$ and $\alpha Z = \widetilde Z \alpha$, and $K$ is a nc kernel on $\Omega$.
Then we compute
\begin{align*}
\alpha K^*(Z, W)(P) & = \alpha K(W,Z)(P^*)^* = \big( K(W,Z)(P^*) \alpha^* \big)^*   \\
& = K^*(\widetilde Z, W)(\alpha P) 
\end{align*}
where the {\em respects intertwinings} property of $K$ with respect to the second argument is used in the last step.
In this way we verify the {\em respects intertwinings} property of $K^*$ with respect to the first argument.
By a similar argument, the {\em respects intertwinings} property of $K$ with respect to the first argument can be used to prove the
{\em respects intertwinings} property of $K^*$ with respect to the second argument.

\smallskip
 \textbf{(2)}   Given $Z \in \Omega_n$, $W \in \Omega_m$,
$P_{12} \in \cA^{n \times m}$, choose $P_{11} \in \cA^{n \times n}$, $P_{12} \in \cA^{n \times m}$
and $P_{22} \in \cA^{n \times n}$ so that $P = \sbm{ P_{11} & P_{12} \\ P_{12}^* & P_{22} }$ is positive in $\cA^{(m+n) \times (m+n)}$.
If $K$ is cp, then we must have
\begin{align*}
0 & \preceq   K\big( \sbm{ Z & 0 \\ 0 & W}, \sbm{ Z & 0 \\ 0 & W} \big)\big( \sbm{ P_{11} & P_{12} \\ P_{12}^* & P_{22} } \big) \\
& = \sbm{ K(Z,Z)(P_{11}) & K(Z,W)(P_{12}) \\ K(W,Z)(P_{12}^*) & K(W,W)(P_{22}) }.
\end{align*}
In particular it follows that the $(2,1)$ entry is the adjoint of the $(1,2)$ entry:
$$
  K(W,Z)(P_{12}^*) = K(Z,W)(P_{12})^*.
$$
where $K(W,Z)(P_{12}^*)^* = K^*(Z,W)(P_{12})$.  As $P_{12}$ is arbitrary, we conclude that $K^* = K$ as claimed.

\smallskip

\textbf{(3)}  Note that the operation of forming the adjoint kernel is an involution and is conjugate linear: 
$$ 
K^{**} = K, \quad (a K_1 + b K_2)^* = \overline{a} K_1^* + \overline{b} K_2^* \text{ for } a, b \in {\mathbb C}.
$$
Hence for any nc kernel,  $K_R:= \frac{1}{2} (K + K^*)$ and $K_I:= \frac{1}{2 i } (K - K^*)$ are Hermitian kernels and we recover
$K$ from $K_R$ and $K_I$ as
$$
 K = \frac{1}{2} (K + K^*) + i  \, \frac{1}{2i} (K - K^*) = K_R + i K_I.
$$
\end{proof}

The next result gives an equivalent formulation on deciding if a Hermitian kernel is decomposable.
These results are modelled on the corresponding results in \cite{BDT} for the case where 
$\Omega = \Omega_1$ (all points in $\Omega$ are $1 \times 1$ matrices over $\cV$).

\begin{proposition}  \label{P:HerKer-criterion}  Let $K$ be a Hermitian nc kernel.  Then the following are equivalent:

\smallskip
\textbf{\rm (1)} $K$ is decomposable, i.e., there exist cp nc kernels $K_1$ and $K_2$ so that $K = K_1 - K_2$.

\textbf{\rm (2) }  $K$ has a Hermitian Kolmogorov decomposition:
for $Z \in \Omega_n$, $W \in \Omega_m$ and $P \in \cA^{n \times m}$, 
\begin{equation}   \label{HerKolDecom}
K(Z,W)(P) = H(Z) \left( ({\rm id}_{n \times m} \otimes \pi)(P) \cdot ({\rm id}_{n\times m} \otimes J)  \right) H(W)^*
\end{equation}
where $\pi \colon \cA \to \cL(\cX)$ is a $*$-representation and $J = J^* = J^{-1} \in \cL(\cX)$ commutes with the range of $\pi$:
$$
a \in \cA \Rightarrow \pi(a) J = J \pi(a).
$$

\smallskip

\textbf{\rm (3)}  There exist cp nc kernels ${\mathbb L}_1$ and ${\mathbb L}_2$ on $\Omega$ with values in
$\cL(\cA_{\rm nc}, \cL(\cY)_{\rm nc})$ so that the kernel
$$
 {\mathbb K} \colon \Omega \times \Omega \to \cL(\cA_{\rm nc}, \cL(\cY)^{2 \times 2}_{\rm nc})
 $$
given by, for $Z \in \Omega_n$, $W \in \Omega_m$, $P \in \cA^{n \times m}$,  
\begin{equation}   \label{bigK}
 {\mathbb K}(Z,W)(P) = \begin{bmatrix} {\mathbb L}_1(Z, W)(P) & K(Z,W)(P) \\ K(Z,W)(P) & {\mathbb L}_2(Z,W)(P) \end{bmatrix}
 \end{equation}
 is cp.
\end{proposition}

\begin{proof}  \textbf{(1) $\Leftrightarrow$ (2)}    Suppose that $K$ is decomposable:  $K = K_1 - K_2$
for cp nc kernels $K_1$ and $K_2$.  Then each of $K_1$ and $K_2$ has a Kolmogorov decomposition:  for $Z \in \Omega_n$, $W \in \Omega_m$ and $P \in \cA^{n \times m}$ we have
$$
K_j(Z,W)(P) = H_j(Z) ({\rm id}_{n \times m} \otimes \pi_j)(P) H_j(W)^*
$$
for a nc function $H_j \colon \Omega \to \cL(\cX_j, \cY)_{\rm nc}$ and a $*$-representation $\pi_j \colon \cA \to \cL(\cX_j)$
for $j=1,2$.  Then we see that
\begin{align*}
& K(Z,W)(P)  = K_1(Z,W)(P) - K_2(Z,W)(P)   \\
& = \begin{bmatrix} H_1(z) & H_2(Z) \end{bmatrix}
\big( {\rm id}_{n \times m} \otimes \sbm{ \pi_1 & 0 \\ 0 & \pi_2 }  \big)(P)
\big( {\rm id}_{n \times m} \otimes \sbm{ I_{\cX_1} & 0 \\ 0 & -I_{\cX_2}} \big)
\begin{bmatrix} H_1(W)^* \\ H_2(W)^* \end{bmatrix}  \\
& = H(Z) ({\rm id}_{n \times m} \otimes \pi)(P) ({\rm id}_{n \times m} \otimes J) H(W)^*,
\end{align*}
where we set 
$$
 \cX = \sbm{ \cX_1 \\ \cX_2 }, \quad H(Z) = \begin{bmatrix} H(Z_1) & H(Z_2) \end{bmatrix}, \quad
 \pi = \sbm{ \pi_1 & 0 \\ 0 & \pi_2 }, \quad J = \sbm{I_{\cX_1} & 0 \\ 0 & -I_{\cX_2} }.
$$
has a Hermitian Kolmogorov decomposition \eqref{HerKolDecom}.

\smallskip

Conversely, if $K$ has a Hermitian Kolmogorov decomposition \eqref{HerKolDecom}, we can choose an orthogonal
decomposition of $\cX$ as $\cS = \sbm{ \cX_1 \\ \cX_2}$ with respect to which $J$ has matrix representation $J =
\sbm{ I_{\cX_1} & 0 \\ 0 & -I_{\cX_2} }$.  As by definition of Hermitian Kolmogorov decomposition  $J$ commutes with $\pi$,
it also follows that $\pi$ has a block diagonal form $\pi(a) = \sbm{\pi_1(a) & 0 \\ 0 & \pi_2(a) }$ for $*$-representations
$\pi_1$ and $\pi_2$ of $\cA$ on $\cX_1$ and $\cX_2$ respectively. Furthermore with respect to this decomposition of $\cX$
as $\cX = \sbm{ \cX_1 \\ \cX_2}$ we have a matrix representation of $H(Z)$ as $H(Z) = \begin{bmatrix} H_1(z) & H_2(Z) \end{bmatrix}$
for noncommutative functions $H_j \colon \Omega \to \cL(\cX_j, \cY)$ for $j=1,2$, we see that
$$
K(Z,W)(P) = H_1(Z) ({\rm id}_{n \times m} \otimes \pi_1)(P) H_1(W)^* 
  - H_2(Z) ({\rm id}_{n \times m} \otimes \pi_2)(P) H_2(W)^*
$$
is decomposable.

\smallskip

\textbf{(2) $\Leftrightarrow$ (3)}  Assume that $K$ has a Hermitian Kolmogorov decomposition \eqref{HerKolDecom}.
Define ${\mathbb L}_1 = {\mathbb L}_2 = : {\mathbb L}$ where ${\mathbb L}$ is given by
$$
  {\mathbb L}(Z,W)(P) = H(Z) \left( ({\rm id}_{n \times m} \otimes  \pi )(P) \right) H(W)^*.
$$
Then ${\mathbb L}$ is defined via a Kolmogorov decomposition and hence is a cp nc kernel..  Define ${\mathbb K}$
as in \eqref{bigK}.  Then
\begin{align}
& {\mathbb K}(Z,W)(P) =    \notag \\
& {\mathbb H}_0(Z)
 \begin{bmatrix} ({\rm id}_{n \times m} \otimes \pi)(P) & ({\rm id}_{n \times m} \otimes \pi)(P)
 \cdot {\rm id}_{n \times m} \otimes J \\
({\rm id}_{n \times m} \otimes \pi)(P) \cdot {\rm id}_{n \times m} \otimes J & ({\rm id}_{n \times m} \otimes \pi)(P)  \end{bmatrix} 
{\mathbb H}_0(W)^*
\label{K-prelim}
\end{align}
where we set
$$
 {\mathbb H}_0(Z) = \begin{bmatrix} H(Z) & 0 \\ 0 & H(Z) \end{bmatrix}.
$$
Let us factor the middle factor as
\begin{align*}
&  \begin{bmatrix} ({\rm id}_{n \times m} \otimes \pi)(P) & ({\rm id}_{n \times m} \otimes \pi)(P)
 \cdot {\rm id}_{n \times m} \otimes J \\
({\rm id}_{n \times m} \otimes \pi)(P) \cdot {\rm id}_{n \times m} \otimes J & ({\rm id}_{n \times m} \otimes \pi)(P)  \end{bmatrix} \\
& =
\begin{bmatrix} {\rm id}_{n \times m} \otimes \pi & 0 \\ 0 & {\rm id} _{n \times m} \otimes \pi \end{bmatrix}( P)
\cdot {\rm id}_{n \times m} \otimes \begin{bmatrix} I & J \\ J & I \end{bmatrix}.
\end{align*}
If we note the factorization
$$
 \begin{bmatrix} I & J \\ J & I \end{bmatrix} = \begin{bmatrix} I \\ J \end{bmatrix} \begin{bmatrix} I & J^* \end{bmatrix}
$$
and use the fact that
$$
 ({\rm id}_{n \times m} \otimes \pi)(P) ({\rm id}_{n \times m} \otimes J) =
({\rm id}_{n \times m} \otimes J)  ({\rm id}_{n \times m} \otimes \pi)(P) 
$$
which we know as a consequence of  \eqref{HerKolDecom} being a Hermitian Kolmogorov factorization for $K$,
we can continue  the computation \eqref{K-prelim} as
$$
 {\mathbb K}(Z,W)(P) =  {\mathbb H}(z) \left( {\rm id}_{n \times m} \otimes \pi)(P) \right) {\mathbb H}(W)^*
$$
where now we set
$$
  {\mathbb H}(z) = \begin{bmatrix} H(Z) \\ H(Z) ({\rm id}_{n \times m} \otimes J) \end{bmatrix},
$$
thereby exhibiting a Kolmogorov decomposition for ${\mathbb K}$.  Thus ${\mathbb K}$ so defined is a cp nc kernel.

\smallskip

Conversely, suppose that one can find two cp nc kernels ${\mathbb L}_1$ and ${\mathbb L}_2$ so that the kernel
\eqref{bigK} is cp.  Then ${\mathbb K}$  has a Kolmogorov decomposition
$$
{\mathbb K}(Z,W)(P) = \begin{bmatrix} H_1(Z) \\ H_2(Z) \end{bmatrix}
\left( ({\rm id}_{n \times m} \otimes \pi)(P)  \right) \begin{bmatrix} H_1(W)^* & H_2(W)^* \end{bmatrix}. 
$$
 Thus $K$ is given by
$$
 K(Z,W)(P) = H_1(Z) \left( ({\rm id}_{n \times m} \otimes \pi)(P) \right) H_2(W)^*.
$$
Rewrite this as
$$
K(Z,W)(P) = \begin{bmatrix} H_1(Z) & H_2(Z) \end{bmatrix} 
\begin{bmatrix} 0 & {\rm id}_{n \times m} \otimes \pi(P) \\ 0 & 0 \end{bmatrix}
\begin{bmatrix} H_1(W)^* \\ H_2(W)^* \end{bmatrix}.
$$
As we are assuming that $K$ is a Hermitian kernel, we also have $K = K^*$ so
$$
K(Z,W)(P) = K^*(Z,W)(P) =
\begin{bmatrix} H_1(Z) & H_2(Z) \end{bmatrix} \begin{bmatrix} 0 & 0 \\  ({\rm id}_{n \times m} \otimes \pi)(P) & 0 \end{bmatrix}
\begin{bmatrix} H_1(W)^* \\ H_2(W)^* \end{bmatrix}.
$$
Taking the average of these two representations for $K$ leads to a third representation
$$
K(Z,W)(P)  = {\mathbb H}'(Z) \left( ({\rm id}_{n \times m} \otimes \sbm{ \pi & 0 \\ 0 & \pi })(P) \cdot
( {\rm id}_{n \times m} \otimes \sbm{ 0 & I \\ I & 0 }) \right) {\mathbb H}'(W)^*
$$
where now
$$
 {\mathbb H}'(Z) = \frac{1}{\sqrt{2}}  \begin{bmatrix} H_1(Z) & H_2(Z) \end{bmatrix}
$$
and where
$$
  \begin{bmatrix} \pi(P) & 0 \\ 0 & \pi(P) \end{bmatrix} \begin{bmatrix} 0 & I \\ I & 0 \end{bmatrix} =
\begin{bmatrix} 0 & I \\ I & 0 \end{bmatrix}     \begin{bmatrix} \pi(P) & 0 \\ 0 & \pi(P) \end{bmatrix}.
$$
We thus see that $K$ has a Hermitian Kolmogorov decomposition as wanted.
\end{proof}

\begin{remark}  \label{R:nc-RKKS}
The study of cp nc kernels was launched only recently (see \cite{BMV1}).  There it is shown that besides characterizations
via the complete positivity condition \eqref{cp-nc-ker'} and via the Kolmogorov decomposition \eqref{Kol-decom},
there is a third characterization:  {\sl there is a Hilbert space $\cH$ consisting of nc functions $f \colon \Omega \to
\cL(\cA, \cY)_{\rm nc}$ such that $K$ serves as the reproducing kernel  for $\cH$ in the following sense:
given $W \in \Omega_m$, $v \in \cA^m$, $y \in  \cY^m$, the function
$K_{W,a,y} \colon \Omega \to \cL(\cA, \cY)_{\rm nc}$ given by
$$
  K_{W, v, y}(Z) \colon u \mapsto K(Z,W)(u v^*) y 
$$
belongs to $\cH$ and reproduces the value of $f \in \cH$ at the point $W \in \Omega_m$ evaluated at $v \in \cA^m$
in direction $y \in \cY^m$ according to the formula}
$$
  \langle f, K_{W,v,y} \rangle_\cH = \langle f(W)v, y \rangle_{\cY^m}
$$
(see \cite[Theorem 3.1]{BMV1}).    We note that the decomposition \eqref{HerKolDecom} can be viewed as an 
indefinite-metric analogue of the Kolmogorov decomposition \eqref{Kol-decom} for a cp nc kernel.
One can then expect that the decomposition \eqref{HerKolDecom} is equivalent to indefinite-metric analogues
of the complete positivity condition \eqref{cp-nc-ker'} and that such a $K$ should be the ``reproducing kernel"
for a nc reproducing kernel Krein space $\cH(K)$. Indeed, 
for the classical Aronszajn setting, investigation
of such reproducing kernel Krein spaces, in particular for the case where the Krein space carries only finitely many
negative squares in which case they are called reproducing kernel Pontryagin spaces, has been an active area of 
research over the past several decades (see e.g.~\cite{ADRdS}).   We leave this topic as a possible direction for future
research in the theory of Krein spaces and the associated operator theory.
\end{remark}

We have seen in Proposition \ref{P:ker/map} that problems concerning nc kernels on a finite-point set $\Omega$ 
with values mapping operators
from $\cA_{\rm nc}$ ($\cA$ equal to a $C^*$-algebra) into $\cL(\cY)_{\rm nc}$ ($\cY$ equal to a Hilbert space) can be reformulated
as problems concerning completely bounded maps from $\cA^{M \times M}$ to $\cL(\cY^M)$ (for some in general large $M$)
which are  $(\cS, \cS^*)$-bimodule maps, for a carefully specified subalgebra $\cS$ of finite complex matrices ${\mathbb C}^{M \times M}$.

Let us translate the reductions for the reducibility problem for cb nc kernels to reductions for the corresponding reducibility problem
for such cb $(\cS, \cS^*)$-bimodule maps.   In general we shall say that a map $\phi \colon \cA \to \cL(\cY)$ is {\em Hermitian} if
$\phi(P)^* = \phi(P^*)$ for all $P \in \cA$.    It is not hard to see that any cb map $\phi$ can be decomposed as a linear combination
$\phi = \phi_R + i \phi_I$ of two Hermitian maps $\phi_R = \frac{1}{2} (\phi + \phi^*) $ and $\phi_I = \frac{1}{2i}(\phi - \phi^*)$, i.e.,
$\phi_R = \phi_R^*$ and $\phi_I = \phi_I^*$ where in general
$$
   \phi^*(P) = \phi(P^*)^*.
$$
Thus the problem of {\em decomposability} for a  cb map $\phi$ (writing $\phi$ as a linear combination of four cp maps) reduces to decomposing a Hermitian map as the difference of two cp maps;  this is less trivial than the Hermitian decomposition just discussed
but has been done in work of Wittstock, Haagerup, and Paulsen (see the discussion in Example \ref{E:1} below).  
Our interest here is to understand these problems
for maps $\phi \colon \cA^{M \times M} \to \cL(\cY^M)$ which are also $(\cS, \cS^*)$-bimodule maps and where we wish to maintain
the $(\cS, \cS^*)$-bimodule structure in the components of the decomposition, a topic also well explored in the book of
\cite{Paulsen}.  The analogue of Proposition \ref{P:HerKer-criterion} is as follows.

\begin{proposition}  \label{P:HerMap-criterion}  Let $\cS$ be a subalgebra of ${\mathbb C}^{M \times M}$.  For $\cA$ a $C^*$-algebra
and $\cY$ a Hilbert space, we may then view $\cA^{M \times M}$ and $\cL(\cY^M) \cong \cL(\cY)^{M \times M}$ as 
$(\cS, \cS^*)$-bimodules.  Let $\phi$ be a cb Hermitian map from $\cA^{M \times M}$ to $\cL(\cY^M)$ which is also a $(\cS, \cS^*)$-bimodule
map.  Then the following are equivalent.

\smallskip

\textbf{\rm (1)} $\phi$ is $(\cS, \cS^*)$-decomposable, i.e., there exists  cp $(\cS, \cS^*)$-bimodule maps $\phi_1$ and $\phi_2$
so that $\phi = \phi_1 - \phi_2$.

\smallskip

\textbf{\rm (2)}  $\phi$ has a Hermitian $(\cS, \cS^*)$-bimodule Stinespring representation, i.e., there exists an operator
$V \colon \cX \to \cY^M$, a $*$-representation $\pi \colon \cA^M \to \cL(\cS)$ and a signature operator $J \in \cL(\cX)$
($J = J^* = J^{-1}$) so that 
$$
 \pi(P) J = J \pi(P), \quad  \phi(P) = V \big( \pi(P) J \big) V^*
$$
for all $P \in \cA^{M \times M}$ such that 
\begin{equation}   \label{extra}
 \pi(\beta^*) V^* - V^* L_\beta^* = 0 \text{ for all } \beta \in \cS.
\end{equation}

\smallskip

\textbf{\rm (3)}  There exists  cp $(\cS, \cS^*)$-bimodule maps $\Phi_{11}$, $\Phi_{22}$ so that the map
$\Phi$ from $\cA^{M \times M}$ to  $\cL(\cY^{2M})$ given by
\begin{equation}   \label{def-hatPhi}
 \widehat  \Phi \colon P \mapsto \begin{bmatrix} \Phi_{11}(P) & \phi(P) \\ \phi(P) & \Phi_{22}(P) \end{bmatrix}
\end{equation}
is cp.
\end{proposition}

\begin{proof}  When $\cS$ is the intertwining algebra associated with a finite nc point-set $\Omega = \{ Z_1, \dots, Z_N\}$
and the map $\phi$ has the form $\phi = \phi_K$ for a nc kernel on $\Omega$,  then the mutual equivalence of 
parts (1) of Propositions \ref{P:HerKer-criterion} and \ref{P:HerMap-criterion} and of parts (3) of Propositions
 \ref{P:HerKer-criterion} and \ref{P:HerMap-criterion}.  For the mutual equivalence of parts (2), one should also note
 the alternative proof of \eqref{Vmod} in the proof of Theorem \ref{T:bimodule-Stinespring}.  
 
 Alternatively, one can avoid the assumption that $\cS$ is an intertwining algebra
 associated with some point set $\Omega$ and prove directly the equivalence of (1), (2), (3) in Proposition \ref{P:HerMap-criterion}
 by paralleling the proofs of the corresponding results in Proposition \ref{P:HerKer-criterion}.  When doing this, when working
 with the Kolmogorov decompositions in part (2), one should bear in
 mind the result of Theorem \ref{T:bimodule-Stinespring} that condition \eqref{Vmod} is automatic in Stinespring representations
 for cp $(\cS, \cS^*)$-bimodule maps, while its counterpart \eqref{extra} is part of the definition for a $(\cS, \cS^*)$-Hermitian
 Stinespring representation.
 \end{proof}

Note that when we apply these criteria to the case where $M = N_0 = \sum_{j=1}^N n_j$ and $\cS$
is the intertwining algebra for the point $Z^{(0)} = \bigoplus_{i=1}^N Z_i$ where $\Omega = \{ Z_1, \dots, Z_N\}$,
then the content of Proposition \ref{P:HerMap-criterion} is just a direct translation of the content of \ref{P:HerKer-criterion}
via the dictionary between nc kernels and cb maps given by  Proposition \ref{P:ker/map}.

Thus the nc kernel-decomposability problem, or equivalently, the cb $(\cS, \cS^*)$-bimodule-map decomposability problem,
is reduced to showing that one of the criteria in parts (2) or (3) of Proposition \ref{P:HerKer-criterion} or in parts (2) or (3)
of Proposition \ref{P:HerMap-criterion} always holds.  The first of our partial results works with criterion (3) in
Proposition  \ref{P:HerMap-criterion} needs an extra hypothesis on the subalgebra $\cS \subset {\mathbb C}^{M \times M}$.
Let us use the notation  $C^*(\cS)$ for the $C^*$-algebra generated by $\cS$ inside ${\mathbb C}^{M \times M}$.

\smallskip

\noindent
\textbf{The Complete Spectral-Factorization Property:}  For any natural number $n \in {\mathbb N}$ and 
$\alpha \in C^*(\cS)^{n \times n}$
with $\alpha \succ 0$, there is a $\beta \in \cS^{n \times n}$ so that $\alpha = \beta \beta^*$.

\smallskip
With this additional hypothesis imposed, we have the following result.

\begin{theorem} \label{T:Decom1}
 Assume the same setup as in the hypotheses of Proposition \ref{P:HerMap-criterion} and assume that $\cS$ has the
 complete spectral-factorization property. 
Then statement (3) in Proposition \ref{P:HerMap-criterion} holds,
and the decomposability problem for Hermitian $(\cS, \cS^*)$-bimodule maps is solvable.  Hence also
the decomposability problem for Hermitian nc kernels on a finite set of nc points $\Omega$ is solvable, as long as
the intertwining algebra $\cS$ associated with $\Omega$ has the complete spectral-factorization property.
\end{theorem} 

For the case where $\cS$ is a $C^*$-algebra, Theorem \ref{T:Decom1} appears in Bhattacharyya-Dritschel-Todd
\cite{BDT} as Theorem 4.1 and is based on the off-diagonal method of Paulsen (see Exercise 8.6 in \cite{Paulsen}).
For completeness we include a complete proof since it also illustrates an application of Corollary \ref{C:bimodule-Stinespring}.

\begin{proof}[Proof of Theorem \ref{T:Decom1}]
By rescaling we may also assume that $\phi$ is cc.  We must construct cp maps $\Phi_{11}$ and $\Phi_{22}$
so that the map $\widehat \Phi$ given by \eqref{def-hatPhi} is cp.  To simplify the notation let us set
$$
{\mathbb A}: = \cA^{M \times M}, \quad {\mathbb L}:= \cL(\cY^M).
$$
For any $\alpha \in C^*(\cS) \subset {\mathbb C}^{M \times M}$, tensoring with the unit $1_\cA$ of $\cA$ gives an element
$\alpha \cdot 1_\cA$ in $\cA^{M \times M} = {\mathbb A}$ and similarly $\alpha \cdot I_\cY \in {\mathbb L}$.
Let us define an operator system ${\mathbb S}$ contained in the $C^*$-algebra ${\mathbb A}^{2 \times 2}$ by
$$
 {\mathbb S} = \begin{bmatrix} C^*(\cS) \cdot 1_\cA  & {\mathbb A} \\ {\mathbb A} & C^*(\cS) \cdot 1_\cA  \end{bmatrix}
$$
and define a map $\widehat \Phi_{\rm pre} \colon {\mathbb S} \to \cL(\cY)^{2M \times 2M}$ by
\begin{equation}   \label{hatPhidef}
  \widehat \Phi_{\rm pre} \colon \begin{bmatrix} \alpha \cdot 1_\cA & P_1 \\ P_2^* & \beta \cdot 1_\cA\end{bmatrix}
  \mapsto \begin{bmatrix} \alpha \cdot I_{\cY}  & \phi(P_1) \\ \phi(P_2^*) & \beta \cdot I_{\cY} \end{bmatrix}.
\end{equation}
We shall prove that this $\widehat \Phi_{\rm pre}$ is cp. To do this, for each $n \in {\mathbb Z}$ we must consider a positive element
in the inflated space ${\mathbb S}^{n \times n}$ and show that $\widehat \Phi^{(n)}$ sends such an element to a positive element in
the corresponding inflated space ${\mathbb L}^{2 \times 2}$.  Using a canonical shuffle procedure, we write elements of 
${\mathbb S}^n$
as
\begin{equation}  \label{element}
 \begin{bmatrix} H \cdot 1_{\cA} & {\mathbb P}_1 \\ {\mathbb P}_2^* & K \cdot 1_{\cA} \end{bmatrix}
\end{equation}
where $H, K \in \cS^{n \times n}$ and $P_1, P_2 \in {\mathbb A}^{n \times n}$, and  via the same canonical shuffle
we view elements of $({\mathbb L}^{2 \times 2})^{n \times n}$ as being elements of $({\mathbb L}^{n \times n})^{2 \times 2}$.
The goal now is to show that if
\eqref{element} is positive in ${\mathbb S}^{n \times n}$, then $\widehat \Phi_{\rm pre}$ is positive in ${\mathbb L}^{2 \times 2}$.
By an approximation argument, it is sufficient to assume that \eqref{element} is strictly positive definite.  It then follows
that $H \succ 0$ and $K \succ 0$ in $\cS^{n \times n}$.  By the assumption that $\cS$ has the complete spectral factorization
property,  we then have factorizations
$$
  H = A A^*, \quad K = B B^*
$$
where $A$ and $B$ are in $\cS^{n \times n}$.  As \eqref{element} is positive, it follows that ${\mathbb P}_1 = {\mathbb P}_2 =:
{\mathbb P}$ in ${\mathbb A}^{n \times n}$.
We then see that  \eqref{element} has a factorization as
$$
\begin{bmatrix} H & {\mathbb P} \\ {\mathbb P}^* & K \end{bmatrix} =
\begin{bmatrix} A & 0 \\ 0 & B   \end{bmatrix} \cdot
\begin{bmatrix}   I_{\cY^{n M}} & A^{-1}   {\mathbb P} B^{*-1} \\
  B^{-1} {\mathbb P}^*  A^{*-1}  &  I_{\cY^{n M}} \end{bmatrix}  \cdot \begin{bmatrix} A^*    & 0 \\ 0 & B^* \end{bmatrix}
$$
As $\phi$ is cc, we have also that
\begin{equation}   \label{phin-cc}
   \| \phi^{(n)} (A^{-1} P B^{-1}) \| < 1.
\end{equation}
Furthermore, as $\phi$ is a $(\cS, \cS^*)$-bimodule map, it follows that $\phi^{(n)}$ is a $(\cS^{n \times n}, (\cS^*)^{n \times n})$-bimodule
map and hence
\begin{equation}   \label{phi-bimodule}
  \phi^{(n)}\big( A^{-1} {\mathbb P} B^{* -1}  \big) = A^{-1} \cdot \phi^{(n)}({\mathbb P}) \cdot B^{*-1}.
\end{equation}
On the other hand let us compute
\begin{align*}
& \widehat \Phi^{(n)}_{\rm pre}  
\bigg(\begin{bmatrix} H \cdot 1_{\cA} & {\mathbb P} \\ {\mathbb P}^* & K \cdot 1_{\cA} \end{bmatrix} \bigg)
= \begin{bmatrix} H \cdot I_\cY  & \phi^{(n)}({\mathbb P})  \\ \phi^{(n)}({\mathbb P})^*  & K \cdot I_\cY \end{bmatrix}  \\
& \quad = \begin{bmatrix} A & 0 \\ 0 &  B \end{bmatrix} \cdot
\begin{bmatrix} I_{\cY^{nM}} & A^{-1} \phi^{(n)}({\mathbb P}) B^{* -1} \\  B^{*-1} \phi^{(n)}({\mathbb P}^* ) A^{*-1} & I_{\cY^{nM}}
\end{bmatrix} \cdot  \begin{bmatrix} A^* & 0 \\ 0 & B^* \end{bmatrix}  \\
& \quad   \begin{bmatrix} A & 0 \\ 0 &  B \end{bmatrix} \cdot
\begin{bmatrix} I_{\cY^{nM}} &  \phi^{(n)}(A^{-1}{\mathbb P} B^{* -1})
\\  \phi^{(n)}(B^{*-1} {\mathbb P}^* A^{*-1}) & I_{\cY^{nM}}
\end{bmatrix} \cdot  \begin{bmatrix} A^* & 0 \\ 0 & B^* \end{bmatrix}  
\end{align*}
where we make use of  \eqref{phi-bimodule} for the last step.  Making use of \eqref{phin-cc}, we see that the middle
factor in this last expression is positive and hence so also is $\widehat \Phi^{(n)}_{\rm pre}\left( \sbm{ H & {\mathbb P}
\\ {\mathbb P}^* & K }\right)$, and it follows that $\widehat \Phi_{\rm pre}$ is cp as wanted.

By the Arveson extension theorem (see e.g.~\cite[Theorem 7.5]{Paulsen}), it follows that there is a cp map
$\widehat \Phi_{\rm ext} \colon {\mathbb A}^{2 \times 2} \to {\mathbb L}^{2 \times 2}$ which 
extends $\widehat \Phi_{\rm pre}$:  $\widehat \Phi_{\rm ext}(X) =  \widehat \Phi_{\rm pre}(X)$ for $X \in {\mathbb S}$.
We claim next that necessarily $\widehat \Phi_{\rm ext}$ has the form
\begin{equation}   \label{hatPhi-ext-form}
\widehat \Phi_{\rm ext} \colon \begin{bmatrix}  P_{11} & P_{12}  \\ P_{21}^* & P_{22}
\end{bmatrix} \mapsto \begin{bmatrix} \Phi_{11}(P_{11})  & \phi(P_{12}) \\ \phi(P_{21}^*) & \Phi_{22}( P_{22}) \end{bmatrix}
\end{equation}
for cp maps $\Phi_{11}$ and $\Phi_{22}$ from ${\mathbb A}$ to ${\mathbb L}$.  To see this, let $\cS_0 \subset {\mathbb C}^{2 \times 2}$
be the subalgebra generated by the matrix $\sbm{ 1 & 0 \\ 0 & -1 }$.   This subalgebra is in fact a $C^*$-subalgebra of 
${\mathbb C}^{2 \times 2}$.  By conventions which we have already used several times, 
we can view both ${\mathbb A}^{2 \times 2}$ as well as  ${\mathbb L}^{2 \times 2}$ as $\cS_0$-bimodules.
It is easily checked that the operator system ${\mathbb S}$ is invariant under multiplication on the left as well as on the right
by elements of $\cS_0$, i.e., ${\mathbb S}$ is a $\cS$-sub-bimodule of ${\mathbb A}^{2 \times 2}$.  Furthermore it is a direct check to 
see that the map $\widehat \Phi_{\rm pre}$ given by \eqref{hatPhidef} is a $(\cS_0, \cS_0)$-bimodule map.  As a consequence of
Corollary \ref{C:bimodule-Stinespring}  it follows that $\widehat \Phi_{\rm ext}$ is also a $\cS_0$-bimodule map.
We claim that this $\widehat \Phi_{\rm ext}$ being a $\cS_0$-bimodule map forces $\widehat \Phi_{\rm ext}$ to have the 
form  \eqref{hatPhi-ext-form}.  Indeed, if $P_{11} \in {\mathbb A}^{2 \times 2}$ has the form 
$P = \sbm{P_{11} & 0 \\ 0 & 0 }$, then $P = \sbm{ 1 & 0 \\ 0 & 0 } \cdot P \cdot \sbm{ 1 & 0 \\ 0 & 0 }$ 
where $\sbm{ 1 & 0 \\ 0 & 0 } = \frac{1}{2}\big( \sbm{ 1 & 0 \\ 0 & 1} + \sbm{ 1 & 0 \\ 0 & -1} \big) \in \cS_0$.   Hence 
$$
\widehat \Phi_{\rm ext}(P) = \widehat \Phi_{\rm ext}\bigg( \sbm{ 1 & 0 \\ 0 & 0 } P \sbm{ 1 & 0 \\ 0 & 0 } \bigg)
 = \sbm{ 1 & 0 \\ 0  & 0 } \widehat \Phi_{\rm ext}(P) \sbm{ 1 & 0 \\ 0 & 0 }  = : \begin{bmatrix}  \Phi_{11}(P_{11}) & 0 \\ 0 & 0 \end{bmatrix}.
$$
for a cp positive map $\Phi_{11}$.
Similarly, if $P = \sbm{ 0 & 0 \\ 0 & P_{22}}$, then there is a cp map $\Phi_{22}$ so that
$\widehat \Phi_{\rm ext}(P) = \sbm{0 & 0 \\ 0 & \Phi_{22}(P_{22})}$.
Finally, if $P \in {\mathbb A}^{2 \times 2}$ has the form $P = \sbm{ 0 & P_{12} \\ 0 & 0 }$,
then $P \in {\mathbb S}$ and hence
$$
\widehat \Phi_{\rm ext}(P) = \widehat \Phi_{\rm pre}(P) = \begin{bmatrix} 0 & \phi(P_{12}) \\ 0 & 0 \end{bmatrix}.
$$
By a similar argument we also have
$$
  \widehat \Phi_{\rm ext}\bigg( \sbm{ 0 & 0 \\ P_{21}^* & 0 } \bigg) = \begin{bmatrix} 0 & 0 \\ \phi(P_{21}^*) & 0 \end{bmatrix}.
$$
By linearity it now follows that $\widehat \Phi_{\rm ext}$ has the form \eqref{hatPhi-ext-form} as claimed.

Let now $\iota \colon {\mathbb A} \to {\mathbb A}^{2 \times 2}$ be the cp map
$$
 \iota \colon P \mapsto \begin{bmatrix} P & P  \\ P & P \end{bmatrix}
$$
and finally let $\widehat \Phi$ be the composition
$$
\widehat \Phi = \widehat \Phi_{\rm ext} \circ \iota \colon {\mathbb A} \to {\mathbb A}^{2 \times 2}.
$$
Then $\widehat \Phi$ is a composition of cp maps and hence is cp itself, and has the form
exactly the form \eqref{def-hatPhi} demanded in part (3) of Proposition \ref{P:HerMap-criterion}.
\end{proof}

There are interesting examples of operator algebras with the complete spectral-factorization property which are not already
$C^*$-algebras themselves, especially in
function theory (e.g., the disk algebra or the Hardy algebra $H^\infty$ over the unit disk (see e.g.~\cite{RR}).  
For the case of finite-dimensional
subalgebras of ${\mathbb C}^{M \times M}$, it is well known that the algebra of upper or lower triangular matrices has the
spectral factorization property (by the $LU$-factorization algorithm), a simple dimension argument shows that
one cannot factor a positive definite block matrix as a product of a block upper triangular and its adjoint with the
block upper triangular matrix having block entries which are again block upper-triangular.  Hence the applicability of the assumption of
the complete spectral-factorization property appears to be limited to subalgebras $\cS$ of ${\mathbb C}^{M \times M}$
which are in fact $C^*$-algebras, in which case the complete spectral-factorization property follows easily from standard
$C^*$-algebra theory.   As in practice there do not appear to be any finite-dimensional operator algebras with
the complete spectral factorization property which are not already $C^*$-algebras, it would appear that the $C^*$-algebra
case is essentially the general case.   For simplicity of terminology, let us make the following definition.

\begin{definition}  \label{D:admissible}  Let $\Omega = \{Z_i \colon 1 \le i \le N\}$ be a finite subset of
$\cV_{\rm nc}$ with $Z_i \in \Omega_{n_i}$.  Set $N_0 = \sum_{i=1}^N n_i$ so $Z^{(0)} = \bigoplus_{i=1}^N Z_i \in
\cV^{N_0 \times N_0}$.  Associate with $\Omega$ the subalgebra $\cS$ of ${\mathbb C}^{N_0 \times N_0}$ given by 
\eqref{algebraS}.  We say that the set of points $\Omega$ is {\em admissible} if $\cS$ is a $C^*$-algebra.
\end{definition}

With this terminology in hand we can state the following immediate corollary of Theorem \ref{T:Decom1}.

\begin{corollary}   \label{C:Decom1}  Suppose that
$$
K \colon \Omega \times \Omega \to \cL(\cA_{\rm nc}, \cL(\cY)_{\rm nc})
$$
is a free nc kernel on an admissible  finite set of points $\Omega \subset \cV_{\rm nc}$.  Then $K$ is decomposable.
\end{corollary}

\smallskip

We now explore some examples where the admissibility hypothesis does hold and thus Problem A is guaranteed to be solvable.

\begin{example} \label{E:1} \textbf{Hermitian decomposition for completely bounded maps.}
Consider the special case where $M=1$ in Theorem \ref{T:Decom1}.  Then of course $\cS = {\mathbb C}
= C^*(\cS)$ is admissible.  Then the result of Corollary \ref{C:Decom1} combined
with the equivalence of (1) $\Leftrightarrow$ (3) in Proposition \ref{P:HerMap-criterion} gives us 
a version of Wittstock's decomposition theorem \cite[Theorem 8.5]{Paulsen}.

\smallskip

\noindent
\textbf{Corollary.}  {\sl  Any cb map $\phi$ from a $C^*$-algebra $\cA$ to the $C^*$-algebra $\cL(\cY)$ ($\cY$ equal to
some Hilbert space) is decomposable.}

\smallskip

Let us also mention that our proof that for this case  $\phi$ being cc implies that $\widehat \Phi_{\rm pre}$ \eqref{hatPhidef}
 is cp appears as Lemma 8.1 in \cite{Paulsen} with the same proof.

For the case where $M=1$ and $\cS = {\mathbb C}$, of course once item (3) in Proposition \ref{P:HerMap-criterion} is
known, then it follows that items (1) and (2) also hold in general.   Item (3) is known as the Wittstock extension theorem
from \cite{Wittstock},  item (1) is known as the Wittstock decomposition theorem
(see \cite{Wittstock81}) and item (2) is a Hermitian version of the generalized Stinespring representation for cb maps
(see Theorem 8.4 in \cite{Paulsen} and the Notes there at the end of the chapter).  All these results are developed
in Chapter 8 of Paulsen's book \cite{Paulsen} but with a somewhat different organization.

It is of interest to specialize all this discussion to the case where $\cA = C(X)$ (continuous functions on a compact Hausdorff space $X$).
In case $\cY = {\mathbb C}$, then a map $\phi \colon \cA \to \cL(\cY)$ is given by a complex measure $\mu$.
 In this case $\mu$ being of finite total variation 
is equivalent to $\phi$ being completely bounded, and the conclusion of Corollary \ref{C:Decom1} for this case 
amounts to the classical Jordan decomposition for a complex measure.  In case we still take $\cA = C(X)$ but we take $\cY$
to be a general Hilbert space, then a linear map $\phi$ from $\cA$ to $\cL(\cY)$ corresponds to an operator-valued measures $\mu$,
but $\mu$ being of finite total variation does not always match up with $\phi$ being completely bounded.  In any case it is the
complete boundedness hypothesis which guarantees a Jordan decomposition for the operator-valued measure $\mu$
(see \cite[pages 104-106]{Paulsen} for a fuller discussion).
\end{example}

\begin{example} \label{E:2}  \textbf{$\Omega$ with block-diagonal selfadjoint intertwining space}
 Let us next consider a finite point set $\Omega = \{ Z_1, \dots, Z_d \}$ (with $Z_i \in \Omega_{n_i}$) contained in the nc universal space
$\cV_{\rm nc}$ such that 
\begin{equation}   \label{trivial-intertwine}
\cI(Z_i, Z_j)  = \begin{cases}  0 & \text{if } i \ne j, \\
\text{a $C^*$-subalgebra $\cD_i$ of ${\mathbb C}^{n_i \times n_j}$}  &\text{if } i = j.
\end{cases}
\end{equation}
where we use the notation 
$$
\cI(Z_i, Z_j) = \{ \alpha \in {\mathbb C}^{n_i \times n_j} \colon Z_i \alpha = \alpha Z_j\}.
$$
Roughly, for $i \ne j$ there is no piece of $Z_i$ similar to a piece of $Z_j$, and for $i=j$, each $Z_i$ has a nice commutant over 
${\mathbb C}^{n_i \times n_i}$. 
For this {\em block-diagonal selfadjoint intertwining space} case, 
the intertwining algebra $\cS = \cI( Z^{(0)})$ is equal to  $\bigoplus_{i=1}^N \cD_i$
and hence $\Omega$ is admissible, and Corollary \ref{C:Decom1} leads to the following result.
\smallskip

\noindent
\textbf{Corollary.}  {\sl Suppose that $K \colon \Omega \times \Omega \to \cL(\cA_{\rm nc}, \cL(\cY)_{\rm nc})$ is a cb nc kernel
and $\Omega$ has the {\em block-diagonal selfadjoint intertwining space} property \eqref{trivial-intertwine}.  Then  $K$ is decomposable.}

\smallskip

A particular instance of this scenario is the case where $\Omega$ consists exclusively of scalar points which are all distinct:
$$
  \Omega = \{ z_1, \dots, z_N \colon z_i \ne z_j \text{ for } i \ne j \} = \Omega_1 \subset \cV_{\rm nc}.
$$
In this case $\cS$ is the $C^*$-algebra consisting of the diagonal matrices
$$
  \cS =  \cD_N:=\{ \operatorname{diag}_{1 \le i \le N_0 = N} [\lambda_i] \text{ with } \lambda_i \in {\mathbb C} \text{ for all } i \}.
$$
A cp nc kernel on a set of distinct scalar points can be viewed as a cp kernel in the sense of Barreto-Bhat-Liebscher-Skeide
\cite{BBLS}.     By the preceding
Corollary we conclude that {\em any BBLS-kernel from $\Omega$ to $\cL(\cA, \cL(\cY))$ is decomposable.}  This essentially
recovers the main result of \cite{BDT} for the finite-point case; there the connection with $\cD_N$-bimodule maps
$\phi \colon \cA^{n \times n} \to  \cL(\cY)^{N \times N}$ plays a prominent role, including the characterization of $\cD_N$-bimodule
maps as entry-wise maps as was used here for the $N=2$ case in the proof of \eqref{hatPhi-ext-form}, but with the
connection with nc kernels not made explicit.

Let us specialize still further by taking $\cA = {\mathbb C}$.  Then we may identify $\phi$ with its value at $1 \in {\mathbb C}$;
thus $\phi$ amount to an operator $\phi(1) = T \in \cL(\cY)$.  Boundedness here is the same as complete boundedness
and Corollary \ref{C:Decom1} gives us the decomposition $T = T_1 - T_2 + i (T_3 - T_4)$ where $T_j$ are positive operators on
$\cH$.
\end{example}

\begin{remark}  \label{R:generic}   For this discussion let us write the set of points of $\Omega$ with superscripts
$$
  \Omega = \{ Z^{(1)}, \dots, Z^{(N)} \} \text{ with } Z^{(i)} \in \Omega_{n_i}.
$$
As $\Omega$ consists of only finitely many points, the linear span of all the coordinates of points in $\Omega$
$$
  \operatorname{span} \{ z^{(k)}_{ij} \colon 1 \le i,j \le n_k, \, 1 \le k \le N \} \text{ where } Z^{(k)} = [ z^{(k)}_{ij}]_{1 \le i,j \le n_k}
  \in \cV^{n_k \times n_k}.
$$
spans a finite-dimensional subspace $\cM$ of $\cV$.  As long as we are working with this finite subset of $\cV_{\rm nc}$,
without loss of generality we may assume that $\cV  = {\mathbb C}^d$ for a sufficiently large $d$.  
By choosing a basis for $\cM$ we may identify $\cM   = \cV = {\mathbb C}^d$ (with $d = \operatorname{dim} \cM$).  
Then we can analyze in terms of this representation
for $\Omega$ when it is the case that $\Omega$ is admissible.  The result is as follows.

\smallskip

\noindent
\textbf{Proposition.} (1) {\sl  If $d=1$, then $\Omega$ is admissible if and only if the matrix
$Z^{(0)} =  Z^{(1)} \bigoplus \cdots \bigoplus Z^{(N)}$ is normal, i.e., there is an orthonormal basis
which diagonalizes $Z^{(0)}$.}

\smallskip

\noindent
(2) {\sl For $d>1$, the case where $\Omega$ is admissible is generic, i.e., the set of $Z^{(0)} \in ({\mathbb C}^{N_0 \times N_0})^d$
such that $\Omega$ is admissible is dense in the space of all $Z^{(0)} \in ({\mathbb C}^{N_0 \times N_0})^d$.}

\smallskip

\noindent
\textbf{Proof.}  Suppose that $d=1$ and $Z^{(0)}$ is normal.  Thus $Z^{(0)} = \sum_{i=1}^M \lambda_i P_i$
where $\lambda_1, \dots, \lambda_M$ are the eigenvalues of $Z^{(0)}$ where where $P_1, \dots, P_M$
is a pairwise-orthogonal family of orthogonal projections summing to the identity.  
Then for $p$ a polynomial, we have
$p(Z^{(0)})  = \sum_{i=1}^M p(\lambda_i) P_i$.  In particular, by solving a Lagrange interpolation problem we can always find
a polynomial $p_i$ with $p_i(\lambda_i) = 1$ and $p_i(\lambda_j) = 0$ for $j \ne i$; with this choice of polynomial we have
$p_i(Z^{(0)}) = P_i$.  We conclude that if $\alpha \in {\mathbb C}^{N_0 \times N_0}$ commutes with $Z^{(0)}$
then $\alpha$ commutes with each $P_i$ and hence must have the form 
\begin{equation}   \label{alpha-form}
\alpha = \sum_{i=1}^M \alpha_i P_i
\end{equation}
where $\alpha_i$ is an operator on $\operatorname{Ran} P_i$.  Conversely, any operator of this form commutes with
$Z^{(0)}$, and hence the form \eqref{alpha-form} characterizes the commutant of $Z^{(0)}$.  If $\alpha$ is of the form
\eqref{alpha-form}, then $\alpha^* = \sum_{i=1}^M \alpha_i^* P_i$ is also of this form, so $ \cS = \cI(Z^{(0)})$ is selfadjoint,
i.e., $\Omega$ is admissible in this case.

Suppose next that $Z^{(0)}$ is not normal.   Then certainly $Z^{(0)}$ commutes with itself, but on the other hand $Z^{(0)*}$
does not commute with $Z^{(0)}$ (i.e., $Z^{(0)^*} Z^{(0)} \ne Z^{(0)} Z^{(0)*}$) since $Z^{(0)}$ is not normal.
Hence the commutant of $\{Z^{(0)}\}$ is not a $C^*$-algebra and $\Omega$ is not admissible. 

Next suppose that $d > 1$ and $Z^{(0)} = ( Z^{(0)}_1, \dots, Z^{(0)}_d\}$ where each $Z^{(0)}_j \in {\mathbb C}^{N_0 \times N_0}$.
In this case it is known that generically the algebra generated by the matrices $Z^{(0)}_1, \dots, Z^{(0)}_d$ in 
${\mathbb C}^{N_0 \times N_0}$ is all of ${\mathbb C}^{N_0 \times N_)}$ (see \cite{Pascoe} for a recent treatment).  
If we assume that we are in this
generic case, then, given any $W \in {\mathbb C}^{N_0 \times N_0}$ there is a nc polynomial $p$ in $d$ arguments so that
$p(Z^{(0)}_1, \dots, Z^{(0)}_d) = W$.  Hence, if $\alpha \in {\mathbb C}^{N_0 \times N_0}$ is in the intertwining space
$\cI(Z^{(0)})$, then 
$$
 \alpha W = \alpha p(Z^{(0)}) = p(Z^{(0)}) \alpha = W \alpha
$$
i.e., $\alpha $ commutes with all of ${\mathbb C}^{N_0 \times N_0}$.  This forces $\alpha$ to be a scalar multiple of the identity
matrix.  Hence $\cI(Z^{(0)})$ consists of scalar matrices, and in particular is selfadjoint.  We conclude that for this case
a generic set of points $\Omega$ is admissible as claimed.   \qed
\end{remark}

We shall discuss possible extensions of the kernel version of Theorem \ref{T:Decom1} to infinite point sets $\Omega$
in Section \ref{S:KuroshDecom}.

\smallskip

Our next partial result concerning the decomposability problem for nc Hermitian kernels/Hermitian maps is based on
criterion (2) in Propositions \ref{P:HerKer-criterion} and \ref{P:HerMap-criterion}.  We shall state only for the setting
of Hermitian maps.

\begin{theorem}  \label{T:Decom2}

Suppose that  $\cS$ is a subalgebra of ${\mathbb C}^{M \times M}$ and
that $\phi$ is a Hermitian map from $\cA^{M \times M}$ to $\cL(\cY)^{M \times M}$.
\smallskip

\noindent
\textbf{\rm (1)} If $\phi$ is decomposable $(\cS, \cS^*)$-bimodule map which can be $(\cS, \cS^*)$-decom\-posable (i.e.,
$\phi$ is the difference of two cp $(\cS, \cS^*)$-bimodule maps), then $\phi$ has a Hermitian Stinespring representation
$\phi(P) = V \pi(P) J V^*$  satisfying the extra condition \eqref{extra} as in part (2) of Proposition \ref{P:HerMap-criterion}.

\smallskip

\noindent
\textbf{\rm (2)}  If $\phi$ is a $(\cS, \cS^*)$-bimodule map, then any Hermitian Stinespring representation
$\phi(P) = V \pi(P) J V^*$ as in part (2) of Proposition \ref{P:HerMap-criterion} (taken with $\cS = {\mathbb C}$ so that
condition \eqref{extra} can be ignored) automatically also satisfies \eqref{extra} but only in the weaker form
\begin{equation}  \label{extra-weak}
\operatorname{span}_{\beta \in \cS} \operatorname{Ran} (\pi(\beta^*) V^* - V L_{\beta^*}) \text{ is a $J$-isotropic subspace
of $\cX$.}
\end{equation} 
\end{theorem}

\begin{proof}  Item (1) is just a restatement of part (2) of Proposition \ref{P:HerMap-criterion}.

As for item (2), we compute, for $y,y' \in \cY^M$ and $\beta, \beta' \in \cS$,  
\begin{align*}
& \langle J ( \pi(\beta^*) V^* - V^* L_{\beta^*}) y, (\pi(\beta'^*) V^* - V^* L_{\beta'^*}) y' \rangle \\
& \quad = \bigg\langle y, \bigg(  V \pi(\beta \beta'^*)V^* - L_\beta V \pi(\beta'^*) V^* - V \pi(\beta) V^*L_{\beta'^*} 
   + L_\beta V^* V L_{\beta'^*} \bigg) y' \bigg\rangle  \\
 & \quad =  \bigg\langle y, \bigg(  \phi(\beta \beta'^*) - L_\beta \phi (\beta'^*)  - \phi(\beta)L_{\beta'^*} 
   + L_\beta \phi(1_{\cA^{M \times M}}) L_{\beta'^*} \bigg) y' \bigg\rangle   \\
    & \quad = \langle y, (\phi(\beta \beta'^*) - \phi(\beta \beta'^*) - \phi(\beta \beta'^*) + \phi(\beta \beta'^*)  ) y' \rangle = 0
 \end{align*}
 where we used that $\phi$ is a $(\cS, \cS^*)$-bimodule map in the penultimate step.
This holding for all $y, y' \in \cY^M$ and $\beta, \beta' \in \cS$ is just the statement that the subspace described
in \eqref{extra-weak} is $J$-isotropic.
\end{proof}

\section{Problem B: The Arveson extension problem for noncommutative positive kernels}  \label{S:Arv-ext}
The goal of this section is to solve affirmatively the Arveson extension problem for cp nc kernels in place of a cp map on a finite set 
$\Omega$.  The case where $\Omega$ is allowed to be infinite is considered in Section \ref{S:KuroshArvExt}.
This proof amounts to combining Corollary \ref{C:bimodule-Stinespring} with the Arveson extension theorem for
completely positive maps.

\begin{theorem} \label{T:Arv-ext-ker}
Suppose that $K$ is a nc kernel on a finite set of nc points $\Omega \subset \cV_{\rm nc}$
 with values in $\cL({\mathbb S}_{\rm nc}, \cL(\cY)_{\rm nc})$, where ${\mathbb S}$ is a operator system and $\cL(\cY)$
 is the $C^*$-algebra of operators on the Hilbert space $\cY$.  Let $\cA$ be a $C^*$-algebra containing ${\mathbb S}$.
 Then there exists a cp nc kernel $\widehat K$ 
 $$
  \widehat K \colon \Omega \times \Omega \to \cL(\cA_{\rm nc}, \cL(\cY)_{\rm nc})
 $$
 such that $\widehat K(Z,W)(P) = K(Z,W)(P)$ for all $Z \in \Omega_n$, $W \in \Omega_m$ whenever $P \in {\mathbb S}^{n \times m}$.
 \end{theorem}

\begin{proof}  Let us assume that the finite set $\Omega \subset \cV_{\rm nc}$ is given by
$$
\Omega = \{ Z_1, \dots, Z_N \} \subset \cV_{\rm nc}.
$$ 
Set $Z^{(0)} = \operatorname{diag}_{1 \le i \le N} [ Z_i ]$, set $N_0 = \sum_{i=1}^N n_i$,  let
$\cS$ be the intertwining algebra of $Z^{(0)}$
$$
\cS = \{ \alpha \in {\mathbb C}^{N_0 \times N_0} \colon \alpha Z^{(0)} = Z^{(0)} \alpha \},
$$
and let $\phi_K$ be the cp $(\cS, \cS^*)$-bimodule map from ${\mathbb S}^{N_0 \times N_0}$ to $\cL(\cY^{N_0})$
associated with $K$ as in \eqref{ker/map}.  By the Arveson extension theorem for completely positive maps 
(see \cite[Theorem 7.5]{Paulsen}), there is an extension of $\phi_K$
to a cp map $\widehat \phi \colon \cA^{N_0 \times N_0} \to \cL(\cY^{N_0})$.  By Corollary  \ref{C:bimodule-Stinespring}, 
any such $\widehat \phi$ is also a $(\cS, \cS^*)$-bimodule map.  Then by Proposition \ref{P:ker/map} applied in the
reverse direction, we conclude that $\widehat \phi = \phi_{\widehat K}$ for a cp nc kernel $\widehat K$ on $\Omega$
with values equal to operators from $\cA_{\rm nc}$ to $\cL(\cY)_{\rm nc}$.  As $\phi_{\widehat K} = \widehat \phi$
is an extension of $\phi_K$, we conclude that $\widehat K$ is an extension of $K$ as wanted.
\end{proof}

\section[Kernel dominance]{Problem C: Kernel-dominance certificates}  \label{S:ker-dom} 
In this section we consider the finite-point case of Problem C for the case where $\Omega \subset \cV_{\rm nc}$ is finite.
Given a point $Z$ in some nc set $\Omega$,  we shall use the notation $n_Z$ for the natural number $n$ such that 
$Z$ has matrix size $n \times n$ (i.e., $Z \in \Omega_n$).

\begin{theorem}  \label{T:ker-dom}
We let $\cE$ and $\cN$ be Hilbert spaces.
 We are given a full nc subset
$\Xi$ of $\cV_{\rm nc}$ together with a  Hermitian nc kernel
$$
\fQ \colon \Xi \times \Xi \to \cL({\mathbb C}_{\rm nc}, \cL(\cN)_{\rm nc}).
$$
Let ${\mathbb P}_\fQ$ be the strict positivity domain for $\fQ$ as defined by
$$
  {\mathbb P}_\fQ = \amalg_{n=1}^\infty \{ Z \in \Xi_n \colon \fQ(Z, Z)(I_n) \succ 0 \}.
$$
Suppose that:

\smallskip

\noindent
\textbf{ \rm (i)} $\Omega$ is a subset of ${\mathbb P}_\fQ$.

\smallskip

\noindent
\textbf{\rm (ii)}
$\fS$ is a Hermitian nc kernel defined on  $\Omega' : = [\Omega]_{\rm full} \cap {\mathbb P}_\fQ$ 
$$
\fS \colon \Omega' \times \Omega' \to \cL({\mathbb C}_{\rm nc}, \cL(\cE)_{\rm nc})
$$
such that
\begin{equation}   \label{ker-dom-hy}
   \fS(Z, Z)(1_{{\mathbb C}^{n_Z \times n_Z}}) \succeq 0 \text{ for all } Z \in \Omega'.
 \end{equation}
 
 \smallskip
 
 \noindent
 \textbf{\rm (iii)} The restriction of $\fQ$ to a finite subset $\Xi_{\rm finite}$ of $\Xi$ is decomposable.
  
 \smallskip
 
 \noindent
 \text{\rm (iv)}  The restriction of $\fS$ to a finite subset $\Omega'_{\rm finite}$ of $\Omega'$ is decomposable.
  
 \smallskip
  
 \noindent
 Then there exist two completely positive nc kernels on $\Omega$
 \begin{equation}   \label{ker-pair}
\Gamma_1 \colon \Omega \times \Omega  \to \cL(\cL(\cN)_{\rm nc}, \cL(\cE)_{\rm nc}),  \quad 
 \Gamma_2 \colon \Omega \times \Omega \to \cL({\mathbb C}_{\rm nc}, \cL(\cE)_{\rm nc})
\end{equation}
so that, for all $Z \in \Omega_n$, $W \in \Omega_m$, $P \in {\mathbb C}^{n \times m}$ we have
\begin{equation}   \label{cone-rep}
\fS(Z,W)(P) = \Gamma_1(Z,W)(\fQ(Z,W)(P)) + \Gamma_2(Z, W)(P).
\end{equation}
\end{theorem}

\begin{proof}
In this section we consider only the case where two additional hypotheses are in force:
\begin{enumerate}
\item[(H1)] The set $\Omega \subset {\mathbb P}_\fQ$ is finite.
\item[(H2)] The coefficient Hilbert space $\cE$ is finite-dimensional.
\end{enumerate}
In applying conditions (ii) and (iv) in the statement of the theorem, we shall write simply $\Omega$ rather than $\Omega_{\rm fin}$.
The reduction of the general case to this special case will be discussed in Section \ref{S:KuroshDom}.

\smallskip

We are assuming that $\Xi$ is a subset of the nc set $\cV_{\rm nc}$ generated by the linear space $\cV$.   
As $\Omega$ is  finite by (H1), the span of all the matrix entries of all the finitely many points in $\Omega$ is contained
 in some finite-dimensional subspace $\cV_0$  of $\cV$.  By choosing some basis for $\cV_{0}$,
we may identify  $\cV_0$ with ${\mathbb C}^{d}$ ($d$ sufficiently large) and thereby
view elements of $\Omega$ as points in $({\mathbb C}^{d})_{\rm nc}$ as well.  Thus we may assume without loss of
generality that $\cV = {\mathbb C}^d$.  To make use of this reduction it is sometimes convenient to spell out
the components of a point $Z \in \Omega$ which exhibits its membership in ${\mathbb C}^d$:  using that
$({\mathbb C}^d)^{n \times n} \cong ({\mathbb C}^{n \times n})^d$, we shall write a typical point of $\Omega$ as
\begin{equation}  \label{Cdnc}
   Z = (Z_1, \dots, Z_d)  \text{ where  each } Z_i \in  {\mathbb C}^{n \times n} \text{ if } Z \in \Omega_{n}.
\end{equation}
Then the action of a pair of complex matrices $\alpha$,  $\beta$ of compatible sizes on a point $Z = (Z_1, \dots, Z_d)$ 
is component-wise:
$$
 \alpha \cdot (Z_1, \dots, Z_d) \cdot \beta = (\alpha Z_1 \beta, \dots, \alpha Z_d \beta)
$$
where the multiplication on the right side of the equation is ordinary matrix multiplication.
We shall use superscripts to distinguish various points of $\Omega$; thus we shall list the points of $\Omega$ as
\begin{equation}   \label{Omega-list}
\Omega = \{ Z^{(1)}, \dots, Z^{(N)} \} 
\end{equation}
 where $Z^{(i)}$ has the form $Z^{(i)} =  ( Z^{(i)}_1, \dots,  Z^{(i)}_d)$ with each component $Z^{(i)}_j$ in
${\mathbb C}^{n_i \times n_i}$ if $Z \in \Omega_{n_i}$.

We let $\fX$ be the linear space of all graded kernels $K$ on
$\Omega$   with values in $\cL({\mathbb C}, \cL(\cE))_{\rm nc}$; 
thus $K \in \fX$ means that $K$ is an operator-valued 
function on $\Omega  \times \Omega$ such that 
\begin{equation}  \label{mathfrakX} 
 K(Z,W) \in \cL({\mathbb C}^{n \times m}, \cL(\cE)^{n \times m})
 \text{ if } Z \in \Omega_{n} \text{ and } W 
\in \Omega_{m}.
\end{equation}
We make $\fX$ a Banach space by endowing
$\fX$ with the supremum norm:
$$
   \| K \|_{\fX} = \max \{ \| K(Z,W) \| \colon Z, W \in \Omega_0\}.
$$
As $\Omega$ is finite and $\cE$ is finite-dimensional by assumption, this Banach space is finite-dimensional.
Hence bounded subsets of $\fX$ are pre-compact in the norm topology.
Thus bounded sequences always have convergent subsequences.
 Furthermore, convergence of a sequence $\{K_k\}_{k \in {\mathbb N}}$ to $K$ in the norm topology of $\fX$ 
concretely just means pointwise that $K_k(Z,W)(P)$ converges to $K(Z,W)(P)$ in the norm topology
of $\cL(\cE^m, \cE^n)$ for each $Z \in \Omega_{n}$, $W \in \Omega_{m}$, and $P \in
{\mathbb C}^{n \times m}$ for all $n,m = 1,2, \dots$. 

We define a subset $\cC$ of $\cX$  to consist of all graded kernels 
$K$ in $\cX$ which have the form
\begin{equation} \label{cone-form}
  K(Z,W)(P) = \Gamma_{1}(Z,W)(\fQ(Z,W)(P)) + 
  \Gamma_{2}(Z,W)(P) 
\end{equation}
for all $Z \in \Omega_{n}$, $W \in \Omega_{m}$, $P \in {\mathbb C}^{n \times 
m}$ and $n,m \in {\mathbb Z}_{+}$ for some cp nc kernels $\Gamma_{1}$ and $\Gamma_2$
on $\Omega$:
$$
\Gamma_1 \colon \Omega \times \Omega \to \cL(\cL(\cN)_{\rm nc}, \cL(\cE)_{\rm nc}), \quad
\Gamma_2 \colon \Omega \times \Omega \to \cL({\mathbb C}_{\rm nc}, \cL(\cE)_{\rm nc}).
$$
Then we have

\begin{lemma}  \label{L:cone}
The subset $\cC$ is a  closed pointed cone in ${\fX}$.
\end{lemma}

\begin{proof} It is easily verified that $\cC$ is a cone. To see that $\cC$ is pointed, we must check that if $K \in \cC$ and $-K \in \cC$ then $K=0$. We see that such a $K$ must have the form
\begin{align*}
K(Z,W)(P)&=\Gamma_1(Z,W) (\fQ(Z,W)(P)) + \Gamma_2(Z,W)(P) \\
&=-\Gamma_3(Z,W) (\fQ(Z,W)(P)) - \Gamma_4(Z,W)(P).
\end{align*}
where $\Gamma_i$ is a cp nc kernel for $i=1,2,3,4$. We know that $\fQ(Z,Z)(I) \succ 0$ for $Z \in \Omega$ which results in $0 \preceq K(Z,Z)(I) \preceq 0$. We conclude that $K(Z,Z)(I)=0$ which results in $\Gamma_2(Z,Z)(I)= \Gamma_4(Z,Z)(I)=0$ as well. By Lemma \ref{L:zeroker}, the cp nc kernels $\Gamma_2$ and $\Gamma_4$ must be zero which means $K$ has the simpler form
\begin{align*}
K(Z,W)(P)=\Gamma_1(Z,W) (\fQ(Z,W)(P)) 
=-\Gamma_3(Z,W) (\fQ(Z,W)(P)).
\end{align*}
As $\Omega$ is a finite subset of ${\mathbb P}_\fQ$, we may choose $\epsilon_0$ with $0 < \epsilon_0 < 1$ so that 
$$
 \fQ(Z,Z)(I_n) \succeq \epsilon_0 I_{\cG^n}  \text{ for all } Z \in \Omega_{n} \text{ (n=1,2,\dots)}.
 $$
Hence, we have
\begin{align*}
	0 = K(Z,Z)(I)& = \Gamma_1 (Z,Z)(Q(Z,Z)(I) \succeq \epsilon_0 \Gamma_1(Z,Z)(I)\\
	0 = K(Z,Z)(I)& = - \Gamma_3 (Z,Z)(Q(Z,Z)(I) \preceq -\epsilon_0 \Gamma_3(Z,Z)(I)
\end{align*}
and it follows that $\Gamma_1 (Z,Z) (I) = \Gamma_3 (Z,Z)(I) = 0$. We apply Lemma \ref{L:zeroker} to conclude that 
$\Gamma_1=\Gamma_2=0$ which gives the desired result, $K=0$.

We now show that the cone is closed. Suppose that $\{K_k \colon k \in {\mathbb N} \}$ is a sequence of elements of $\cC$ such 
that $\| K - K_k \|_\fX \to 0$ as $k \to \infty$ for some $K \in \fX$.  The problem is to show that $K$ is again in $\cC$.

By definition for each $k$ there are  cp nc kernels $\Gamma_{1,k}$
from $\Omega \times \Omega$ to $\cL(\cL(\cN), \cL(\cE)))_{\rm nc}$ and $\Gamma_{2,k}$ from 
$\Omega \times \Omega$ to $\cL({\mathbb C}, \cL(\cE))_{\rm nc}$ so that 
$$
K_k(Z,W)(P) = \Gamma_{1,k}(Z,W)\left( \fQ(Z,W)(P)\right) + \Gamma_{2,k}(Z,W)(P)
$$
for all $Z \in \Omega_{n}$, $W \in \Omega_{m}$, $P \in {\mathbb C}^{n \times m}$ for all $m, n \in {\mathbb N}$.
As $\Omega$ is a finite subset of ${\mathbb P}_\fQ$, we may choose $\epsilon_0$ with $0 < \epsilon_0 < 1$ so that 
$$
 \fQ(Z,Z)(I_n) \succeq \epsilon_0 I_{\cN^n}  \text{ for all } Z \in \Omega_{n} \text{ (n=1,2,\dots)}.
 $$
 Hence
 \begin{align*}
 K_k(Z,Z)(I_n)  & = \Gamma_{1,k}(Z,Z)(\fQ(Z,Z)(I_{\cG^n}) ) + \Gamma_{2,k}(Z,Z)(I_n)\\
 & \succeq \epsilon_0  \Gamma_{1,k}(Z,Z) (I_{\cN^n}) + \Gamma_{2,K}(Z,Z)(I_n)
 \end{align*}
from which we get the estimates
\begin{equation}\label{est}
\Gamma_{1,k}(Z,Z)(I_{\cG^n})  \preceq (1/\epsilon_0) K_k(Z,Z)(I_n), \quad
\Gamma_{2,k}(Z,Z)(I_n) \preceq  K_k(Z,Z)(I_n)
\end{equation}
for $k = 1,2, 3, \dots$.  As the sequence $K_k$ is converging to $K$  in $\fX$-norm, it follows that
$\| K_k(Z,Z)(I_n)\|$ is uniformly bounded in norm with respect to $k=1,2,\dots$.  We conclude that both
collections $\Gamma_{1,k}(Z,Z)(I_n)$ and $\Gamma_{2,k}(Z,Z)(1_{\cA^{n \times n}})$ are 
uniformly bounded in norm (with respect to $k = 1, 2, \dots$).  
 As $\Gamma_{1,k}$ and $\Gamma_{2,k}$ are cp, 
it is known that
\begin{align*}
&  \| \Gamma_{1,k}(Z,Z) \|_{\cL( \cL(\cN)^{n \times n}, \cL(\cE)^{n \times n})} = 
  \| \Gamma_{1,k}(Z,Z)(I_{\cN^n})\|_{\cL(\cE^n)}   \\
&    \| \Gamma_{2,k}(Z,Z) \|_{\cL({\mathbb C}^{n \times n}, \cL(\cE)^{n \times n})} = 
  \| \Gamma_{2,k}(Z,Z)(I_n)\|_{\cL(\cE^n)}
\end{align*}
We conclude that the collection of numbers
\begin{equation}  \label{est1}
\{ \| \Gamma_{1,k}(Z,Z)\|,\, \| \Gamma_{2,k}(Z,Z)\| \colon Z  \in \Omega,\, k=1,2,\dots\}
\end{equation}
is  bounded.

 A consequence of the {\em respects direct sums} property for nc kernels is that
$$
 \fQ\left( \sbm{ Z & 0 \\ 0 & W }, \sbm{ Z & 0 \\ 0 & W } \right) 
 \left( \sbm{I_n & 0 \\ 0 & I_m }\right)
  =  \sbm{\fQ(Z,Z)(I_n) & 0 \\ 0 & \fQ(W,W)(I_m) }
 $$
 for $Z \in \Omega_{n}$, $W \in \Omega_{m}$.  Consequently, if $Z, W \in {\mathbb P}_{\fQ}$, then also
 $\sbm{ Z & 0 \\ 0 & W } \in {\mathbb P}_\fQ$.  Then the same argument as given above leads to the conclusion
 that  the set of numbers, a priori larger than that in \eqref{est1},
 \begin{equation}  \label{est2}
\{ \| \Gamma_{1,k}\left( \sbm{Z & 0 \\ 0 & W}, \sbm{ Z & 0 \\ 0 & W} \right)\|,\,
\| \Gamma_{2,k}\left( \sbm{Z & 0 \\ 0 & W}, \sbm{ Z & 0 \\ 0 & W}  \right )\| \colon Z, W  \in \Omega_0,\, N \in 
{\mathbb N}\}
\end{equation}
is bounded.

 Another consequence of the {\em respects direct sums}  property of nc kernels is that,  for any nc kernel $\Gamma$,
 $$
 \Gamma \left( \sbm{ Z & 0 \\ 0 & W}, \sbm{ Z & 0 \\ 0 & W} \right) \left( \sbm{ 0 & P \\ 0 & 0} \right) =
 \sbm{ 0 & \Gamma(Z,W)(P) \\ 0 & 0 }.
 $$
 From this identity combined with the boundedness of the set \eqref{est2} we can read off that the set
 \begin{equation}  \label{est3}
 \{ \| \Gamma_{1,k}(Z,W) \|, \, \| \Gamma_{2,k}(Z,W) \| \colon Z,W \in \Omega, \, k \in {\mathbb N} \}
 \end{equation}
 is uniformly bounded as well.

By the aforementioned pre-compactness in the norm-topology of norm-bounded subsets in the Banach space
$\fX$, we conclude that there exists a subsequence $\{k_j\}_{j \in {\mathbb N}}$ so that
$\Gamma_{2, k_j}$ converges to some kernel $\Gamma_2 \in \fX$  in the norm topology of $\fX$.

Note that the kernels $\Gamma_{1,k}$
sit in the space $\fX'$ consisting of graded kernels from $\Omega \times \Omega$ to 
$\cL(\cL(\cN)_{\rm nc}, \cL(\cE)_{\rm nc})$, i.e., just the space $\fX$ but with the  $C^*$-algebra ${\mathbb C}$ replaced
by the $C^*$-algebra $\cL(\cN)$.  This space is infinite-dimensional whenever $\cN$ is infinite-dimensional. 
To handle this term, we need the following more sophisticated argument. 
It can be shown that {\em $\fX'$  is a 
dual Banach space and hence is equipped with a weak-$*$ topology} and   furthermore, {\em convergence of a net $K_\alpha$
to $K$ in the weak-$*$ topology of $\cX$ concretely just means pointwise convergence of  $K_\alpha(Z,W)(P)$ to $K(Z,W)(P)$
in the weak-$*$ topology of $\cL(\cE^m, \cE^n)$}
 (in fact  the norm topology for our setting here since $\cE$ is finite-dimensional).
For additional information on this point we refer to \cite[Corollary 2 page 230]{Paulsen} and \cite[Section IV.2]{Takesaki} 
as well as \cite[Section 4]{BMV2} for an application very similar to what is being done here.
Therefore essentially the same analysis as above applies equally well to the sequence
$\Gamma_{1,N}$ of kernels in $\fX'$ but we may need to drop down to subnets rather than to subsequences.
By dropping down to a subnet  $\{\Gamma_{1, \alpha}\}_{\alpha \in {\mathfrak A}}$ if necessary, we may assume that the net
$\Gamma_{1,\alpha}$ converges weak-$*$ to a kernel $\Gamma_1$ in $\fX'$.  
By using the identification of 
weak-$*$ convergence of the kernel-net $\Gamma_\alpha$ as pointwise weak-$*$ convergence of the values-net
 $\Gamma_\alpha(Z,W)(P)$ 
in $\cL(\cE^{n_W}, \cE^{n_Z})$  for each $Z,W \in \Omega$
and $P \in {\mathbb C}_{\rm nc}$, one can see that the weak-$*$ limit of a net of cp nc kernels is again a cp nc kernel.
In this way we see that the limit kernels $\Gamma_1$ and $\Gamma_2$ are again cp nc kernels.  Furthermore,
again making use of the concrete identification of weak-$*$ convergence in $\fX$ and taking the limit along the 
net $\alpha$ in the identity
$$
 K_\alpha(Z,W)(P) = \Gamma_{1, \alpha}(Z,W)\left( \fQ(Z,W)(P)\right) + \Gamma_{2, \alpha}(Z,W)(P),
 $$
 we see that $\Gamma_1$ and $\Gamma_2$ give the desired representation \eqref{cone-form} for membership
 of the limit kernel $K$ in the cone $\cC$, as needed to complete the proof of Lemma \ref{L:cone}.
\end{proof}

\textbf{The cone separation argument.}
Recall that $\mathfrak{X}$ is the finite dimensional Banach space  of graded kernels given by \eqref{mathfrakX}. 
We established in Lemma \ref{L:cone} that the cone of nc kernels $\mathcal{C} \subset \mathfrak{X}$ defined by
 \eqref{cone-form} is both closed and pointed.

We now show by contradiction that the kernel $\mathfrak{S}$ lies in the cone $\mathcal{C}$. Suppose that $\mathfrak{S}$ 
does not lie in the cone $\mathcal{C}$. Then there exists a separating linear functional $\ell_1$ on $\fX$ such that 
$\ell_1(\mathfrak{S})<0$ and $\ell_1(\mathcal{C}) \geq 0$ \cite{Rudin}. Since $\cC$ is a closed pointed cone, there exists 
another linear functional $\ell_2$ on $\fX$ such that $\ell_2(\cC \backslash \{0\} ) > 0$ \cite[Theorem 3.38]{Jahn}. 
As $\ell_1(\mathfrak{S}) < 0$, we may choose $\epsilon > 0$ so that 
$$
  \ell_1({\mathfrak S}) + \epsilon \ell_2({\mathfrak S}) < 0.
$$
Let us then set $\ell = \ell_1 + \ell_2$.  Then $\ell$  is a linear function on $\cX$ with the property that
\begin{equation} \label{ell}
\ell(\fS)<0 \text{ and } \ell(\cC \backslash \{0\} ) > 0. 
\end{equation}

For $f: \Omega \to \cE_{\rm nc}$ a graded function on $\Omega$, given $Z \in \Omega_n$, we see that 
$f(Z) \in \cE^{n \times n} \simeq \cL(\bbC,\cE)^{n \times n}$. We define the function $f^*: \Omega \to \cL(\cE, \bbC)^{n \times n}$ by setting $f^*(Z)=f(Z)^*$, and we construct vector spaces $\cG$ and $\cH$ by 
\begin{align*}
& \cG = \{ f^*  \, | \,  f\colon \Omega \to \cE_{\rm nc}  \text{ a graded function} \} \\
&\cH = \{ f^*  \, | \, f \colon \Omega \to \cE_{\rm nc}  \text{ a nc function} \}.
\end{align*} 
As any nc function is graded by definition, we have the vector-space inclusion $\cH \subset \cG$.  Since $\Omega$
is a finite set of points and $\cE$ is finite-dimensional, it is easily checked that both $\cG$ and $\cH$ are
finite-dimensional.

We show that the linear functional $\ell$ induces an inner product on $\cH$.
 Notice that any $f^*,g^* \in \cH$ induce a nc kernel $K_{f,g}: \Omega \times \Omega \to \cL(\bbC_{\rm nc},\cL(\cE)_{\rm nc})$ by
\begin{equation*}
K_{f,g}(Z,W)(P)=f(Z)P g(W)^*.
\end{equation*}
With the linear functional $\ell$, we define an inner product on $\cH$ by
\begin{equation}\label{Hip}
\ip{f^*,g^*}_{\cH} = \ell(K_{f,g}).
\end{equation}
Note that $\ip{f^*,f^*}_{\cH} = \ell(K_{f,f}) \geq 0$ by the fact that $K_{f,f}$ is in the cone $\cC$ \eqref{cone-form} and 
$\ip{f^*,f^*}_{\cH}=0 \Rightarrow f^*=0$ by \eqref{ell}.  Then $\cH$ equipped with this inner product becomes a
Hilbert space, still denoted as $\cH$.

As in \eqref{Omega-list},
let us denote the finite set of points $\Omega$ by $\Omega = \{ Z^{(1)}, \dots, Z^{(N)}\}$ with $Z^{(i)} \in \Omega_{n_i}$.  
We let $Z^{(0)}$ denote the direct sum  
\begin{equation} \label{Z0}
Z^{(0)} = \bigoplus_i^N  Z^{(i)} \in (\bbC^d)^{N _0 \times N_0} \text{ where } N_0 = \sum_{i-=1}^N n_i.
\end{equation}
We define a mapping
\begin{equation*}
\cI: \cG \to \cL(\cE^{N_0}, \bbC^{N_0})
\end{equation*}
by 
\begin{equation*}
\cI: f^* \mapsto f(Z^{(0)})^*.
\end{equation*}
Let $\widehat \cG := \cI(\cG) \subset \cL(\cE^{N_0}, \bbC^{N_0 \times N_0})$.
We make $\widehat \cG$ a Hilbert space by equipping it with the inner product
\begin{equation}   \label{Gip}
\ip{\cI(f^*),\cI(g^*)}=\textrm{tr}(f(Z^{(0)})^* g(Z^{(0)}))
\end{equation}

The map $\cI$ is a vector space isomorphism between $\cG$ and $\widehat{\cG}$. We let $\widehat{\cH}=\cI(\cH)$, and we 
let $\cI_0$ be the restriction of $\cI$ to $\cH$.  Then $\cI_0$ is a vector-space isomorphism between the
two finite-dimensional Hilbert spaces $\cH$ and $\widehat{\cH}$.

We shall need higher-multiplicity versions of the spaces $\cG$ and $\cH$ defined as follows.  For $\cX$ an auxiliary Hilbert space, we let
\begin{align*}
\cG_\cX = \{ f^* \, | \, f \colon \Omega \to \cL(\cX, \cE)_{\rm nc} \text{ graded}\}, \\
\cH_\cX = \{ f^* \, | \, f \colon \Omega \to \cL(\cX, \cE)_{\rm nc} \text{ a nc function} \}.
\end{align*}
The space $\cH_\cX$ can be given an inner-product again defined by \eqref{Hip} where now the associated kernel $K_{f,g}$ is equal to 
the nc kernel with values in $\cL(\cL(\cX)_{\rm nc}, \cL(\cE)_{\rm nc})$.  As explained on pages 78-79 of \cite{BMV2}, we have
the identifications
\begin{equation}   \label{tensor-identification}
\cG_\cX \cong \cG \otimes \cX, \quad \cH_\cX \cong \cH \otimes \cX.
\end{equation}

Given a nc function $Q \colon \Omega \to \cL(\cR, \cS)_{\rm nc}$ for two auxiliary Hilbert spaces $\cR$ and $\cS$, we let $M_{Q}$
be the operator with adjoint $M^*_{Q}$ equal to the  multiplication operator given by
$$
  M^*_Q \colon f(W)^* \mapsto Q(W)^* f(W)^*.
$$
It is routine to check from the definitions that
\begin{itemize}
\item if $Q$ is graded, then $M^*_Q$ maps $\cG_\cS$ to $\cG_\cR$, and
\item if $Q$ is a nc function, then $M^*_Q$ maps $\cH_\cS$ to $\cH_\cR$.
\end{itemize}
The restriction of $M^*_Q$ to $\cH_\cS$ with target space taken to be $\widehat \cH_\cR$ will be denoted by
$M^*_{Q,0}$.  Note that here $\cH_\cS$ and $\cH_\cR$ have inner products induced by \eqref{Hip} while $\cG_\cS$ and $\cG_\cR$
have their own quite different inner products induced by \eqref{Gip}.  In particular, the term subspace here is used loosely:  $\cH_\cS$ is a vector subspace
of $\cG_\cS$ but not a Hilbert subspace as the inner product on $\cH_\cS$ is not inherited from $\cG_\cS$ as a subset of $\cG_\cS$.

We define higher-multiplicity versions of the spaces $\widehat \cG$ and $\widehat \cH$ as follows.  
Given an auxiliary Hilbert space $\cX$ together with a function $f \in \cG_\cX$ and the point $Z^{(0)} = \bigoplus_{i=1}^N Z^{(i)}$ 
as as in \eqref{Z0}, define 
$$
\cI_\cX \colon f^* \mapsto \bigoplus_{1}^N f(Z^{(i)})^* = f(Z^{(0)})^* \in \bigoplus_{i=1}^N \cL(\cE^{n_i}, \cX^{n_i}).
$$
Note that as $f$ sweeps the space of all graded functions on $\Omega$ with values in $\cL(\cE, \cX)_{\rm nc}$, the 
resulting space of values $\{\bigoplus_!^N f(Z^{(i)})^* \colon f \in \cG_\cX\}$ sweeps exactly the space 
\begin{equation}   \label{hatGX}
\widehat \cG_\cX:= \bigoplus_{i=1}^N \cL(\cE^{n_i}, \cX^{n_i}).
\end{equation}
 We introduce the Hilbert-Schmidt inner product on $\widehat \cG_\cX$  given by
\begin{equation}   \label{hat-ip'}
\langle f^*, g^* \rangle = \operatorname{tr} (f(Z^{(0)})^* g(Z^{(0)})).  
\end{equation}
to make $\widehat \cG_\cX$ a Hilbert space.
Note that 
$$
 \operatorname{tr} (f(Z^{(0)}) g(Z^{(0)})^*) =  \operatorname{tr} (g(Z^{(0)})^*f(Z^{(0)}))
 $$
 where $g(Z^{(0)})^* f(Z^{(0)}) \in \cL(\cE^{N_0})$ is a finite-dimensional operator even if $\cX$ is infinite-dimensional, so indeed 
 $f(Z^{(0)}) g(Z^{(0)})^*$ is in the trace class and the inner product is well-defined.  Note also that the direct sum decomposition
 in \eqref{hatGX} is orthogonal and that the map $\cI_\cX$  between $\cG_\cX$ and $\widehat \cG_\cX$ is unitary.
 Making use of the
 tensor-product identifications mentioned in \eqref{tensor-identification}, the map $\cI_\cX$ can be identified with
 $$
  \cI_\cX \cong \cI \otimes I_\cX.
 $$
 We let $\cI_{\cX,0}$ denote the restriction of $\cI_\cX$ to $\cH_\cX$.  Then $\cI_{\cX, 0}$ is a linear isomorphism from
 $\cH_\cX$ onto its image $\cI_{\cX, 0}(\cH_\cX)$, which we denote  as $\widehat \cH_\cX$.  Again making the
 identification of spaces in \eqref{tensor-identification}, one can arrive at the operator identification
 $$
    \cI_{\cX,0} \cong \cI_0 \otimes I_\cX.
 $$

Given two auxiliary Hilbert spaces $\cS$ and $\cR$ and $Q$ a function from $\Omega$ to $\cL(\cR, \cS)$,
we write $\widehat M^*_Q$ for the multiplication operator of the form
\begin{equation}  
\widehat M^*_Q \colon f(Z^{(0)})^* \mapsto Q(Z^{(0)})^* f(Z^{(0)})^*.
\end{equation}
Then it is routine to verify from the definitions that
\begin{itemize}
\item if $Q \colon \Omega \to \cL(\cR, \cS)_{\rm nc}$ is graded, then $ \widehat M^*_Q$ maps $\widehat \cG_\cS$ to 
$\widehat \cG_\cR$, and

\item if $Q \colon \Omega \to \cL(\cR, \cS)_{\rm nc}$ is a nc function, then $ \widehat M^*_Q$ maps $\widehat \cH_\cS$ to 
$\widehat \cH_\cR$.
\end{itemize}

For $Q \colon \Omega \to \cL(\cR, \cS)_{\rm nc}$ a nc function, we let $\widehat M^*_{Q, 0}$ be the restriction of the operator
$\widehat M^*_Q$ to $\widehat \cH_\cS$.  Notice that we have the intertwining relations
\begin{equation}   \label{intertwining}
\cI_\cR \, M^*_Q = \widehat M^*_Q \, \cI_\cS, \quad \cI_{\cR,0} \, M^*_{Q,0} = \widehat M^*_{Q,0}\,  \cI_{\cS,0}.
\end{equation}

Let us identify the operator $\widehat M^*_{Q}$ from $\widehat \cG_\cS$ to $\widehat \cG_\cR$ more explicitly as follows.
First let us note that a nc function $Q \colon \Omega \to \cL(\cR, \cS)_{\rm nc}$,  can always be extended uniquely to a nc function,
again denoted as $Q$,  
on the nc envelope $[\Omega]_{\rm nc}$ by using the {\em respects direct sums} property as a definition 
(see \cite[Proposition 2.17]{BMV1}).  In particular we can extend $Q$ to a nc function on $\Omega \cup \{ Z^{(0)}\}$ by defining
$Q(Z^{(0)})  = \bigoplus_{i=1}^N Q(Z_i)$.
By definition we then have
\begin{equation}   \label{M*Q}
\widehat M_Q^* \colon f(Z^{(0)})^* \mapsto Q(Z^{(0)})^* f(Z^{(0)})^*.
\end{equation}   
Viewing $\widehat \cG_\cS$ as the direct-sum space $\bigoplus_{j=1}^N \cL(\cE^{n_j}, \cS^{n_j})$ and writing
a generic element of $\bigoplus_{j=1}^N \cL(\cE^{n_j}, \cS^{n_j})$ as $\bigoplus_{j=1}^N f^*_j$, we can rewrite  \eqref{M*Q}
as 
$$
\widehat M_Q^* \colon \operatorname{col}_{1 \le j \le N}[ f_j^*]  \mapsto \operatorname{diag}_{1 \le j \le N} 
[Q(Z^{(j)})^* ]  \cdot \operatorname{col}_{1 \le j \le N} [f_j^*].
$$
where the generic element  $f(Z^{(0})^*$ of $\widehat \cG_\cS$ is parametrized by a graded function $f \colon \Omega
\to \cL(\cE, \cS)_{\rm nc}$.  If we write $\widehat \cG_\cS$ and $\widehat \cG_\cR$ in column form
\begin{equation}  \label{col-form}
\widehat \cG_\cS = \operatorname{col}_{1 \le j \le N} \cL(\cE^{n_j}, \cS^{n_j}), \quad
\widehat \cG_{\cR} = \operatorname{col}_{1 \le j \le N} \cL(\cE^{n_j}, \cR^{n_j}),
\end{equation}
we see from the formula \eqref{M*Q} that the action \eqref{M*Q} from $\widehat \cG_\cS$ to $\widehat \cG_\cR$    with 
$\widehat \cG_\cS$ and $\widehat \cG_\cR$ given in column form \eqref{col-form} is given by
$$ 
\operatorname{col}_{1 \le j \le N} f_j^* \mapsto \operatorname{diag}_{1 \le j \le N} Q(Z^{(j)})^* \cdot
\operatorname{col}_{1 \le j \le N} f_j^*.
$$
As by assumption $Q$ is a nc function on $\Omega \cup \{Z^{(0)}\}$, we then have
$$
  \operatorname{diag}_{1 \le j \le N} Q(Z^{(j)})^* = Q(Z^{(0)})^*
$$
and the formula above can be written as
$$
 \operatorname{col}_{1 \le j \le N}[ f_j^*] \mapsto L_{Q(Z^{(0)})^*} \cdot [\operatorname{col}_{1 \le j \le N} f_j^*]
 $$
 where here $f_j^*$ runs through the space $\cL(\cE^{n_j}, \cS^{n_j})$ and where $L_{Q(Z^{(0)})^*}$ denotes the operator
 of left multiplication by $Q(Z^{(0)})^*$ on the space \\
 $\operatorname{col}_{1 \le j \le N} \cL(\cE^{n_j}, \cS^{n_j})$
 mapping into the space $\operatorname{col}_{1 \le j \le N} \cL(\cE^{n_j}, \cR^{n_j})$.
 
 Note the space $\cL(\cE^{n_j}, \cS^{n_j}) \cong \cL(\cE, \cS)^{n_j \times n_j}$ can be thought of as $n_j \times n_j$-block
 matrices with matrix entries equal to operators from $\cS$ to $\cE$.   To get still more granular, let us fix a basis for
$\{e_1, \dots, e_{n_\cE}\}$ for  $\cE$ (where $n_\cE = \dim \cE < \infty$ by (H2))
and represent an element $T$ of  $\cL(\cE, \cS)^{n_j \times n_j}$ as a 
$n_j \times n_j \cdot n_\cE$ block
 matrix with matrix entries in $\cL({\mathbb C}, \cS^{n_j}) \cong \cS^{n_j}$.
   We can decompose $\cL(\cE^{n_j}, \cS^{n_j})$ into
 $n_j \cdot n_\cE$ pairwise orthogonal (in the Hilbert-Schmidt inner product) subspaces, corresponding to the subspaces of
 $\cL(\cE, \cS)^{n_j \times n_j}$ with nonzero matrix entries supported in a single given column.  These subspaces are all invariant for
 $L_{Q(Z^{(j)})^*}$ and each of these $n_j \cdot n_\cE$ restricted operators amounts to the same operator $Q(Z^{(j)})^*$ in the 
 standard presentation of a matrix as an operator acting on a column space $\cS^{n_j}$.  Therefore, if we apply the \textbf{vec}
 operator of stacking the columns of the matrix $T \in \cL({\mathbb C}, \cL(\cS))^{n_j \times n_j \cdot n_\cE}$ down from each other
 to form a single column
 $$
 {\rm \textbf{vec}} ([ T_{ij}]) = \operatorname{col}_{1 \le j \le N} \operatorname{col}_{\sum_{\ell = 1}^{i-1} \le j \le \sum_{\ell=1}^i n_\ell
 \cdot n_\cE} [ T_{ij}],
 $$
 we see that the left multiplication operator $L_{Q(Z^{(0)})^*}$ acting between spaces of Hilbert-Schmidt matrices
 $$
 L_{Q(Z^{(0)})^*} \colon \oplus_{j=1}^N  \cL(\cE^{n_j}, \cS^{n_j})  \mapsto \oplus_{j=1}^N  \cL(\cE^{n_j}, \cR^{n_j})
 $$
 is unitarily  equivalent to a single matrix $Q(\widehat Z^{(0)})^*$ (to be defined in a moment)
 acting on a standard complex column space.
 Here we set
 \begin{equation}  \label{hatZ0}
   \widehat Z^{(0)} = \bigoplus_{j=1}^N  \bigoplus_{1}^{n_j \cdot n_\cE} Z^{(j)}  \in [\Omega]_{\rm nc},
 \end{equation}
 extend the nc function $Q$ on $\Omega$ uniquely to a nc function  on $\Omega \cup \{ \widehat Z^{(0)} \}$ via
 $$
   Q(\widehat Z^{(0)})^* = \bigoplus_{j=1}^N \bigoplus_{1}^{n_j \cdot n_]\cE} Q(Z^{(j)})^* \colon 
   \bigoplus_{j=1}^N \bigg( \bigoplus_1^{n_j \cdot n_\cE} \cS^{n_j}  \bigg) \mapsto 
    \bigoplus_{j=1}^N \bigg( \bigoplus_1^{n_j \cdot n_\cE} \cR^{n_j}   \bigg).
 $$
 The final conclusion is:  
 \begin{equation}  \label{M*Q'}
 \widehat M^*_Q = Q(\widehat Z^{(0)})^*.
 \end{equation}
 where $\widehat Z^{(0)}$ is given by \eqref{hatZ0}.

One particular example of a multiplication operator is given by the nc function $\chi_k: \Omega \to \bbC_{\rm nc}$ ($1 \le k \le d$)
defined by  $\chi_k(Z)=Z_k$ where $ Z=(Z_1,  \dots , Z_d) $ is the decomposition of $Z \in \Omega \subset {\mathbb C}^d$ as an 
element of ${\mathbb C}^d$ (see \eqref{Cdnc}). Everything said above for a general nc function 
$Q \colon \Omega \to \cL(\cR, \cS)_{\rm nc}$ applies in particular to the case $Q = \chi_k$ for each $k=1, \dots, d$, with
now $\cR = \cS = {\mathbb C}$.  Thus, for each $k=1, \dots, d$ we associate the operators
\begin{align*}
& M^*_{\chi_k,0} \colon f(W)^* \mapsto W_k^* f(W)^* \text{ on } \cH,  \\
& M^*_{\chi_k} \colon f(Z^{(0)})^* \mapsto Z_k^{(0) *} f(Z^{(0)})^* \text{ on } \cG
\end{align*}
where $M^*_{\chi_k,0}$ is just the restriction of $M^*_{\chi_k}$ to $\cH$, and from \eqref{intertwining} we have
\begin{equation}  \label{intertwining-chi}
\cI M^*_{\chi_k} = \widehat M^*_{\chi_k} \cI, \quad \cI_0 M^*_{\chi_k, 0} = \widehat M^*_{\chi_k,0} \cI_0.
\end{equation}
and from \eqref{M*Q'} we have
\begin{equation}    \label{M*chi'}
\widehat M^*_{\chi_k} = \big( \widehat Z^{(0)}_k \big)^*.
\end{equation}

 In general suppose that $\cV$ is a finite-dimensional vector space (say of dimension $n_\cV$) and $T = (T_1, \dots, T_d)$
 is a $d$-tuple of linear transformations on $\cV$ and suppose that $f \colon \Omega \to ({\mathbb C}^N)_d$ is a nc function,
 where $\Omega$ is a subset of ${\mathbb C^d}_{\rm nc} \cong  \amalg_{n=1}^\infty ({\mathbb C}^{n \times n})^d$
 ($d$ tuples of matrices of size $n \times n$ over all possible $n \in {\mathbb N}$).  Next choose a vector space isomorphism
 ${\mathbb S} \colon \cV \to {\mathbb C}^{n_\cV}$.  We choose a fixed basis of ${\mathbb C}^{n_\cV}$ to identify
 linear transformations on ${\mathbb C}^{n_\cV}$ with $n_\cV \times n_\cV$ matrices.    Given a $d$-tuple of linear transformations
 $({\mathbb T}_1, \dots, {\mathbb T}_n)$ acting on $\cV$ such that it happens that the matrix $d$-tuple
 $$
 (T_1, \dots,  T_d) = ({\mathbb S} {\mathbb T}_1{\mathbb S}^{-1}, \dots, {\mathbb S} {\mathbb T}_d {\mathbb S}^{-1})
  \in ({\mathbb C}^{n_\cV \times n_\cV})^d
 $$
 turns out to be in $\Omega$ (the domain of $f$),  let us define $f(T_1, \dots, T_d) \colon \cV \to \cV$ by
 $$
   f(T_1, \dots, T_d) = {\mathbb S}^{-1} f({\mathbb S} {\mathbb T}_1 {\mathbb S}^{-1}, \dots, {\mathbb S} {\mathbb T}_d {\mathbb S}^{-1})
   {\mathbb S}.
 $$ 
 To show that $f(T_1, \dots, T_d)$ is well-defined (i.e., independent of the choice of vector-space isomorphism
 ${\mathbb S} \colon \cV \to {\mathbb C}^{n_\cV}$, we proceed as follows.  Suppose that we had instead used
 vector-space isomorphism ${\mathbb S}' \colon \cV \to {\mathbb C}^{n_\cV}$ and had defined $f(T_1, \dots, T_d)$
 to be instead
 $$
  f'(T_1, \dots, T_d) =  {\mathbb S}'^{-1} f({\mathbb S}' {\mathbb T}_1 {\mathbb S}'^{-1}, \dots, {\mathbb S}' {\mathbb T}_d {\mathbb S}'^{-1})
   {\mathbb S}'.
 $$ 
 To  show that $f'(T_1, \dots, T_d) = f(T_1, \dots, T_d)$ we must show that
 $$
 {\mathbb S}' {\mathbb S}^{-1} f({\mathbb S} {\mathbb T}_1 {\mathbb S}^{-1}, \dots, {\mathbb S} {\mathbb T}_d {\mathbb S}^{-1})
  {\mathbb S}{\mathbb S}'^{-1} =  
 f({\mathbb S} '{\mathbb T}_1 {\mathbb S}'^{-1}, \dots, {\mathbb S}' {\mathbb T}_d {\mathbb S}'^{-1}).
 $$
 But this is an immediate consequence of the {\em respects similarities} property of  $f$ as a function defined on $d$-tuples of
 matrices.
Therefore, as far as the nc functional calculus is concerned, we are free to
 identify a point in $({\mathbb C}^d )_{\rm nc}$ at level $n$ with a $d$-tuple of linear transformations on a finite-dimensional 
 vector space of dimension $n$.

Thus we may identify $\widehat Z^{(0)}$ not only as the operator $d$-tuple on $\widehat \cG$ but also as $d$-tuple of matrices of size 
$M \times M$
($M = \dim \cG)$ and similarly for $\widehat M_{\chi^*}^*$.    
From \eqref{intertwining-chi} we see that 
\begin{equation}   \label{statement0}
\widehat M_\chi = (\widehat M_{\chi_1}, \dots, \widehat M_{\chi_d}) = \widehat Z^{(0)}.
\end{equation}
From the definition \eqref{hatZ0} of $Z^{(0)}$ and the fact that $\Omega = \{ Z^{(1)}, \dots, Z^{(N)} \}$, \eqref{statement0} implies that
\begin{equation}  \label{statement1}
\widehat M_\chi = (\widehat M_{\chi_1}, \dots, \widehat M_{\chi_d}) \in [ \Omega]_{\rm nc}.
\end{equation}
From the first of relations \eqref{intertwining-chi} we read off that
$$
M_{\chi,0} = (M_{\chi_1,0}, \dots, M_{\chi_d,0}) \text{ is similar to } \widehat Z^{(0)}|_{\widehat \cH^*}
$$
where $\widehat \cH^*$ is the result of taking adjoints pointwise on $\widehat \cH$ and then applying the \textbf{vec} operation
to get an invariant subspaces for $\widehat Z^{(0)}$ rather than for $\widehat Z^{(0)*}$.  This last statement then implies that
\begin{equation}   \label{statement2}
M_{\chi, 0} \in [\Omega]_{\rm full},
\end{equation}
where the {\em full envelope} $[\Omega]_{\rm full}$ of $\Omega$ is defined as in Subsection \ref{S:tax}.  We shall see below that in fact it is also the case
that $M_{\chi,0} \in {\mathbb P}_\fQ$ implying that $M_{\chi,0} \in \Omega'$ (see \eqref{Mchipos} below), but at this point we have to 
work with only the knowledge that $M_{\chi,0} \in [\Omega]_{\rm full}$.

Let us now suppose that $Q$ is any nc function defined on the set $\Omega \cup \{ M_{\chi,0}^* \} \subset \Omega_{\rm full}$ with values in $\cL(\cR, \cS)_{\rm nc}$.
We seek  to show that
\begin{equation} \label{eqQ}
Q(M_{\chi, 0})^*  = M^*_{Q,0} \text{ on } \cH.
\end{equation}
We already know a related version of this result
\begin{equation}   \label{know}
Q(\widehat M_\chi)^* = \widehat M^*_Q
\end{equation}
as a consequence of \eqref{M*Q'} combined with \eqref{statement0}.
This then implies that
\begin{equation} \label{know1}
  Q(\widehat M_\chi)^* \cI_{\cS,0} = \widehat M^*_Q \cI_{\cS,0} = \widehat M^*_{Q,0} \cI_{\cS,0}.
\end{equation}
A consequence of the {\em respects intertwining} property for nc functions is a {\em respects invariant subspaces} property
(see \cite{KVV-book}):
$$
Q \bigg( \begin{bmatrix}  T & * \\ 0 & * \end{bmatrix} \bigg) = \begin{bmatrix} Q(T) & * \\ 0 & * \end{bmatrix}.
$$
From this general principle we can see that $Q(\widehat M_\chi)^* \cI_{\cS,0} = Q(\widehat M_{\chi,0})^* \cI_{\cS,0}$ 
and \eqref{know1} becomes
\begin{equation}  \label{know2}
Q(\widehat M_{\chi,0})^* \cI_{\cS,0} = \widehat M^*_{Q,0} \cI_{\cS,0}.
\end{equation}
From the second of relations \eqref{intertwining} we know that 
$\widehat M^*_{Q,0} \cI_{\cS,0} = \cI_{\cR,0} M^*_{Q,0}$.  From the second set of intertwining relations
\eqref{intertwining-chi} combined with the fact that nc functions respect intertwining conditions we get
$ Q(\widehat M_{\chi,0})^* \cI_{\cS,0} = \cI_{\cR,0}  Q(M_{\chi,0})^*$.  Plugging this information back into
\eqref{know2} leaves us with
\begin{equation} \label{know3}
\cI_{\cR,0} Q(M_{\chi,0})^* = \cI_{\cR,0} M^*_{Q,0}.
\end{equation}
Cancelling off the injective factor $\cI_{\cR,0}$ finally gets us to \eqref{eqQ} as wanted.

Since we have now established that $M_{\chi,0} \in [\Omega]_{\rm full}$ (see \eqref{statement2}), assumption (iii) in the statement of
Theorem \ref{T:ker-dom} tells us that $\fQ$ is decomposable  on $\widetilde \Omega:=
\Omega \cup \{ M_{\chi,0} \}$.  Hence, 
for $Z \in \widetilde  \Omega_n$, $W \in \widetilde \Omega_m$, $P \in {\mathbb C}^{n \times m}$ we can write
\begin{equation}   \label{fQ}
\fQ(Z,W)(P) = Q_+(Z) (P \otimes I_{\cM_+}) Q_+(W)^* - Q_-(Z) (P \otimes I_{\cM_-}) Q_-(W)^*
\end{equation}
for some nc functions $Q_{\pm} \colon \widetilde \Omega  \to \cL(\cM_\pm, \cN)_{\rm nc}$ for some auxiliary Hilbert spaces
$\cM_\pm$,  

 Let us use this decomposition for $\fQ$ to show that
\begin{equation} \label{MchifQ}
\fQ(M_{\chi,0}, M_{\chi,0})(I_\cH) \succ 0.
\end{equation}
Indeed, for any nc function $ f \colon \Omega \to \cL({\mathbb C}, \cH)_{\rm nc}$ generating a generic element $f^*$ of $\cH$,
let us compute
\begin{align}
\langle \fQ(M_{\chi,0}, M_{\chi,0}) f^*, f^* \rangle & =
\langle Q_+(M_{\chi,0})^* f^*, \, Q_+(M_{\chi,0})^* f^* \rangle  \notag  \\
& \quad  - \langle Q_-(M_{\chi,0})^* f^*, Q_-(M_{\chi,0})^* f^* \rangle  \notag  \\
& = \langle M^*_{Q_+,0} f^*, M^*_{Q_+,0} f^* \rangle - \langle M^*_{Q_-,0} f^*, M^*_{Q_-,0} f^* \rangle   \notag \\
& = \ell(K_{\fQ,f,f})
\label{fQ-quadform}
\end{align}
where, now for $Z, W$ only in $\Omega$, 
\begin{align*}
K_{\fQ,f,f}(Z,W)(P) &= f(Z) \bigg(Q_+(Z) (P \otimes I_{\cM_+}) Q_+(W)^*   \\
& \quad \quad - Q_-(Z) (P \otimes I_{\cM_-} )Q_-(W)^* \bigg)f(W)^* \\
& = f(Z) \fQ(Z,W)(P) f(W)^*.
\end{align*}
From this representation we see that $K_{\fQ,f,f}$ has the form of the first term on the right-hand side of \eqref{cone-form}
and hence is in the cone $\cC$.  Hence either the kernel $K_{\fQ,f,f}$ is identically zero or $\ell(K_{\fQ,f,f}) > 0$.
In the first case, it then follows in particular that
\begin{equation}    \label{fQf}
  f(Z) \fQ(Z,Z)(I_n) f(Z)^* = 0 \text{ for all } Z \in \Omega_n, \, n \in {\mathbb N}.
\end{equation}
By assumption $\Omega \subset {\mathbb P}_\fQ$, so $\fQ(Z,Z)(I_n)$ is strictly positive definite.  
Then \eqref{fQf} implies that $f(Z) = 0$ for all $Z \in \Omega$, and hence $f^*$ is the zero element of $\cH$.
Thus, if $f^*$ is not the zero element of  $\cH$ we have that $\ell(K_{\fQ, f, f}) > 0$ which translates to positivity of
the quadratic form:
$\langle \fQ(M_{\chi,0}, M_{\chi,0})(I_\cH) f^*,  f^* \rangle_\cH > 0  \text{ for } 0 \ne f \in \cH$.
As $\cH$ is finite-dimensional, this is just the statement that 
\begin{equation}   \label{Mchipos}
\fQ(M_{\chi,0}, M_{\chi,0})(I_\cH) \succ 0
\end{equation}
as an operator on $\cH$, i.e., $M_{\chi,0} \in {\mathbb P}_\fQ$.
 As we have already
established that $M_{\chi,0} \in [\Omega]_{\rm full}$ (see \eqref{statement2}), we now know that
$M_{\chi,0} \in \Omega':= [\Omega]_{\rm full} \cap {\mathbb P}_\fQ$.

Consequently, we can now invoke condition (iv) in the statement of the theorem 
to conclude that $\fS$ is also decomposable on $\Omega \cup\{ M_{\chi,0}\}$.
Hence we can find
nc functions $S_\pm \colon  \Omega \cup \{ M_{\chi, 0}\} \to \cL(\cD_\pm, \cE)_{\rm nc}$ for some additional
auxiliary Hilbert spaces $\cD_\pm$ so that
\begin{equation}   \label{fS}
\fS(Z,W)(P) = S_+(Z)(P \otimes I_{\cD_+}) S_+(W)^* - S_-(Z) (P \otimes I_{\cD_-}) S_-(W)^*
\end{equation}
for $Z, W \in \widetilde \Omega$.
As we have already observed that $\fQ(M_{\chi,0}, M_{\chi,0})(I_\cH) \succ 0$, we now can apply the standing kernel-dominance 
hypothesis \eqref{ker-dom-hy} to conclude that
$$
  \fS(M_{\chi,0}, M_{\chi,0})(I_\cH) \succeq 0.
 $$
  Using the assumed decomposition \eqref{fS} for $\fS$, by a repeat of the computation \eqref{fQ-quadform} with $\fS$
 in place of $\fQ$ we see that
 $$
0 \le  \langle \fS(M_{\chi,0}, M_{\chi,0}) f^*, f^* \rangle = \ell(K_{\fS,f,f})
 $$
 where, now again  for $Z, W$ only in $\Omega$, 
 $$
 K_{\fS,f,f}(Z,W)(P) = f(Z) \fS(Z,W)(P) f(W)^*.
 $$
 In particular, we make take $1_\cE^*\in \cH_\cE$ with $1_\cE(Z) = I_{\cE^n}$ for $Z \in \Omega_n$ to conclude that
 $$
 \ell(\fS) = \langle \fS(M_{\chi,0}, M_{\chi,0}) 1_\cE^*, 1_\cE^* \rangle_\cH \ge 0,
 $$
 in contradiction with  our earlier conclusion that $\ell(\fS) < 0$.  The existence of an $\ell$ with $\ell(\fS) < 0$ was 
 immediate from the supposition that $\fS$ was not in the cone $\cC$.  We conclude that in the presence of the
 kernel-dominance condition \eqref{ker-dom-hy}, it must be the case that $\fS$ is in $\cC$, i.e., that $\fS$ has a representation as
 in \eqref{cone-rep}.  This finally completes the proof of Theorem \ref{T:ker-dom} for the case where $\Omega$ consists of 
 only finitely many points.  The proof of the general case will be completed in Section \ref{S:KuroshDom} to come.
\end{proof}

We note that Theorem \ref{T:ker-dom} handles only the case where $\Omega$ consists of finitely many points.
The general case can be handled by using a theorem of Kurosh to reduce the general case to the finite-point case.
This will be taken up in the next Section.

 \section{Extensions of results for a finite set of points to a general set of points via the Kurosh method}  
 \label{S:Kurosh}

The theorem of Kurosh  (see \cite[Theorem 2.56]{AMcC-book} as well as \cite[pages 73-75]{AP}) asserts that the limit of an inverse spectrum of nonempty compact sets is a nonempty compact set.  In practice one is given a family of compact sets $\fX_\mu$
indexed by a directed set $\fA$ ($\mu \in \fA$). Here $\fA$ being a {\em directed set} means that $\fA$ is equipped with a
partial order $\preceq$ satisfying reflexivity and transitivity
\begin{align*}
& \alpha \preceq \beta \text{ and } \beta \preceq \alpha \Rightarrow \alpha = \alpha,  \\
& \alpha \preceq \beta \text{ and } \beta \preceq \gamma \Rightarrow \alpha \preceq \gamma.
\end{align*}
as well as
\begin{equation}  \label{dir-sys-axiom}
\text{ given } \alpha, \beta \in \fA, \, \exists \,  \gamma \in \fA \text{ so that } \alpha \preceq \gamma \text{ and } \beta \preceq \gamma.
\end{equation}
We suppose that we are given an {\em inverse spectrum}, i.e., a family of nonempty compact subsets ${\mathbb K}_\alpha
\subset \fX_\alpha$ for each $\alpha \in \fA$ together with a collection of continuous {\em restriction maps} 
$\pi^\alpha_\beta \colon {\mathbb K}_\alpha \to {\mathbb K}_\beta$ for each $\alpha, \beta \in \fA$ with $\beta \preceq \alpha$ such that
\begin{equation}  \label{inv-sys-axiom}
 \pi^\beta_\gamma \circ \pi^\alpha_\beta = \pi^\alpha_\gamma \text{ for } \gamma \preceq \beta \preceq \alpha,
\quad \pi^\alpha_\alpha = {\rm id}_{{\mathbb K}_\alpha} \text{ for all } \alpha \in \fA.
\end{equation}
Then an element $\bGamma = \{\Gamma_\alpha\}_{\alpha \in \fA}$ of the Cartesian product set
 $\Pi_{\alpha \in \fA} {\mathbb K}_\alpha$ is said to be a limit point of the inverse spectrum $\Pi_{\alpha \in \fA} {\mathbb K}_\alpha$
 if it is the case that that $\pi^\alpha_\beta \Gamma_\alpha = \Gamma_\beta$ for all $\alpha, \beta \in \fA$ with $\beta \preceq \alpha$.
 The assertion of the Kurosh theorem is that, with all the compactness assumptions and definitions as listed above, it is always
 the case that the set of limit points of such an inverse spectrum is nonempty.

\subsection{Decomposition of kernels}  \label{S:KuroshDecom}
Let us consider the decomposability problem for classical Hermitian kernels on an infinite point set.
Thus we let $\Omega$ be an abstract point set (considered as distinct $1 \times 1$ matrices over some vector space
to fit into the noncommutative kernel theory) and we consider a complex-valued function $K \colon \Omega \times \Omega \to 
{\mathbb C}$ with $K(z, w)^*  = K(w, z)$.  This can be considered as a special case of the theory of nc kernels by
taking $\cA$ to be ${\mathbb C}$, $\cY = {\mathbb C} = \cL(\cY)$, and identifying $K(z,w)$ with $K(z,w)(1)$.
As a consequence of the Corollary in Example 2 (or by elementary linear algebra), we know that the restriction of $K$ to
any finite subset of $\Omega$ ($K|_F$) is decomposable (i.e., can be represented as the difference of two positive kernels).
One can use the Kurosh theorem to try to find a global such decomposition for all of $K$ as follows.

We wish to denote elements of $\fA$ by lower case Greek letters, e.g., $\mu$.  We identify $\mu \in \fA$ with a finite
subset of $\Omega$, denoted as $\Omega_\mu$, so $\fA$ consists of all finite subsets of $\Omega$.  We say that
$\mu \preceq \nu$ in $\fA$ if $\Omega_\mu \subset \Omega_\nu$. 
It is easily checked that $\fA$ so defined is a directed set. Given a $\mu \in \fA$, we let $\fX_\mu$ be the set of all 
pairs of positive kernels $(K_+, K_-)$ defined on $\Omega_\mu \subset \Omega$, and we let ${\mathbb K}_\mu$ consist of all 
positive-kernel pairs $(K_+, K_-)$ defined on $\Omega_\mu \times \Omega_\mu$ so that the given kernel $K$ restricted to 
$\Omega_\mu$ has the form 
$$
   K(z,w) = K_+(z,w) - K_-(z,w) \text{ for } z,w \in \Omega_\mu.
$$
By our remarks above, each ${\mathbb K}_\mu$ so defined is non-empty.  For $\mu \preceq \nu$,  we define the restriction maps 
$\pi^\nu_\mu$ by $\pi^{\nu}_{\mu} \colon (K_+, K_-) \mapsto (K_+|_{\Omega_\mu}, K_[|_{\Omega_\mu})$ for $(K_+, K_-) \in 
{\mathbb K}_{\nu}$.
Then it is easily checked that  $\{{\mathbb K}_\mu \colon \mu \in \fA\}$ with the system of maps $\{ \pi^{\nu}_{\mu} \colon
\mu  \subset \nu \text{ in } \fA\}$ is an inverse spectrum.  If one can find a limit point $\{K_{+, \mu}, K_{-, \mu}\}_{\mu \in \fA}$ of this
inverse spectrum, then one can construct a global decomposition $K = K_+ - K_-$ of $K$ as the difference of positive
kernels on all of $\Omega$ as follows.    Given any two points $z,w$ in $\Omega$, define $K_+(z,w)$ and $K_-(z,w)$ by
$$
  K_+(z,w) = K_{+, \mu}(z,w), \quad K_-(z,w) = K_{-, \mu}(z,w)
$$
where $\mu$ is any choice of element of $\fA$ such that $z,w \in \Omega_\mu$.  One can check that this definition 
is well-defined, that the global functions $K_+$ and $K_-$ so defined are positive kernels (since the condition of 
positivity involves checking at only
finitely many points at a time, and each $K_+^F$ and $K_-^F$ are positive kernels on $\Omega$, and furthermore
$K(z,w) = K_+(z,w) - K_-(z,w)$.  However there is an example due to Schwartz \cite{Schwartz} (see also \cite{Alpay} for
additional information) that this is not always possible, so the Kurosh theorem fails in this case.  Of course the reason 
is that there is a hypothesis missing:  the set of all pairs $(K_+, K_-)$ of cp nc kernels providing a representation
of the given Hermitian kernel $K$ as a difference of cp nc kernels
can fail to be bounded.  Consequently the weak-$*$ topology 
(the ``right" topology for this problem as we shall see in the examples to follow) restricted to this fiber fails to be 
compact, and the Kurosh theorem does not apply.

\subsection{Arveson extension theorem for kernels}  \label{S:KuroshArvExt}

The goal of this section is to prove general Theorem \ref{T:Arv-ext-ker} to the general case where $\Omega$
is allowed to consist of infinitely many points.

\begin{theorem} \label{T:Arv-ext-ker-gen}
Suppose that $K$ is a nc kernel on a set of nc points $\Omega \subset \cV_{\rm nc}$
 with values in $\cL({\mathbb S}_{\rm nc}, \cL(\cY)_{\rm nc})$, where ${\mathbb S}$ is a operator system and $\cL(\cY)$
 is the $C^*$-algebra of operators on the Hilbert space $\cY$.  Let $\cA$ be a $C^*$-algebra containing ${\mathbb S}$.
 Then there exists a cp nc kernel $\widehat K$ 
 $$
  \widehat K \colon \Omega \times \Omega \to \cL(\cA_{\rm nc}, \cL(\cY)_{\rm nc})
 $$
 such that $\widehat K(Z,W)(P) = K(Z,W)(P)$ for all $Z \in \Omega_n$, $W \in \Omega_m$ whenever $P \in {\mathbb S}^{n \times m}$.
 \end{theorem}
 
 \begin{proof}  By Theorem \ref{T:Arv-ext-ker} we know that the result is true in case $\Omega$ is a finite set of points.
 We wish to reduce the general case to the finite-point case by applying the Theorem of Kurosh.
 
Toward this goal, we set up the directed set $\fA$ exactly as in Section \ref{S:KuroshDecom}.  We let $\fX_\mu$ be the linear space
of all kernels $\Gamma \colon \Omega_\mu \times \Omega_\mu \to \cL(\cA_{\rm nc}, \cL(\cY)_{\rm nc})$.  
For $Z \in \Omega_{\mu, n_Z}$ (i.e., $Z \in \Omega_\mu$ has size $n_Z \times n_Z$), we have that $\Gamma(Z,W)$
is an element of the Banach space $\cL(\cA^{n_Z \times n_W}, \cL(\cY)^{n_Z \times n_W})$, and the fiber $\fX_\mu$
can be identified with the direct sum of Banach spaces
$$
   \fX_\mu = \bigoplus_{Z,W \in \Omega_\mu}  \cL(\cA^{n_Z \times n_W}, \cL(\cY)^{n_Z \times n_W}).
$$
If we endow $\fX_\mu$ with the supremum norm
\begin{equation}   \label{sup-norm}
\| \Gamma \| = \underset{Z,W \in \Omega_\mu} \sup \| \Gamma(Z,W) \|,
\end{equation}
then $\fX_\mu$ is a Banach space.  To finish the proof, we wish to make explicit the following additional information.

\begin{remark} \label{R:weak*}
A Banach space of the form $\cL(\cA^{n \times m}, \cL(\cY)^{n \times m})$ has a predual
which we denote as $\cL(\cA^{n \times m}, \cL(\cY)^{n \times m})_*$ such that the weak-$*$ topology
on   $\cL(\cA^{n \times m}, \cL(\cY)^{n \times m})$ is the same as the pointwise weak-$*$ topology, i.e.,
convergence of a net $\lambda \to \phi_\lambda$ with $\phi_\lambda \in \cL(\cA^{n \times m}, \cL(\cY)^{n \times m})$ 
to $\phi$ in the weak-$*$ topology is equivalent to pointwise weak-$*$ convergence  on $\cL(\cY)^{n \times m}$:
$\phi_\lambda(T)$ converges to  $\phi(T)$  in the weak-$*$ topology of $\cL(\cY)^{n \times m}$ 
for each $T \in \cA^{n \times m}$.  As weak and weak-$*$ topologies agree on bounded sets, this topology is sometimes called the bounded-weak topology or BW-topology (see \cite{Paulsen} for more details).  It is possible to identify the pre-dual spaces
$\cL(\cA^{n \times m}, \cL(\cY)^{n \times m})_*$ more explicitly, but we shall not need this.
\end{remark}

Since the direct summands of $\fX_\mu$ have a predual $\cL(\cA^{n \times m}, \cL(\cY)^{n \times m})_*$, so also does $\fX_\mu$ itself:
$$
(\fX_\mu)_* = \bigoplus_{Z,W \in \Omega_\mu }\bigg(\cL(\cA^{N_Z \times N_W}, \cL(\cY)^{n_Z \times n_W})\bigg)_*
$$
with pairing
$$
 \langle \Gamma, \tau \rangle = \sum_{z,W \in \Omega_\mu} \langle \Gamma(Z,W), \tau(Z,W) \rangle
$$
for $\Gamma = \bigoplus_{Z,W \in \Omega_|mu} \in \fX_\mu$ and $\tau = \bigoplus_{Z,W \in \Omega_\mu} \tau(Z,W) \in
(\fX_\mu)_*$.

For each $\mu \in \fA$, we let ${\mathbb K}_\mu$ be the subset of $\fX_\mu$ consisting of all cp nc kernels 
$$
\Gamma \colon \Omega_\mu \times \Omega_\mu \to \cL(\cA_{\rm nc}, \cL(\cY)_{\rm nc})
$$
so that, for $Z, W \in \Omega_\mu$, 
\begin{equation}  \label{ext-prop}
  K_\mu(Z,W)(P) = K(Z,W)(P) \text{ for } P \in {\mathbb S}^{n_Z \times n_W}.
\end{equation}
Note that since $\Omega_\mu$ is a finite set it is a consequence of Theorem \ref{T:Arv-ext-ker} that each
set ${\mathbb K}_\mu$ is nonempty.

The crucial next point is to show that each ${\mathbb K}_\mu$ is compact (the part that failed for the 
kernel-decomposition problem in Section \ref{S:KuroshDecom}).  By the Banach-Alaoglu Theorem, norm-closed 
and bounded sets are pre-compact in the
weak-$*$ topology on the dual of a Banach space.  Thus, to show that ${\mathbb K}_\mu$ is compact in the 
weak-$*$ topology,
it suffices to show that ${\mathbb K}_\mu$ is weak-$*$ closed (which implies norm closed) and bounded.

For boundedness, it suffices to show that $\| \Gamma(Z,W) \|$ is uniformly bounded for $Z, W \in \Omega_\mu$ and $\Gamma
\in {\mathbb K}_\mu$.  Let $\Omega_\mu = \{ Z_1, \dots, Z_N \}$ where $Z_i \in \Omega_{n_i}$.  Consider the direct sum
$Z^{(0)} = \bigoplus_{i-1}^N Z_i \in \Omega_{N_0}$ where $N_0 = \sum_{i=1}^n n_i$.  We have
\begin{align*}
\| \Gamma(Z_i, Z_j) \|  & \le \bigg\| \Gamma \bigg( \begin{bmatrix} Z_i & 0 \\ 0 & Z_j \end{bmatrix},
\begin{bmatrix} Z_i & 0 \\ 0 & Z_j \end{bmatrix} \bigg) \bigg\|  \le \| \Gamma(Z^{(0)}, Z^{(0)} ) \|  \\
&  = \| \Gamma(Z^{(0)}, Z^{(0)}) \| =\|  \Gamma(Z^{(0)}, Z^{(0)})(I_{N_0}) \|   \\
& = \| K(Z^{(0)}, Z^{(0)} )(I_{N_0}) \|
\end{align*}
where the first two inequalities follow from the fact that $\Gamma$ respects direct sums and the next equality
is a consequence of $\Gamma(Z^{(0)}, Z^{(0)})$ being a completely positive map.

To analyze closedness  of ${\mathbb K}_\mu$ in the weak-$*$ topology,  let $\{ \Gamma_\lambda \}$ be a net in ${\mathbb K}_\mu$
which is weak-$*$ convergent to $\Gamma \in \fX_\mu$.  We must show that in fact $\Gamma$ is back in ${\mathbb K}_\mu$.
By Remark \ref{R:weak*}, it suffices to work with the weak-$*$ topology on $\cL(\cY)$, i.e., we can say that
 $\Gamma_\lambda(Z,W) \underset{ {\rm weak-}*} \to \Gamma(Z,W)$ if and only if 
 $$
 \langle y, \Gamma_\lambda(Z,W)(P) y' \rangle \to \langle y, \Gamma(Z,W)(P) y' \rangle
 $$
 for every $P \in \cA^{n \times m}$, $y \in \cY^n$, $y' \in \cY^m$.  To show that $\Gamma \in {\mathbb K}_\mu$, we must show that
 
  \smallskip
 
 \noindent
 \textbf{(i)} $\Gamma$ is an extension of $K$ from ${\mathbb S}_{\rm nc}$ to $\cA_{\rm nc}$, i.e.
 $K_\mu = \Gamma$ satisfies \eqref{ext-prop}, and
 \smallskip
 
 \noindent
  \textbf{(ii)}  $\Gamma$ is a cp nc kernel.
\smallskip

As for (i), note that $\Gamma_\lambda(Z,W)(P) = K(Z,W)(P)$ for every $Z \in \Omega_{\mu,n}$ and $W \in \Omega_{mu,m}$,
and $P \in {\mathbb S}^{n \times m}$, so in the weak-$*$  limit we also have $\Gamma(Z,W)(P)  = K(Z,W)(P)$.

As for (ii), $\Gamma$ is graded by construction, so it remains to check that $\Gamma$ respects intertwinings and is
completely positive.  Let $Z \in \Omega_{\mu,n}$, $\widetilde Z \in \Omega_{\mu, \widetilde n}$, 
$W \in \Omega_{\mu,m}$, $\widetilde W \in \Omega_{\mu, \widetilde m}$ and suppose that $\alpha Z = \widetilde Z \alpha$
and $\beta W = \widetilde W \beta$ where $\alpha \in {\mathbb C}^{\widetilde n \times n}$ and $\beta \in {\mathbb C}^{\widetilde m
\times m}$.  For $P \in \cA^{n \times m}$, $y \in \cY^n$, $y' \in \cY^{\widetilde m}$, we have
$$
\langle y, (\alpha \Gamma_\lambda(Z,W)(P) \beta^*) y' \rangle = \langle y, \Gamma_\lambda(\widetilde Z, \widetilde W)
(\alpha P \beta^*) y' \rangle.
$$
Taking the limit in $\lambda$ and using weak-$*$ convergence of $\Gamma_\lambda$ to $\Gamma$ in the 
equivalent sense explained in Remark \ref{R:weak*}, we see that the limit in $\lambda$ of this last expression gives us
$$
\langle y, \alpha \Gamma(Z,W)(P) \beta^*) y' \rangle =
\langle y, \Gamma(\widetilde Z, \widetilde W)(\alpha P \beta^*) y' \rangle.
$$
As $y \in \cY^{\widetilde n}$ and $y' \in \cY^{\widetilde m}$ are arbitrary, we conclude that
$\Gamma$  respects intertwining and hence is a noncommutative kernel.
This completes the proof that ${\mathbb K}_\mu$ is compact for each $\mu \in \fA$.

Let us next define the system of restriction mappings $\pi^\nu_\mu \colon {\mathbb K}_\mu \to {\mathbb K}_\nu$ for
$\nu \preceq \mu$ in $\fA$ by
$$
(\pi^\nu_\mu \colon K)(Z,W)(P) = K(Z,W)(P) \text{ for } Z,W \in \Omega_\nu, \, P \in \cA^{n_Z \times n_W}, i.e.,
$$
given an element $K$ of ${\mathbb K}_\nu$ which is kernel on $\Omega_\nu$,  we get a kernel $\pi^\nu_\mu (K)$
on $\Omega_\nu$ simply by restricting the arguments $Z,W \in \Omega_\nu$ to be in the subset $\Omega_\mu
\subset \Omega_\nu$.    Then one can check that that inverse-system-axiom \eqref{inv-sys-axiom} is satisfied.
Thus we can apply the theorem of Kurosh to conclude that there exist limit points for this inverse system, i.e.,
an element $\bGamma = \{\Gamma_\mu\}_{\mu \in \fA}$ in ${\mathbb K}:= \Pi_{\mu \in \fA} {\mathbb K}_\mu$ 
such that $\pi^\nu_\mu \Gamma_\nu = \Gamma_\mu$
for all $\mu \preceq \nu$ in $\fA$.  Given such a $\bGamma$, define a kernel $K$ on all of $\Omega \times \Omega$ 
mapping into $\cL(\cA_{\rm nc}, \cL(\cY)_{\rm nc})$ by
\begin{equation}   \label{K-formula}
  K(Z,W)(P) =  \Gamma_\mu(Z,W)(P)
\end{equation}
for $Z, W \in \Omega$ and $P \in \cA^{n_Z \times n_W}$, where we choose $\mu$ to be any element of $\fA$ such that
both $Z$ and $W$ are in the finite set $\Omega_\mu$.  Using properties of directed sets, one can show that
the formula \eqref{K-formula} defining $K(Z,W)(P)$ is independent of the choice of $\mu$ satisfying $Z, W \in \Omega_\mu$.
As checking whether a given kernel is a cp nc kernel involves only working with finitely many points at a time,
it follows that $K$ is a cp nc kernel on $\Omega$ since each $\Gamma_\mu$ is a cp nc kernel on the finite set $\Omega_\mu$
for each $\mu \in \fA$.  Furthermore, for $P \in {\mathbb S}^{n_Z \times n_W}$,  it is the case that \eqref{ext-prop}
holds with $K_\mu = K$ with $Z, W \in \Omega$ since this is the case with $K_\mu = \Gamma_\mu$ with $Z,W \in \Omega_\mu$
for each $\mu$.  The proof of Theorem \ref{T:Arv-ext-ker-gen} is now complete.
\end{proof}

\begin{remark}  \label{Arv-ext-Paulsen}  We note that the idea of the preceding proof was to use the Kurosh argument to reduce 
the general case to the finite-point case which in turn is handled in Section \ref{S:Arv-ext} by using the results of Subsection 
\ref{S:encoding} to reduce to the Arveson extension theorem for a cp map from an operator space ${\mathbb S}$ to a 
$C^*$-algebra $\cL(\cY)$.
However the proof of the Arveson extension theorem in Paulsen's book \cite[Theorem 7.5]{Paulsen} actually 
proceeds by first
reducing to the case where $\dim \cY < \infty$ and then using finite-dimensional analysis to handle this simpler case
\cite[Theorem 6.2]{Paulsen}.  Our comment here is that it is also possible to reduce the kernel-version of the 
Arveson extension theorem Theorem \ref{T:Arv-ext-ker-gen} directly to the finite-dimensional version of the 
Arveson extension theorem (\cite[Theorem 6.2]{Paulsen}
by applying the Kurosh argument to an inverse-system based on a net $\fA$ indexed by both finite point-sets 
$\Omega_\mu$ and 
finite-dimensional subspaces $\cY_\mu \subset \cY$ ($\mu = (\Omega_\mu, \cY_\mu)$).   This in fact is what we 
do in the next Subsection \ref{S:KuroshDom} coming up, where we show how  the general case of the 
kernel-dominance theorem 
(Theorem \ref{T:ker-dom})  can be reduced to the special case involving both finite-point sets 
$\Omega_\mu \subset \Omega$ and finite-dimensional subspaces $\cE_\mu \subset \cE$.
\end{remark}

\subsection{Kernel dominance theorem}  \label{S:KuroshDom}
In this Section we use the theorem of Kurosh to show how the general case of Theorem \ref{T:ker-dom} can 
be completed by reducing the general case to the special case with hypotheses (H1) and (H2) in force which was completed in
Section \ref{S:ker-dom}.

\begin{proof}
Let us assume that we have the setup of Theorem \ref{T:ker-dom} but without assuming hypotheses (H1) and (H2).
To use the Kurosh theorem we must define the inverse spectrum of interest.   We denote elements of the directed set $\fA$
by lower case Greek letters ($\mu$, $\nu$, $\gamma$, etc.).  To specify an element $\mu$ of $\fA$, we specify
a finite subset $\Omega_\mu$ of $\Omega$ together with a finite-dimensional subspace $\cE_\mu$ of $\cE$.
We then say that $\mu \preceq \nu$ exactly when both
$$
   \Omega_\mu \subset \Omega_\nu \text{ and } \cE_\mu \subset \cE_\nu.
$$
Then it is easily checked that this is a partially ordered set satisfying the directed set axiom \eqref{inv-sys-axiom}.  
For each $\mu \in \fA$,
we let $\fX_\mu$ be the linear space consisting of pairs of cp nc kernels  $(\Gamma_1, \Gamma_2)$
\begin{align}
& \Gamma_1 \colon \Omega_\mu \times \Omega_\mu \to \cL(\cL(\cN)_{\rm nc}, \cL(\cE_\mu)_{\rm nc}),  \notag \\
& \Gamma_2 \colon \Omega_\mu \times \Omega_\mu \to \cL({\mathbb C}_{\rm nc}, \cL(\cE_\mu)_{\rm nc}).
\label{ker-pair-mu}
\end{align}
The set of such pairs forms a Banach space in the norm
$$
\| (\Gamma_1, \Gamma_2) \| = \sup \{ \| \Gamma_1(Z,Z)(I_{\cN^{n_Z}})\|, \,
\| \Gamma_2(Z,Z)(I_{n_Z}) \| \colon Z \in \Omega_\mu \}.
$$
Just as was the case in  Section \ref{S:KuroshArvExt} where $\fX_\mu$ was defined using only one kernel,
the space $\fX_\mu$ is a dual Banach space and the weak-$*$ topology is given  via pointwise weak-$*$ convergence in
$\cL(\cE_\mu)$ after evaluation at a fixed element $T$ in the pre-dual space.

We define the inverse spectrum $\{ {\mathbb K}_\mu \colon \mu \in \fA\}$ with ${\mathbb K}_\mu \subset \fX_\mu$ as follows.
We let ${\mathbb K}_\mu$ consist of all cp nc kernel-pairs $(\Gamma_1, \Gamma_2) $ on $\Omega_\mu$
as in \eqref{ker-pair-mu} so that
\begin{equation}  \label{cone-rep-mu}
P_{\cE_\mu^{n_Z}} \fS(Z,W)(P)|_{\cE_\mu^{n_W}} =
\Gamma_1(Z,W)(\fQ(Z,W)(P)) + \Gamma_2(Z,W)(P).
\end{equation}
The fact that ${\mathbb K}_\mu$ is nonempty for each $\mu$ is the content of the proof of the special case in Section
\ref{S:ker-dom} with hypotheses (H1) and (H2) assumed;   let us note that hypotheses (iii) and (iv) in the statement of 
Theorem \ref{T:ker-dom} are used only to prove the finite-point special case as handled in Section \ref{S:ker-dom} above.

For $\mu, \nu \in \fA$ with $\mu \preceq \nu$, we define the restriction map $\pi^\nu_\mu$ by
$$
\pi^\nu_\mu(\Gamma_1, \Gamma_2) = (\pi_\mu^\nu \Gamma_1, \pi_\mu^\nu \Gamma_2)
$$
where
$$
(\pi_\mu^\nu \Gamma_1)(Z,W)(P) = P_{\cE_\mu^{n_Z}} \Gamma_1(Z,W)(P)|_{\cE_\mu^{n_W}}
$$
for $Z,W \in \Omega_\mu$, $P \in \cL(\cN)^{n_Z \times n_W}$, and 
$$
(\pi_\mu^\nu \Gamma_2)(Z,W)(P) = P_{\cE_\mu^{n_Z}} \Gamma_1(Z,W)(P)|_{\cE_\mu^{n_W}}
$$
for $Z,W \in \Omega_\mu$, $P \in {\mathbb C}^{n_Z \times n_W}$, where here
$P_{\cE_\mu^{n_Z}}$ is the orthogonal projection of $\cE_\nu^{n_Z}$ onto its subspace $\cE_\mu^{n_Z}$.
It is straightforward to check that this projection system $\{ \pi^\nu_\mu \colon \mu \preceq \nu \text{ in } \fA\}$
satisfies the inverse-spectrum axioms \eqref{inv-sys-axiom}

The fact that each ${\mathbb K}_\mu$ is bounded can be seen from the estimate \eqref{est} for the special case 
where  the sequence of kernels $\{K_k\}_{k \in {\mathbb N}}$ is taken to be the fixed kernel 
$P_{\cE_\mu^{n_Z}} \fS(Z,W)(P)|_{\cE_\mu^{n_W}}$.
  Just as in the proof of the Arveson extension theorem in
Section \ref{S:KuroshArvExt}, it follows that ${\mathbb K}_\mu$ is pre-compact in the weak-$*$ topology on $\fX_\mu$.
One can check that the ${\mathbb K}_\mu$ is weak-$*$ closed in $\fX_\mu$ and hence is itself compact
in the weak-$*$ topology inherited from $\fX_\mu$.  Hence we are now in a position to apply the theorem of Kurosh
to conclude that limit points of this inverse system exist.

Let $\boldsymbol{\Gamma} = \{ \Gamma_{\mu, 1}, \Gamma_{\mu, 2} \}$ be any such inverse-system limit point.
We define two kernels $\Gamma_1$ and $\Gamma_2$ defined on all of $\Omega$ via the quadratic form
$$
 \langle \Gamma_i(Z,W)(P)e', e \rangle_{\cY_{n_Z}} \text{ where } e \in \cE^{n_Z}, \, e' \in \cE^{n_W}
$$
for points $Z,W \in \Omega$, $P \in \cN^{n_Z \times n_W}$ if $i=1$ and $P \in {\mathbb C}^{n_Z \times n_W}$ for $i=2$,
$e \in  \cE^{n_Z}$ as follows.  Given $Z,W \in \Omega$ and 
$$
e = \sbm{ e_1 \\ \vdots \\ e_{n_Z}} \in \cE^{n_Z},  \quad   e' = \sbm{ e'_1 \\ \vdots \\ e'_{n_W}} \in \cE^{n_W},
$$
choose any $\mu \in \fA$ so that $Z,W \in \Omega_\mu$, $e_1, \dots, e_{n_Z}, e'_1, \dots, e'_{n_W}  \in \cE_\mu$ and then define
\begin{equation} \label{Gamma-formula}
\langle \Gamma_i(Z,W)(P)e', e \rangle_{\cE^{n_Z}}  = 
 \langle \Gamma_{\mu, i}(Z,W)(P) e', e \rangle_{\cE_\mu^{n_Z}} \text{ for } i=1,2.
\end{equation}
It is a routine check to see that the inverse-system axioms \eqref{inv-sys-axiom} satisfied by the system of 
restriction mappings $\{ \pi^\nu_\mu \colon \mu \preceq \nu \text{ in } \fA \}$ imply that this formula \eqref{Gamma-formula}
is well-defined and uniquely specifies a kernel pair $(\Gamma_1, \Gamma_2)$ having all the desired properties
yielding a global representation \eqref{cone-rep} for $\fS$ as wanted.
\end{proof}

\begin{remark} \label{R:nc-int}
\textbf{Noncommutative interpolation and Schur-Agler class revisited.}  
The special case of Theorem \ref{T:ker-dom} where the kernels $\fQ$ and $\fS$ are assumed to have the special
Hermitian-decomposed form 
\begin{align}   
&\fQ(Z,W)(P) = P \otimes I_\cN - Q(Z) (P \otimes I_\cM) Q(W)^*,  \notag \\
& \fS(Z,W)(P) = a(Z)(P \otimes I_\cE) a(W)^*  - b(Z) ( P \otimes I_\cU) b(W)^*
\label{special}
\end{align}
 is a key piece in one of the main results of our paper \cite{BMV2}, specifically, the implication (1$^\prime$) $\Rightarrow$ (2)
in Theorem 3.1 there.  The cone separation argument as presented here in Section \ref{S:ker-dom}
largely follows the arguments in \cite{BMV2} but with some improvements in the exposition.
There are also some distinctive differences which we would like to point out:

\smallskip

\noindent
\textbf{(1)}  For the case where the kernel $\fQ$ has the special form in \eqref{special},  Lemma 4.1 in \cite{BMV2}
shows that {\em any cp nc kernel $\Gamma_2(Z,W)(P)$ can be represented as
$\Gamma_1(Z,W)(\fQ(Z,W)(P))$ for a cp nc kernel $\Gamma_1$.  Hence the second term 
$\Gamma_2(Z,W)(P)$ in the sought-after cone representation 
\eqref{cone-rep} can be absorbed into the first term and thus can be dropped in the formula \eqref{cone-rep}.}
The idea for this lemma in fact goes back to the seminal paper of Agler \cite{Agler-Hellinger} which handles the
commutative polydisk situation: $\fQ(z,w) =  \sbm{ 1 - z_1 \overline{w}_1 & & \\ & \ddots & \\ & &1-  z_d \overline{w}_d }$.

\smallskip

\noindent
\textbf{(2)}  The fact that the cone $\cC$ \eqref{cone-form} is not only closed and bounded but also pointed 
(Lemma \ref{L:cone}) was missed in \cite{BMV2}.  The fact that the cone is also pointed leads to a sharper cone-separation principle:
{\em if $\fS$ is not in the cone $\cC$, then there is a linear functional $\ell$ so that
$\operatorname{Re} \ell(\fS) < 0$ and $\operatorname{Re} \ell(\fS)(\cC \setminus \{0\}) >  0$  (not just 
$\operatorname{Re} \ell(\cC) \ge 0$.)}
Consequently the cone-separation argument in \cite{BMV2}
is more involved.  The inner product on the spaces $\cH$ in \cite{BMV2} is seen to only be a semi-inner product;
to form a Hilbert space one has to form equivalence classes.  Then it is more complicated to define the operators
$M^*_{\chi,0}$  and $M^*_{Q,0}$ on $\cH$ and $\cH_\cN$;  to make sense of such formal operators not necessarily
well defined on equivalence classes, it is necessary to introduce extra regularization terms depending on a parameter
$\rho < 1$ to form an approximate cone $\cC_\rho$, and then eventually let $\rho$ tend up to $1$.   

In this more complicated setting,  it should be explained that the definition of $\widehat M^*_{\chi,0}$ involves 
restriction to an invariant subspace followed by a quotient map, and in the end one only knows that
$M_{\chi,0}$ is similar to a point in the bi-full envelope $[\Omega]_{\rm bi-full}$ of $\Omega$ (points corresponding
to a semi-invariant subspace $\cN$ of a point $\widehat Z^{(0)}$ in $[\Omega]_{\rm nc}$).  In the present exposition
in  Section \ref{S:ker-dom}, there is no quotient operation involved in the definition of $\widehat M_{\chi,0}$
and it turns out that $M_{\chi,0}$ is in $[\Omega]_{\rm full}$ after all.

\smallskip

\noindent
\textbf{(3)}  The extension to the general case of the result for the special case where $\Omega$ consists of only finitely many points and
$\dim \cE < \infty$ as sketched here is exactly the same as the extension from the special case to the general case as
carried out in \cite{BMV2} for the case where $\fQ$ has the special form \eqref{special}.
\end{remark}

\begin{remark}  \label{R:PSS}
Theorem \ref{T:ker-dom} is somewhat unsatisfying since the kernel decomposition \eqref{cone-rep}
holds only for $Z,W \in \Omega$ even if we assume that the kernel $\fS$ is given on all of $\Xi$, unlike the case for Positivstellens\"atze
as appearing in classical commutative algebraic geometry (see e.g.~\cite{BCR, Marshall} as well as \cite{Drexel3Vin2016} for the  matrix-valued case).
To get a version of Theorem \ref{T:ker-dom} closer to the form of a classical Positivstellensatz, let us consider 
the following variation of Problem C:

\smallskip

\noindent
\textbf{Problem C$^\prime$:}  {\sl Given nc Hermitian kernels $\fQ$ and $\fS$ 
$$
\fQ \colon \Xi \times \Xi \to \cL({\mathbb C}_{\rm nc}, \cL(\cN)_{\rm nc}), \quad
\fS \colon \Xi \times \Xi \to \cL({\mathbb C}_{\rm nc}, \cL(\cE)_{\rm nc})
$$
(for auxiliary coefficient Hilbert spaces $\cN$ and $\cE$) such that
$$
Z \in \Xi_n \text{ and } \fQ(Z,Z)(I_n) \succ 0 \Rightarrow \fS(Z,Z)(I_{\cE^n} \succeq 0,
$$
find cp nc kernels
$$
\Gamma_1 \colon \Xi \times \Xi \to \cL(\cL(\cN)_{\rm nc},\cL(\cE)_{\rm nc}), \quad
\Gamma_2 \colon \Xi \times \Xi \to \cL({\mathbb C}_{\rm nc}, \cL(\cE)_{\rm nc})
$$
so that}
\begin{equation}   \label{PSS-ker-dec}
\fS(Z,W)(P) = \Gamma_1(\fQ(Z,W)(P)) + \Gamma_2(Z,W)(P)?
\end{equation}
The result of Theorem \ref{T:ker-dom} is that Problem C$^\prime$ has an affirmative solution under some additional hypotheses, but in a weaker form.  The weaker form is that 
we get cp kernels $\Gamma_1$, $\Gamma_2$ defined only for $Z,W \in {\mathbb P}_\fQ$ (the set where 
$\fQ(Z,Z)(I) \succ 0$), and then the decomposition \eqref{cone-rep} or \eqref{PSS-ker-dec} is guaranteed to hold only for $Z,W
\in {\mathbb P}_\fQ$.  To apply Theorem \ref{T:ker-dom} to arrive at this result, we take $\Omega = {\mathbb P}_\fQ$ in which case also 
$\Omega' = {\mathbb P}_\fQ$.  Then hypotheses (i) and (ii) in Theorem \ref{T:ker-dom} are satisfied, but the decomposability assumptions
(iii) and (iv) do not appear to hold in general (see Theorem \ref{T:Decom1} and the ensuing discussion).  On the other hand,
the hereditary case of Theorem 1.2 in the paper of Helton-McCullough \cite{HMcC2004} gives an affirmative answer to 
Conjecture C$^\prime$ under the assumptions
that $\cV = {\mathbb C}^d$,  $\Xi = \cV_{\rm nc}$, the coefficient Hilbert spaces $\cN$ and $\cE$ are finite-dimensional, and
that $\fQ$ and $\fS$ are matrix-polynomial nc kernels, i.e., 
$$
\fQ(Z,W)(P) = q_1(Z) P q_2(W)^*, \quad \fS(Z,W)(P) = s_1(Z) P s_2(W)^*
$$
for nc matrix polynomials $q_1$, $q_2$, $s_1$, $s_2$ of appropriate sizes, and also the positivity set ${\mathbb P}_\fQ$
should be bounded in norm (the {\em Archimedean hypothesis}); here we make use of the observations of Section 3.6 of
\cite{BMV2} that the general kernel decomposition  \eqref{PSS-ker-dec} is equivalent to the special case
with $P=I$ and $Z = W$ as long as \eqref{PSS-ker-dec} holds over a set $\Xi$ which is finitely open.
This Positivstellensatz has now been extended to the case where
$\fQ$ and $\fS$ are {\em rational matrix Hermitian kernels} (i.e., the functions $q_1$, $q_2$, $s_1$, $s_2$ as above are 
nc rational matrix functions with intersection of domains equal to some generic nc subset $\cD$ of  ${\mathbb C}^d_{\rm nc}$).
Furthermore, in case the $\fQ$-positivity set ${\mathbb P}_\fQ$ is convex, results of Helton-Klep-McCullough \cite{HKMcC2012}
provide an explicit finite algorithm for actually computing the Kolmogorov decompositions for the kernels
$\Gamma_1$ and $\Gamma_2$ in the kernel decomposition \eqref{PSS-ker-dec}.
\end{remark}

\end{document}